\documentclass[10pt]{article} 
\usepackage{vmargin} 
\setmarginsrb{3.2cm}{2.5cm}{3.2cm}{2.5cm}{1cm}{1cm}{1cm}{1.6cm} 
\usepackage{amsmath} 
\usepackage{amssymb}

\title{Tensor Products of Principal Unitary Representations of Quantum Lorentz Group and Askey-Wilson Polynomials}

\author{E. Buffenoir\thanks{e-mail:buffenoi@lpm.univ-montp2.fr},\cr   
Laboratoire de Physique Math\'ematique et Th\'eorique\thanks{Laboratoire du CNRS UMR 5825},  
 \cr Universit\'e Montpellier 2,  Place Eug\`ene Bataillon \cr 34000 Montpellier France \cr  
 \cr  
Ph. Roche\thanks{e-mail:roche@math.mit.edu},\cr Dpt of  
Mathematics, MIT, \cr Cambridge, MA 02139 USA  
\cr {\em on leave from} \cr  
 CPT Ecole Polytechnique\thanks{Laboratoire   
 Propre du CNRS UPR 14}\cr  
91128 Palaiseau Cedex France}   
  
\date{\today}  
\begin{document}

\maketitle  
\begin{abstract}  
We study the tensor product of principal unitary representations of the quantum Lorentz Group, prove a decomposition theorem  
 and compute the associated intertwiners. We show that these  
intertwiners can be expressed in terms of complex continuations of $6j$ symbols of ${\mathfrak U}_q(su(2)).$ These intertwiners are expressed in terms of q-Racah polynomials and Askey-Wilson polynomials. The orthogonality of these intertwiners imply some relation mixing these two families of polynomials. The simplest of these relations is the orthogonality of Askey-Wilson Polynomials. 
\\Work supported by CNRS.  
\end{abstract}  
  
\def\kE#1#2#3#4#5#6{{\buildrel {#1} \over k}{}^{#2}_{\! #3} \!\otimes \!{\buildrel {#4} \over E}{}^{#5}_{#6}}  
  
\def\coefalambda#1#2#3{\Lambda^{#1#2}_{#3}}  
 \def\coefblambda#1#2#3#4{\Lambda^{#1#2}_{#3#4}}  
\def\coefclambda#1#2#3#4#5#6{\Lambda^{#1#2}_{#3#4}(#5,#6)}

\def\ontop#1#2{\buildrel {{}_{#1}}\over #2}  
\def\twoontop#1#2#3{\buildrel {{}_{(#1,#2)}}\over #3}  
  
\def\End{\rm End}  
\def\Hom{\rm Hom}  
\def\Irr{\rm Irr}  
\def\dim{\rm dim}  
\def\id{\rm id}  
\def\Mat{\rm Mat}  
  
\def\Rff#1#2#3#4{{\buildrel #1 #2 \over {\bf R}}_{{} \atop {\!\!\!#3#4}}}  
\def\Rppff#1#2#3#4#5#6{{\buildrel #1 #2 \over R}\!{}^{(+)}{}^{#3#4}_{#5#6}}  
\def\Rmmff#1#2#3#4#5#6{{\buildrel #1 #2 \over R}\!{}^{(-)}{}^{#3#4}_{#5#6}}  
\def\Rpff#1#2{\buildrel #1#2 \over {R'}}  
\def\Rmff#1#2{\buildrel #1#2 \over {R^{-1}}}  
  
\def\GPhin#1#2#3{\Phi^{{}_{(#1)}{}_{(#2)}}_{\;\;\;\;\;{}_{(#3)}}}  
\def\GPhic#1#2#3{{\hat{\Phi}}^{{}_{(#1)}{}_{(#2)}}_{\;\;\;\;\;{}_{(#3)}}}  
\def\GPhip#1#2#3{\Phi^{{}_{(#1)}{}_{(#2)}}_{{}_{[+]}\;\;{}_{(#3)}}}  
\def\GPhim#1#2#3{\Phi^{{}_{(#1)}{}_{(#2)}}_{{}_{[-]}\;\;{}_{(#3)}}}  
\def\GPhie#1#2#3{\Phi^{{}_{(#1)}{}_{(#2)}}_{{}_{[\epsilon]}\;\;{}_{(#3)}}}  
\def\GPhis#1#2#3{\Phi^{{}_{(#1)}{}_{(#2)}}_{{}_{[\sigma]}\;\;{}_{(#3)}}}  
  
\def\CGphi#1#2#3{\phi^{{#1}{#2}}_{#3}}  
\def\CGpsi#1#2#3{\psi_{{#1}{#2}}^{#3}}  
\def\Clebphi#1#2#3#4#5#6{\left(\begin{array}{ll}   
\!#4 & \!#5 \\ \!\!\!#1 & \!\!\!#2 \end{array}\vline \begin{array}{ll} \!#3 \\ #6 \end{array}  
 \!\!\right)}  
\def\Clebpsi#1#2#3#4#5#6{\left(\!\begin{array}{ll}   
#6 \\ \!#3 \end{array}\vline \begin{array}{ll} #1 & \!\!\!#2 \\ #4 & \!\!\!#5 \end{array}\!\!  
 \right)}  
\def\sixj#1#2#3#4#5#6{\left\{\!\!\begin{array}{ll}   
#4 & \!\!\!#5 \\ #1 & \!\!\!#2 \end{array}\!\vline \begin{array}{ll} #3 \\ #6 \end{array}  
 \!\!\right\}}  
\def\sixjn#1#2#3#4#5#6#7{\left\{\!\!\begin{array}{ll}   
#4 & \!\!\!#5 \\ #1 & \!\!\!#2 \end{array}\!\vline \begin{array}{ll} #3 \\ #6 \end{array}  
 \!\!\right\}_{\!\!{}_{(#7)}}}  
\def\valsixj#1#2#3#4#5#6{\left\{\!\!\left\{\!\!\begin{array}{ll}   
#4 & \!\!\!#5 \\ #1 & \!\!\!#2 \end{array}\!\vline \begin{array}{ll} #3 \\ #6 \end{array}  
 \!\!\right\}\!\!\right\}}  
\def\nor6j#1#2#3#4#5#6{{\cal N}\left(\!\!\begin{array}{ll}   
#4 & \!\!\!#5 \\ #1 & \!\!\!#2 \end{array}\!\vline \begin{array}{ll} #3 \\ #6 \end{array}  
 \!\!\right)}  
\def\noraskey6j#1#2#3#4#5#6{\Gamma\left(\!\!\begin{array}{ll}   
#4 & \!\!\!#5 \\ #1 & \!\!\!#2 \end{array}\!\vline \begin{array}{ll} #3 \\ #6 \end{array}  
 \!\!\right)}  
  
\def\ElemRed#1#2#3#4#5#6{\left[\!\!\begin{array}{ll} #4 & \!\!#5 \\   
\;\;\mbox{\small \it #1} & \!\mbox{\small \it #2} \end{array}  
\!\vline \!\begin{array}{ll} \;\;\mbox{\small \it #3} \\  #6 \end{array}  
 \!\!\!\right]}  
\def\Entrela#1#2#3#4#5#6#7#8#9{\left[\!\!\begin{array}{lll}  
 \; \mbox{\small \it #1 } & \; \mbox{\small \it #2 } \\   
#7 & #8\\  
\;\mbox{\small \it #4} & \; \mbox{\small \it #5 } \end{array}\!\vline \!\begin{array}{lll}   
\; \mbox{\small \it #3 }  \\   
#9\\  
\;\mbox{\small \it #6}   
 \end{array} \!\!\!\right]}

\def\qphit43#1#2#3#4#5#6#7{{}_4{\tilde\Phi}_{3}{\left[\!\!\begin{array}{cccc}   
 #1 &  #2 &  #3& #4\\   
 #5 &  #6 &  #7& {} \end{array}  
 \!\!\right]}}  
  
\def\qphi43#1#2#3#4#5#6#7#8#9{{}_4{\Phi}_{3}{\left[\!\!\begin{array}{cccc}   
 #1 &  #2 &  #3& #4\\   
 #5 &  #6 &  #7& {} \end{array};#8,#9  
 \right]}}

\def\qn#1{[#1]}  
\def\qd#1{[d_{#1}]}

\def\halfinteger{\frac{1}{2}\mathbb{Z}^{+}}  
\def\onehalf{\frac{1}{2}}  
\def\ZZ{\mathbb Z}  
\def\NN{\mathbb N}  
\def\RR{\mathbb R}  
\def\CC{\mathbb C}  
\def\SS{\mathbb S}  
\def\ZZ{\mathbb Z}  
\def\SS+{{\mathbb S}^{+}}  
\def\SSc{{\mathbb S}_c}  
\def\XX{\mathbb X}  
\def\Proof{\underline{\sf Proof:}}  
\def\eoProof{$\square$}  
  
\newtheorem{proposition}{Proposition}  
\newtheorem{lemma}{Lemma}  
\newtheorem{theorem}{Theorem}  
\newtheorem{definition}{Definition}  
  
\newpage  
  
\section{Introduction}  
In \cite{BR1} we have pursued the work of Podles-Woronowicz \cite{PW} and Pusz \cite{Pu}: we have  
shown that  the unitary representations of ${\mathfrak U}_q(sl(2,\mathbb{C})_{\RR})$ can be nicely
 expressed in  
terms of one variable complex continuation of  
$6j$ symbols of ${\mathfrak U}_q(su(2))$. Using this result we were able to construct the characters of these  
unitary representations and prove a Plancherel theorem for $L^2$ functions.  
One  has to work in the category of  $C^{*}$  
multiplier Hopf algebras \cite{PW,Wo2} in order to handle functional analysis problems.
 \\ 
Let $G$ be the complex Lie group $SL(2,\mathbb{C})$ and denote by   
$\{{\buildrel {\lambda} \over {\Pi}}, \lambda\in  {\frac{1}{2}}\mathbb{Z}\times \RR\}$   
the set of principal unitary representations of $G$. The Plancherel measure is given by   
$P(\lambda)d\lambda = {\frac{1}{2}}(m^2+\rho^2)d\rho$ where $\lambda=(m, \rho)$  
and ${\buildrel {\lambda} \over {\Pi}}$ is equivalent to ${\buildrel {-\lambda} \over {\Pi}}$.  
Naimark \cite{N} has shown the following decomposition theorem:  
 
$${\buildrel {\!{}_{(m,\rho)}} \over {\Pi}}\otimes {\buildrel 
{\!{}_{(m',\rho')}} \over {\Pi}}=  \bigoplus_{m''\in J_{m,m'} }  \int^{\oplus} 
\!\!\!d\rho'' \;\; {\buildrel {{}_{(m'',\rho'')}} \over {\Pi}}$$ 
 
 where 
$J_{m,m'}=\{m\in \onehalf\ZZ, m+m'+m''\in \ZZ\}.$    The aim of the present 
article is to prove the quantum analog of this theorem and to give   explicit 
formulae for the Clebsch-Gordan coefficients associated to this decomposition. 
 
 Let us also denote by $\{ {\buildrel {\lambda} \over {\Pi}},   \lambda\in 
\frac{1}{2}\mathbb{Z}\times ]-\frac{\pi}{\hbar},\frac{\pi}{\hbar}] \}$, where  
 $q=e^{-\hbar},$ the set of  principal representations of  ${\mathfrak 
U}_q(sl(2,\mathbb{C})_{\RR}).$ The representation ${\buildrel {\lambda} \over 
{\Pi}}$ is a unitary representation of  ${\mathfrak 
U}_q(sl(2,\mathbb{C})_{\RR})$ with domain $V_{\lambda}.$  We have shown in 
\cite{BR1}, that these representations can be constructed using complex 
continuation   in one variable (namely $\rho$) of $6j$ of ${\mathfrak 
U}_q(su(2)).$ Let us denote by $6j(1)$ the complex continuation   of these 
coefficients to distinguish them from the  $6j$ of finite dimensional 
representations of ${\mathfrak U}_q(su(2)),$ which  will be denoted by 
$6j(0).$   We have shown that the Plancherel measure is   $P(\lambda)d\lambda 
= (q-q^{-1})^2\frac{\hbar}{4 \pi}[m+i\rho][m-i\rho]d\rho$ where $\lambda=(m, 
\rho).$   

 In section 2 we recall definitions and properties of $6j(0)$ and 
$6j(1).$ We then introduce and study   the basic properties of complex 
continuations, in three independent continuous spins, of $6j$  symbols   that 
we call $6j(3).$ We review the main theorems of harmonic analysis on 
$SL_q(2,\CC)_{\RR}$ which are   needed in the sequel.    

 In section 3 we 
study the space of intertwiners   $\Phi^{\lambda\lambda'}_{\lambda''} :  
V_{\lambda}\otimes V_{\lambda'}\rightarrow  V_{\lambda''},$ and   give an 
expression for them in terms of $6j(1)$ and $6j(3).$    Even in the classical 
case, such a simple expression of $\Phi^{\lambda\lambda'}_{\lambda''}$ was not 
known.  We can define a linear map, denoted ${\hat \Phi},$ from 
$V_{\lambda}\otimes V_{\lambda'}$  to $\int^{\oplus} \!d\lambda''{\buildrel 
{\;\lambda''} \over {\Pi}},$ by associating to each   $u\in  
V_{\lambda}\otimes V_{\lambda'}$ the map   $\lambda''\mapsto\Phi^{\lambda 
\lambda'}_{\lambda''}(u).$  It remains to show that we can find a 
normalization $N(\lambda,\lambda',\lambda'')$   of $\Phi^{\lambda 
\lambda'}_{\lambda''}$ in such a way   that ${\hat \Phi}$ is an isometry.  
This cannot easily be obtained from the definition of the 
$\Phi^{\lambda\lambda'}_{\lambda''}$   in terms of $6j(3).$ 

 We therefore use 
another construction of the intertwiners using the quantum analogue of the 
operator   $\int dh(g) \; {\buildrel {\lambda} \over {\Pi}}(g)\otimes 
{\buildrel {\lambda'} \over {\Pi}}(g)  \otimes {\buildrel {\;\lambda''} \over 
{\Pi}}(g^{-1})$ where $dh$ is the Haar measure.  By using the relations 
between these operators and the $\Phi^{\lambda \lambda'}_{\lambda''}$ we are 
able   to compute the normalization factor.   This is the content of section 
4.    

 In section 5 we show that with this choice of normalization ${\hat 
\Phi}$ is an isometry.  We then show that $\Phi^{\lambda 
\lambda'}_{\lambda''}$ is expressed in terms of $q-$ Racah polynomials and  
Askey-Wilson polynomials.  As a result this last property implies non trivial 
identities  which mix q-Racah polynomials and Askey-Wilson polynomials.  

   It 
is important to keep in mind the following hierarchy of complex continuations 
of $6j$ symbols of ${\mathfrak U}_q(su(2))$:  \vskip -0.6cm \begin{itemize}  
\item $6j(0)$ are defined as being Racah coefficients (i.e recoupling 
coefficients) of finite dimensional   representations of ${\mathfrak 
U}_q(su(2)).$   \end{itemize} \begin{itemize} \item $6j(1)$ are building 
blocks of matrix elements of unitary representations of   ${\mathfrak 
U}_q(sl(2,\mathbb{C})_{\RR}).$ They are matrix elements of the universal 
shifted  cocycle \cite{BBB} and are equivalent to   the Fusion Matrix of 
\cite{EV}.   \item $6j(3),$ as we will show, are building blocks of the 
Clebsch-Gordan coefficients associated to   the tensor product of principal 
representations of  ${\mathfrak U}_q(sl(2,\mathbb{C})_{\RR}).$ For the   
moment there is no real understanding of these $6j(3)$ in terms   of Fusion 
matrix or as matrix elements   of some universal element.  \end{itemize}  
There also exists a final level of this hierarchy, called $6j(6)$, where the 6 
spins are arbitrary complex numbers. They are the building blocks of the Racah 
coefficients of principal representations of ${\mathfrak 
U}_q(sl(2,\mathbb{C})_{\RR}).$ Their expressions in terms of basic 
hypergeometric functions, as well as their properties will be given in 
\cite{BR2}.      

  \bigskip      \section{Definition and properties of various 
continuations of 6j symbols of ${\mathfrak U}_q(su(2))$}  

  Let us first 
recall some notations and results of \cite{BR1} which will be used throughout  
this work. The reader is also invited to read the first subsection of the 
appendix of the present article for definitions of basic hypergeometric  
functions.   

 \subsection{Intertwiners and 6j(0) of ${\mathfrak U}_q(su(2))$ } 
   Let  $q=e^{-\hbar}$ with $\hbar\in \RR^{+*},$  ${\mathfrak U}_{q}(su(2))$ 
is the star Hopf algebra generated by $q^{\pm H}, J^{(\pm)}$ with the defining 
relations:  \begin{eqnarray}  && q^{H}J^{(\pm)}q^{-H}=q^{\pm 
1}J^{(\pm)},\;\;\;\;\;\;\;\;[J^{(+)},J^{(-)}]=\frac{q^{2H}-q^{-2H}}{q-q^{-1}},\;\;\;\;\;\;\;\;\;\;\;\;\;\;\;\;\;\;\;\;(J^{(\pm)})^*=q^{\mp 1} J^{(\mp)},\nonumber\\  &&\Delta(J^{(\pm)})=q^{-H}\otimes J^{(\pm)}+J^{(\pm)}\otimes q^{H}, \;\;\;\;\;\;\;\;\Delta(q^{\pm H})=q^{\pm H} \otimes q^{\pm H},\;\;\;\;\;\;\;\;\;\;\;\;  (q^{H})^*=q^{H}.  
\end{eqnarray}  
This Hopf algebra is a ribbon quasi-triangular Hopf algebra, with a universal $R-$matrix denoted $R.$ We will define as usual  $R^{(+)}=R, R^{(-)}=R_{21}^{-1}$ and $\mu=q^{2H}.$ 
  
We will denote in the rest of this article $Irr({\mathfrak U}_{q}(su(2)))$ the set of all  
equivalence  classes of finite dimensional irreducible unitary representations with 
$Sp(q^{H}) \in \RR^+.$ They are completely classified by a half-integer  $K$  and we will denote
 by ${\buildrel K \over\pi}$ the representation of spin $K.$   The tensor product of elements of $Irr({\mathfrak U}_{q}(su(2)))$ is completely reducible in elements of $Irr({\mathfrak U}_{q}(su(2))).$ Let us define ${\buildrel K \over V}$ as being the vector space, of  dimension  $d_K=2K+1,$ associated to the representation of spin $K.$ Let   $({\buildrel {{}_{K}} \over e}_m)_{m=-K \cdots K},$ be an orthonormal basis of ${\buildrel K \over V}$ such that ${\buildrel {{}_{K}} \over e}_m$ is of weight $q^m$ for the action of $q^H$.  
The central ribbon element that we choose is such that $ {\buildrel K \over\pi}(v)=v_{K}id,$ where $v_K=e^{2i \pi K}q^{-2K(K+1)}.$ Note that we have included a sign in order to satisfy relation (\ref{twistclebsch}).

Let us introduce  the following notation: $\forall I,J,K \in \halfinteger, \forall m \in \onehalf\ZZ$ we define  
\vskip -0.5cm  
\begin{eqnarray}  
Y^{(0)}_{(I,J,K)}&=&  
\left\{\begin{array}{ll} 1\;\;\;\;\;\;  
\mbox{if\;\;\; $I\!+\!J\!-\!K, \;J\!+\!K\!-\!I,\; K\!+\!I\!-\!J \in \ZZ^+$} \\0 \;\;\;\;\;\;\mbox{otherwise}   
\end{array}\right.\\  
Y^{(1)}_{(I,m)}&=&  
\left\{\begin{array}{ll} 1\;\;\;\;\;\;  
\mbox{if\;\;\; $I+m,\; I-m \in \ZZ^+$} \\0 \;\;\;\;\;\; 
\mbox{otherwise.}\end{array}\right.  
\end{eqnarray}  
 Let us now recall some properties about braiding matrices, $3j$ and $6j$ symbols of the finite dimensional representations of ${\mathfrak U}_q(su(2)).$   
 We will always refer to conventions as well as explicit expressions given in  
 \cite{BR1}, see also \cite{KR,Ge,KV}.  
The  group-like element ${\buildrel { {}_{K}} \over \mu},$ the square root of the ribbon element $v^{1/2}_{ {}_{K}}$, the quantum Weyl element ${\buildrel  { {}_{K}}\over w},$ the braiding matrices ${\buildrel { {}_{I}{}_{J}} \over R}{}^{(\pm)}$, and the  Clebsch-Gordan coefficients satisfy the following relations:  
\vskip -0.4cm  
\begin{eqnarray}   
&&\hskip -0.4cm\sum_{m,n}\!\Clebpsi{I}{\!J}{L}{m}{n}{r}\!\!\Clebphi{I}{\!\!J}{K}{m}{\!n}{p}\!=\!Y^{(0)}_{(I,J,K)}\delta_{K,L}\delta_{p}^{r}Y^{(1)}_{(K,p)},\nonumber\\  
&&\hskip -0.4cm \sum_{K,p}\!\Clebphi{I}{\!\!J}{K}{m}{\!\!n}{p}\!\!\Clebpsi{I}{J}{K}{i}{j}{p}\!\!=\delta_{i}^{m}\delta_{j}^{n}Y^{(1)}_{(I,m)}Y^{(1)}_{(J,n)},\label{ortho1clebsch}
\end{eqnarray}  
\vskip -0.4cm   
\begin{eqnarray}   
&&\hskip -1cm\sum_{k,l}\Clebpsi{I}{J}{K}{k}{l}{p} {\buildrel { {}_{J}{}_{I}} \over R}{}^{(\pm)}{}_{ji}^{lk}=(\frac{v^{1/2}_{I}v^{1/2}_{J}}{v^{1/2}_{K}})^{\pm 1} \Clebpsi{J}{I}{K}{j}{i}{p},\label{twistclebsch}  
\end{eqnarray}  
\vskip -0.4cm  
\begin{eqnarray}   
&&\hskip -2cm\sum_{i'j'm}\Clebpsi{A}{B}{C}{i'}{j'}{k} {\buildrel { {}_{A}{}_{D}} \over R}{}^{(\pm)}{}_{im}^{i'l}   
{\buildrel { {}_{B}{}_{D}} \over R}{}^{(\pm)}{}_{jn}^{j'm}= \sum_p{\buildrel { {}_{C}{}_{D}} \over R}{}^{(\pm)}{}_{pn}^{kl}  \Clebpsi{A}{B}{C}{i}{j}{p},\label{quasitriang}  
\end{eqnarray}  
\vskip -0.4cm  
\begin{eqnarray}  
&&\hskip -1cm \sum_{a'}v^{1/2}_{I} \; {\buildrel { {}_{I}} \over w}{}_{aa'} \Clebphi{I}{J}{K}{a'}{b}{k}=e^{i \pi (J-K)} (\frac{\qd{K}}{\qd{J}})^{1/2}\Clebpsi{I}{K}{J}{a}{k}{b}= \sum_{k'}v^{1/2}_{I} \;{\buildrel { {}_{K}} \over w}{}_{k'k} \Clebphi{J}{K}{I}{b}{k'}{a},\nonumber\\  
&&\hskip -1cm\sum_{i',k'} \ontop{J}{w}_{ii'} {\buildrel { {}_{J}{}_{I}} \over R}{}^{(\pm)}{}_{ji'}^{lk'} \ontop{J}{w}{}^{k'k}={\buildrel { {}_{I}{}_{J}} \over R}{}^{(\mp)}{}_{kj}^{il} \;,\;\;\;\;\;\;\;\;\;\;\;\;\;\;\;\;\;\;{v_{J}}\sum_k \ontop{J}{w}{}^{ik} \ontop{J}{w}_{kj}=\delta^i_j Y^{(1)}_{(J,j)}, \label{contra} 
\end{eqnarray}  
\vskip -0.5cm  
\begin{eqnarray}  
&&\hskip -2cm \Clebpsi{I}{J}{0}{a}{b}{0}\!\!=\!\frac{\delta_{I,J} \ontop{J}{w}_{ab} e^{i \pi J} }{\sqrt{\qd{I}}},  
\;\;\;\;\;\;\;\;\;\;\;\;{v_{J}}\sum_b \ontop{J}{w}^{ab} \ontop{J}{w}_{kb}=e^{2i \pi J}\ontop{J}{\mu}{}^a_k, \;\;\;\;\;\;\;\;\;\;\qd{K}\!=\!\sum_a \ontop{K}{\mu}{}^a_a, \label{normaclebsch}  
\end{eqnarray}  
\vskip -0.5cm  
\begin{eqnarray}  
&&\hskip -2cm \Clebpsi{A}{B}{C}{i}{j}{k}=\Clebphi{A}{B}{C}{i}{j}{k}\in {\mathbb R}  
.\end{eqnarray}   
  
The Clebsch-Gordan coefficients also satisfy the relations \cite{KV}: 
\begin{eqnarray}  
 \Clebpsi{A}{B}{C}{i}{j}{k}(q)=\Clebpsi{B}{A}{C}{-j}{-i}{-k}(q)=(-1)^{A+B-C} \Clebpsi{A}{B}{C}{-i}{-j}{-k}(q^{-1})\label{symmetry3j}. 
\end{eqnarray} 
  
\bigskip 
  
Let us now recall basic facts about $6j(0)$ that we will use  extensively in our work. 
 
 6j(0) of ${\mathfrak U}_q(su(2))$ are defined as follows:  
\begin{eqnarray*}  
&&\hskip -1cm\sixjn{C}{F}{E}{A}{B}{D}{0}\delta_{F,H}\delta^{p}_{n}Y^{(1)}_{(H,p)}=\sum_{i,j,k,l,m}\Clebpsi{A}{B}{E}{i}{j}{m}\Clebpsi{E}{C}{H}{m}{k}{p}\Clebphi{B}{C}{D}{j}{k}{l}\Clebphi{A}{D}{F}{i}{l}{n}.\label{vertexIRF}  
\end{eqnarray*}   
  
The properties of the Clebsch-Gordan coefficients, recalled above, imply the following relations on 6j(0):  
\begin{eqnarray}  
&&\sixjn{C}{D}{E}{A}{B}{F}{0}=\sixjn{D}{C}{E}{B}{A}{F}{0}=\sixjn{A}{B}{E}{C}{D}{F}{0}=\sixjn{C}{B}{F}{A}{D}{E}{0}\mbox{(Symmetries)},\label{trivsymm6j0}  
\end{eqnarray}  
\vskip -0.4cm  
\begin{eqnarray}  
&&\sixjn{B}{E}{F}{A}{C}{D}{0}=e^{i \pi(C-F+E-D)}(\frac{[d_F][d_D]}{[d_C][d_E]})^{1/2}\sixjn{B}{D}{C}{A}{F}{E}{0}\mbox{(Racah-Wigner)},\label{racahwigner6j0}  
\end{eqnarray}  
\vskip -0.4cm  
\begin{eqnarray}  
&&\hskip 0.7cm\sum_{C}\sixjn{G}{H}{C}{A}{B}{I}{0}\sixjn{G}{H}{C}{A}{B}{J}{0}=\delta_{I,J}Y^{(0)}_{(A,H,I)}Y^{(0)}_{(B,G,I)}\mbox{(Orthogonality)},\label{ortho6j0}\\  
&&\hskip 0.4cm\sum_{C}\sixjn{H}{G}{C}{A}{B}{I}{0}\sixjn{G}{H}{C}{A}{B}{J}{0}(\frac{v_{J}^{1/2}v_{I}^{1/2}v_{C}^{1/2}}{v_{G}^{1/2}v_{H}^{1/2}v_{A}^{1/2}v_{B}^{1/2}})^{\pm 1}=\sixjn{B}{H}{I}{A}{G}{J}{0},\label{racah6j0}\\  
&&\sum_{A}\sixjn{I}{G}{A}{D}{F}{J}{0}\!\!\sixjn{E}{B}{A}{D}{F}{C}{0}\!\!\sixjn{G}{H}{B}{E}{A}{I}{0}\!\!=\sixjn{J}{H}{C}{E}{F}{I}{0}\!\!\sixjn{H}{G}{B}{D}{C}{J}{0}\label{penta6j0}.  
\end{eqnarray}  
The relation (\ref{racah6j0}) is called the ``Racah-relation'', whereas the relation  (\ref{penta6j0}) is usually referred to as the  ``Pentagonal equation'' or the ``Biedenharn-Elliot''equation.  
The  relations  (\ref{racah6j0})(\ref{penta6j0}) imply the Yang-Baxter equation on $6j(0),$ also called ``Hexagonal relation''.  
  
 \medskip 
  
Explicit expressions, as well as  special values for low spins, of $6j(0)$ symbols  are given in \cite{BR1}, but we want to add here the following asymptotic formulae which will be of importance in the rest of this paper (for details see \cite{KV}):  
\vskip -0.4cm  
\begin{eqnarray}  
&&\hskip -1.2cm \lim_{K \rightarrow +\infty}\sixjn{K}{K\!+\!m\!+\!n}{C}{A}{B}{K\!+\!n}{0}=e^{i \pi (A+B-C)}\Clebpsi{A}{B}{C}{m}{n}{m\!+\!n},\label{asymptot1}  
\end{eqnarray}  
\vskip -0.65cm  
\begin{eqnarray}  
&&\hskip -1.2cm \lim_{K \rightarrow +\infty}\sixjn{D}{K}{K\!+\!n_2}{A}{K\!+\!n_1\!+\!n_2}{K\!+\!n'_1}{0}\!  
(\frac{v^{1/2}_{K+n_2}v^{1/2}_{K+n'_1}}{v^{1/2}_{K}v^{1/2}_{K+n_1+n_2}})^{\pm 1} \delta_{n_1+n_2,n'_1+n'_2}=\!  
{\buildrel {{}_A {}_D} \over R}{}^{(\pm)}{}^{n'_1 n'_2}_{n_1 n_2}.\label{asymptot2}  
\end{eqnarray}  
  
\subsection{Definitions and  properties of continuation of  quantum $6j$ symbols}  
  
In order to understand the problem of continuation of quantum $6j$ symbols the reader is invited to read   
the appendix of \cite{BR1} where definitions, proofs, as well as references on $6j(1)$ are given.   
  
Let us give a definition which is essential in order to describe the complex continuations of quantum $6j$ symbols (the notations are explained in the appendix of the present work).  
  
\begin{definition}  
Let $\XX=\{(T,U,V,X,Y,Z)\in \CC^6, (2T, T+U-V, T+V-U, T+Y-Z, T+Z-Y)\in \NN^5\},$  we define for any $(T,U,V,X,Y,Z)\in \XX,$ 
\begin{eqnarray}  
\hskip -0.5cm\nor6j{X}{Y}{V}{T}{U}{Z} \!\!\!\!&=&\!\!\!\!  
\;\;e^{i\pi (T+Y-Z)} \;q^{(T+Y-Z)(U-T+Y-X+1)+(V+Y-X)(T+V-U)}\nu_{1}(d_Z)\nu_{1}(d_V)  
\times\nonumber\\  
&&\times\frac{(U\!-\!T\!+\!X\!-\!Y\!+\!1)_{\infty}\;\;\omega(Y;V,X)\omega(T;Y,Z)\omega(T;V,U)}{\!\!  
(T\!-\!U\!+\!X\!+\!Y\!+\!1)_{\infty}(2T+1)_{\infty}(1)_{\infty}\;\;\omega(U;X,Z)  
}\nonumber\\  
\hskip -0.5cm\valsixj{X}{Y}{V}{T}{U}{Z}\!\!\!\!&=&\!\!\!\!\nor6j{X}{Y}{V}{T}{U}{Z}\!  
{}_4\Phi_3\!  
\left[\!\!\begin{array}{cccc} \mbox{\it \small U$-\!$V$-\!$T} &\; \mbox{\it \small U$+\!$V$-\!$T$+\!$1} & \;\;\mbox{\it \small Z$-\!$Y$-\!$T} & \;\;\;\mbox{\it \small $-\!$Z$-\!$Y$-\!$T$-\!$1}\\  
& \hskip -2cm\mbox{\it \small $-\!2$T} & \hskip -1.8cm\mbox{\it \small $\!-\!$Y$-\!$X$-\!$T$+\!$U} & \hskip -1cm \mbox{\it \small U$+\!$X$-\!$Y$-\!$T$+\!$1} \end{array} \!\!\right].  
\end{eqnarray}  
\end{definition}  
\medskip  
It is important to understand that for $(T,U,V,X,Y,Z)\in \XX,$ $\nor6j{X}{Y}{V}{T}{U}{Z}$ and $\valsixj{X}{Y}{V}{T}{U}{Z}$ are both square roots of rational functions in the variables $q^{2T}, q^{2U},q^{2V}, q^{2X},q^{2Y}, q^{2Z}.$ This is a simple consequence of the fact that the hypergeometric series is of terminating type.  

 \medskip

We have been very cautious with the determination of the signs. This annoying problem already appeared in the case of $6j(0)$ but can be handled quite easily. This problem is strengthned  in the case of $6j(1)$ and $6j(3)$ because in that case we really have to take square roots of complex numbers.  This is the reason why we introduced   a particular square root, denoted $\nu_{\infty}(x)$  of the Eulerian product $(x)_{\infty}$ for $x\in \CC.$  
 
\medskip

It will  sometimes be useful to have  the explicit value of the $6j$ symbols when  
 one of the spins is equal to $0$ or $\onehalf$. From the last definition we easily get:   
\begin{eqnarray}  
&&\forall A,B,C \in \CC, \valsixj{A}{B}{C}{0}{C}{B}=1\label{sixjzerocont}  
\end{eqnarray}  
\vskip -0.6cm  
\begin{eqnarray}  
&&\hskip -1.5cm\valsixj{A}{B}{C}{\onehalf}{C+\onehalf}{B+\onehalf}=\frac{\nu_1(B+C-A+1)\nu_1(A+B+C+2)}{\nu_1(2C+2)\nu_1(2B+1)}\nonumber  
\end{eqnarray}  
\vskip -0.5cm   
\begin{eqnarray}  
&&\valsixj{A}{B}{C}{\onehalf}{C+\onehalf}{B-\onehalf}=\;-\;\frac{\nu_1(A+B-C)\nu_1(A+C-B+1)}{\nu_1(2C+2)\nu_1(2B+1)}\;q^{C+B-A+1}\nonumber  
\end{eqnarray}  
\vskip -0.5cm    
\begin{eqnarray}  
&&\hskip -0.75cm\valsixj{A}{B}{C}{\onehalf}{C-\onehalf}{B+\onehalf}=\frac{\nu_1(A+C-B)\nu_1(A+B-C+1)}{\nu_1(2C)\nu_1(2B+1)}\;q^{C+B-A}\nonumber  
\end{eqnarray}  
\vskip -0.5cm 
\begin{eqnarray}  
&&\hskip -2cm\valsixj{A}{B}{C}{\onehalf}{C-\onehalf}{B-\onehalf}=\frac{\nu_1(A+B+C+1)\nu_1(C+B-A)}{\nu_1(2C)\nu_1(2B+1)}.\label{sixjonehalfcont}  
\end{eqnarray}  
\medskip  
  
The usual  $6j(0)$ symbols which  properties were described in the last section and  denoted $\sixjn{D}{E}{C}{A}{B}{F}{0}$ where $A,B,C,D,E,F \in \halfinteger,$ are given by:  
  
\begin{eqnarray}  
\sixjn{D}{E}{C}{A}{B}{F}{0}=\valsixj{D}{E}{C}{A}{B}{F}  Y^{(0)}_{(A,B,C)} Y^{(0)}_{(A,E,F)}  Y^{(0)}_{(D,B,F)} Y^{(0)}_{(D,E,C).}  
\end{eqnarray}  
This expression can  easily be obtained from the usual expressions \cite{KR,Ge} using the inversion relation and the Sears identities  
(\ref{inversionrelation},\ref{sears})  recalled in the appendix.

 In the previous  formula, and in the following ones, we will make a distinction between the  first part of the right handside, which is called "explicit value", and the second part, products of $Y$ functions  and called  "selection rules". 
  
  In order to  describe  properties of continuation of 6j symbols, we will make in the sequel a convenient abuse of notation, which greatly simplifies  formulae: \\  
\indent if $X_1\in \CC$ is fixed, if $k$ is a positive integer  and $f:\CC^{k+1}\rightarrow \CC$ is a function, a  series  of the type $\sum_{X_2,X_3,...,X_k}f(X_1,X_2,...,X_k)$ will always be defined as  
\begin{equation} 
\sum_{X_2,X_3,...,X_k}f(X_1,X_2,...,X_k)=\sum_{n_1,n_2,...,n_k \in \onehalf\ZZ}f(X_1, X_1+n_1,...,X_k+n_k).   
\end{equation}  
 
 Note that  we  
will only use the notation $X_0, X_1$ to denote a couple of complex numbers  such that  
$X_0=-{\overline {X_1}}-1,$ and $X_0-X_1 \in\onehalf \ZZ,$   
(this notation will be explained in the next part.)  
  
Let us define the involutive endomorphism of the complex line: $\forall X\in \CC, X\mapsto \underline{X}=-X-1,$ as we will see it is important to understand the action of this symmetry on the $6j$ symbols.

In \cite{CGR,BR1}  two types of $6j(1)$ symbols were defined. The first type is a family of numbers denoted $\sixjn{X_1}{X_2}{C}{A}{B}{X_3}{1}$ and  the second one is denoted  $\sixjn{B}{X_1}{X_3}{A}{X_2}{X_4}{1},$ where both are defined for   $A,B,C \in \halfinteger, \forall i,j \in \{1,\cdots,4\}, X_i \in \CC-\halfinteger, X_i\!-\!X_j \in \onehalf\ZZ.$

 Their explicit expressions are given by  
\begin{eqnarray}  
&&\hskip -1.35cm\sixjn{X_1}{X_2}{C}{A}{B}{X_3}{1}\!\!\!\!\!=\!\valsixj{X_1}{X_2}{C}{A}{B}{X_3}  Y^{(1)}_{(A,X_2-X_3)} Y^{(1)}_{(B,X_1-X_3)}  Y^{(1)}_{(C,X_1-X_2)} Y^{(0)}_{(A,B,C)},\nonumber\\ 
&&\hskip -1.35cm\sixjn{B}{X_1}{X_3}{A}{X_2}{X_4}{1}\!\!\!\!\!=\!\valsixj{B}{X_1}{X_3}{A}{X_2}{X_4}  Y^{(1)}_{(A,X_2-X_3)} Y^{(1)}_{(A,X_1-X_4)}  Y^{(1)}_{(B,X_1-X_3)} Y^{(1)}_{(B,X_2-X_4)}.  
\end{eqnarray}  
These formulae deserve some remarks,    
concerning the very word of "continuation". The explicit value of the $6j(1)$ is a regular function when the $X_i$ approach half-integer values, but  the selection rules appears to be different, i.e the support of the continuated 6j(1) are larger than that of the 6j(0) symbols. Nethertheless, for fixed $I,J\in \halfinteger,$ and for a sufficiently large half-integer $K$ we have $Y^{(1)}_{(I,J)}= Y^{(0)}_{(I,K,J+K)}.$ It is in this sense that the term  ``continuation'' to complex spins has been used. 
  
 From  arguments developped in \cite{CGR,BR1}, we can check different polynomial identities which are the continuation of the  properties satisfied by the 6j(0). Rather than being exhaustive we  will just mention some of them which will be important in  our present work:  
  
\begin{proposition}  
\begin{eqnarray}  
&&\hskip -1cm\sixjn{X_1}{X_2}{C}{A}{B}{X_3}{1}\!\!\!\!=\!\sixjn{X_2}{X_1}{C}{B}{A}{X_3}{1}\!\!\!\!=\!(e^{i \pi}q)^{(X_1-X_3+A-C)}  
\frac{\nu_1(d_C)\nu_1(d_{X_3})}{\nu_1(d_A)\nu_1(d_{X_1})}\sixjn{X_3}{X_2}{A}{C}{B}{X_1}{1},\nonumber\\  
&&\hskip -1cm\sixjn{B}{X_1}{X_3}{A}{X_2}{X_4}{1}\!\!\!\!=\!\sixjn{B}{X_2}{X_4}{A}{X_1}{X_3}{1}=\sixjn{A}{X_2}{X_3}{B}{X_1}{X_4}{1}\!\!\!\!=\!\!\nonumber\\  
&&\hskip 3cm=\!(e^{i \pi}q)^{(X_1-X_3+X_2-X_4)}  
\frac{\nu_1(d_{X_4})\nu_1(d_{X_3})}{\nu_1(d_{X_2})\nu_1(d_{X_1})} \sixjn{B}{X_4}{X_2}{A}{X_3}{X_1}{1},\label{symcont}  
\end{eqnarray}    
\begin{eqnarray}  
&&\hskip -2cm\sum_{C}\sixjn{X_1}{X_2}{C}{A}{B}{X_3}{1}\!\sixjn{X_1}{X_2}{C}{A}{B}{X_4}{1}=\delta_{X_3,X_4}Y^{(1)}_{(A,X_2-X_3)}Y^{(1)}_{(B,X_1-X_3)},\label{ortho1cont}\\ 
&&\hskip -2cm\sum_{X_3}\!\sixjn{X_1}{X_2}{C}{A}{B}{X_3}{1}\!\sixjn{X_1}{X_2}{D}{A}{B}{X_3}{1}=\delta_{C,D}Y^{(0)}_{(A,B,C)}Y^{(1)}_{(C,X_1-X_2)},\label{ortho2cont}\\ 
&&\hskip -2cm\sum_{X_3}\!\sixjn{B}{X_1}{X_3}{A}{X_2}{X_4}{1}\!\sixjn{B}{X_1}{X_3}{A}{X_2}{X_5}{1}=\delta_{X_4,X_5}Y^{(1)}_{(A,X_1-X_4)}Y^{(1)}_{(B,X_2-X_4)},\label{ortho3cont} 
\end{eqnarray}    
\begin{eqnarray}  
&&\hskip -1.5cm\sum_{C}\!\sixjn{X_2}{X_1}{C}{A}{B}{X_3}{1}\!\!\sixjn{X_1}{X_2}{C}{A}{B}{X_4}{1}\!(\frac{v_{X_4}^{1/2}v_{X_3}^{1/2}v_{C}^{1/2}}{v_{X_1}^{1/2}v_{X_2}^{1/2}v_{A}^{1/2}v_{B}^{1/2}})^{\pm 1}\!=\!\sixjn{B}{X_2}{X_3}{A}{X_1}{X_4}{1},\label{racahcont}  
\end{eqnarray}    
\begin{eqnarray}  
&&\hskip -1.3cm\sum_{A}\!\sixjn{X_3}{X_1}{A}{D}{F}{X_4}{1}\!\!\!\sixjn{E}{B}{A}{D}{F}{C}{0}\!\!\!\sixjn{X_1}{X_2}{B}{E}{A}{X_3}{1}\!\!\!\!\!=\!\sixjn{X_4}{X_2}{C}{E}{F}{X_3}{1}\!\!\!\sixjn{X_2}{X_1}{B}{D}{C}{X_4}{1},\label{penta1cont}  
\end{eqnarray}    
\begin{eqnarray}  
&&\hskip -1cm \sum_{X_6}\!\sixjn{D}{X_6}{X_4}{C}{X_1}{X_2}{1}\!\!\!\sixjn{A}{X_5}{X_6}{C}{X_2}{X_3}{1}\!\!\!\sixjn{X_4}{X_5}{B}{A}{D}{X_6}{1}\!\!\!\!\!=\!\sixjn{B}{X_5}{X_4}{C}{X_1}{X_3}{1}\!\!\!\sixjn{X_1}{X_3}{B}{A}{D}{X_2}{1}\!\!,\label{penta2cont}  
\end{eqnarray}   
\begin{eqnarray}  
&&\hskip -1.3cm\sum_{X_5}\!\sixjn{X_2}{X_1}{P}{B}{A}{X_5}{1}\!\!\!\sixjn{X_4}{X_5}{M}{A}{Q}{X_2}{1}\!\!\!\sixjn{X_4}{X_1}{N}{B}{M}{X_5}{1}\!\!\!\!\!=  
\!\sixjn{Q}{N}{P}{B}{A}{M}{0}\!\!\!\sixjn{X_4}{X_1}{N}{P}{Q}{X_2}{1},\label{penta3cont}  
\end{eqnarray}    
\begin{eqnarray}  
&&\hskip -1.3cm\sum_{M}\!\sixjn{X_4}{X_3}{M}{A}{P}{X_2}{1}\!\sixjn{B}{N}{M}{A}{P}{C}{0}\!\sixjn{X_4}{X_1}{N}{B}{M}{X_3}{1}  
(\frac{v_P^{1/2}v_{X_3}^{1/2}}{v_M^{1/2}v_{X_2}^{1/2}})^{\pm 1}=\nonumber\\  
&&\hskip 1cm=\sum_{X_5}\!\sixjn{X_4}{X_5}{C}{B}{P}{X_2}{1}\!\sixjn{B}{X_1}{X_3}{A}{X_2}{X_5}{1}\!\sixjn{X_4}{X_1}{N}{A}{C}{X_5}{1}  
(\frac{v_C^{1/2}v_{X_1}^{1/2}}{v_N^{1/2}v_{X_5}^{1/2}})^{\pm 1}.\label{hexacont}  
\end{eqnarray}  
\end{proposition}  
  
It is important to stress  that every sums in the previous formulae are finite sums  because of the selection rules   entering into the  definitions of $6j(1).$

 We also have the following symmetry:  
\begin{proposition} \label{propsym6j1} 
\vskip-0.2cm  
\begin{eqnarray}  
&&\hskip -1cm  
\sixjn{\underline{X_1}}{\underline{X_2}}{C}{A}{B}{\underline{X_3}}{1}=  
\sixjn{X_1}{X_2}{C}{A}{B}{X_3}{1}  
(-1)^{A+B-C}  
\label{symbiz6j1}  
\end{eqnarray}  
\end{proposition}  
\vskip -0.2cm  
 
\Proof  
 
{}From the Sears identity and the inversion relation (\ref{sears}, \ref{inversionrelation}) we easily obtain the relation:  
  
$ 
\sixjn{\underline{X_1}}{\underline{X_2}}{C}{A}{B}{\underline{X_3}}{1}=  
\sixjn{X_1}{X_2}{C}{A}{B}{X_3}{1}  
f(A,B,C,X_1,X_1,X_3)  
$ 
where   
$$f(A,B,C,X_1,X_1,X_3)=(-1)^{2A}\frac{\varphi_{(A+X_2+X_3,2A+1)}  
\varphi_{(C+X_1+X_2,2C+1)}}{\varphi_{(2X_3+1,1)}\varphi_{(B+X_1+X_3,2B+1)}}.$$  
Using the explicit expression of $\phi(\alpha,n)$ for $n\in \ZZ$ explained  in the appendix, we conclude that $f=(-1)^{A+B-C}.$  
\eoProof  
  
\medskip  
  
In order to obtain neat expressions for the  Clebsch-Gordan coefficients of principal representations of the  quantum Lorentz group, it is necessary to introduce a new level in this  hierarchy of continuations of $6j$ symbols.  
  
\begin{definition}  
 We will define  $6j(3)$ symbols  to be the  family of numbers denoted $\sixjn{Z_1}{Y_1}{X_2}{A}{X_1}{Y_2}{3}$ where $A \in \halfinteger, X_1,X_2,Y_1,Y_2,Z_1 \in \CC-\halfinteger, X_1\!-\!X_2 \in \onehalf\ZZ,Y_1\!-\!Y_2 \in \onehalf\ZZ$ and defined by  
\begin{eqnarray}  
\sixjn{Z_1}{Y_1}{X_2}{A}{X_1}{Y_2}{3}=\valsixj{Z_1}{Y_1}{X_2}{A}{X_1}{Y_2}\;Y^{(1)}_{(A,X_1-X_2)}Y^{(1)}_{(A,Y_1-Y_2)}.  
\end{eqnarray}  
\end{definition}

These $6j(3)$ satisfy properties, which will be used in the rest of  
the paper and which are continuation of the identities satisfied by $6j(0)$  
and $6j(1).$  We  prefered to give  combinatorial proofs of these properties  
rather than continuation arguments, in order to prepare to the 
 $6j(6)$ case. Here again, the sums are  finite which is a consequence 
 of the selection rules  in    
the definition of  $6j(3)$. Note that these $6j(3)$  satisfy other 
 relations, which will be  explained in the section 5,  as consequences of the  study of  
the tensor product of principal  representations of the quantum Lorentz group.    
  
\begin{proposition}  
 The following symmetry properties are satisfied:  
\vskip -0.2cm  
\begin{eqnarray}  
&&\hskip -1.4cm\sixjn{Z_1}{Y_1}{X_2}{A}{X_1}{Y_2}{3}\!\!\!\!=\!  
\sixjn{Z_1}{X_1}{Y_2}{A}{Y_1}{X_2}{3}\!\!\!\!=\!(e^{i \pi}q)^{(X_1-X_2+Y_1-Y_2)}  
\frac{\nu_1(d_{X_2})\nu_1(d_{Y_2})}{\nu_1(d_{X_1})\nu_1(d_{Y_1})}  
\sixjn{Z_1}{Y_2}{X_1}{A}{X_2}{Y_1}{3}.\label{sym6j3}  
\end{eqnarray}  
\end{proposition}  
 
\Proof  
 
 Simple use of  Sears identity (\ref{sears}).  
\eoProof  
  
\begin{proposition} We also have a discrete orthogonality property 
\begin{eqnarray}  
\sum_{X_2}\sixjn{Z_1}{Y_1}{X_2}{A}{X_1}{Y_2}{3}\sixjn{Z_1}{Y_1}{X_2}{A}{X_1}{Y_3}{3}=\delta_{Y_2,Y_3}Y^{(1)}_{(A,Y_1-Y_2)}.\label{ortho16j3}  
\end{eqnarray}  
\end{proposition}  
\vskip -0.4cm 
\Proof  
 
This  is equivalent to the  orthogonality of q-Racah polynomials  (see our appendix for notations (\ref{defracahpol}) and \cite{W}\cite{AW1} for  proofs):  
\begin{eqnarray}  
\sum_{x=0}^{N} w^{(R)}(x;a,b,c,d) p^{(R)}_n(\mu(x);a,b,c,d) p^{(R)}_m(\mu(x);a,b,c,d)=\delta_{n,m}h^{(R)}_n(a,b,c,d),  
\end{eqnarray}  
\vskip -0.3cm 
\begin{eqnarray}  
&&\hskip -4cm\mbox{because }\;\;\;\;\;\nor6j{Z_1}{Y_1}{X_2}{T}{X_1}{Y_2}^2 =\frac{w^{(R)}(x;a,b,c,d)}{h^{(R)}_n(a,b,c,d)}, \mbox{and}\nonumber 
 \end{eqnarray}  
\vskip -0.3cm 
\begin{eqnarray}  
&&\hskip -1.5cm{}_4\Phi_3\!\left[\!\!\begin{array}{cccc} \mbox{\it \small $X_1\!-\!X_2-\!$T} &\; \mbox{\it \small $X_1\!+\!X_2-\!$T$\;+$1} & \;\;\;\;\;\mbox{\it \small $Y_2\!-\!Y_1-\!$T} & \;\;\;\mbox{\it \small $-Y_2\!-\!Y_1\!-\!$T$-\!$1}\\ & \hskip -3.5cm\mbox{\it \small $-\!2$T} & \hskip -3.5cm\mbox{\it \small $\!-\!Y_1\!-\!Z_1-\!$T$\;+\!X_1$} & \hskip -1.8cm \mbox{\it \small $X_1+Z_1\!-\!Y_1\!-\!$T$\;+$1} \end{array} \!\!\right]\!\!=p^{(R)}_n(\mu(x);a,b,c,d),\nonumber\\  
&&\mbox{with}\;\;\;\;\;n=T\!+\!X_2\!-\!X_1, \;x=T\!+\!Y_1\!-\!Y_2,\; N=2T,\nonumber\\  
&& a=q^{-4T-2},\; b=q^{4X_1+2}, \; c=q^{2(Z_1+X_1-Y_1-T)},\; d=q^{-2(Z_1+Y_1+X_1+T+2)}.  
\end{eqnarray}  
\vskip -0.3cm  
\eoProof  
  
These $6j(3)$ symbols satisfy the pentagonal relations:  
\begin{proposition}  
\begin{eqnarray}  
&&\hskip -1.2cm\sum_C\!\sixjn{Y_3}{Y_2}{C}{A}{B}{Y_1}{1}\!\!\sixjn{X_1}{X_2}{C}{B}{A}{X_3}{1}\!\!\sixjn{Z_1}{Y_3}{X_2}{C}{X_1}{Y_2}{3}\!\!\!\!=\sixjn{Z_1}{Y_1}{X_3}{A}{X_1}{Y_2}{3}\!\!\sixjn{Z_1}{Y_3}{X_2}{B}{X_3}{Y_1}{3},\label{penta1sixj3}\\  
&&\hskip -1.2cm\sum_{X_3}\!\sixjn{Z_1}{Y_3}{X_3}{B}{X_2}{Y_1}{3}\!\!\sixjn{B}{X_3}{X_2}{A}{X_1}{X_4}{1}\!\!\sixjn{Y_3}{Z_1}{X_3}{A}{X_4}{Z_2}{3}\!\!\!\!=\sixjn{Y_1}{Z_1}{X_2}{A}{X_1}{Z_2}{3}\!\!\sixjn{Z_2}{Y_3}{X_4}{B}{X_1}{Y_1}{3}.\label{penta2sixj3}  
\end{eqnarray}  
\end{proposition}  
 
\Proof  
 
These identities are simply proved by induction on $A.$ Indeed, for $A=0$ these identities are trivial and  for $A=\onehalf$ they are easily checked on the exact expression of the $6j(3)$ using the  finite difference equations verified by basic hypergeometric functions.  
 The induction  easily follows from  the use of orthogonality and pentagonal equations on $6j(1).$\\  
\eoProof  
  
Finally, we have a Racah relation:  
\begin{proposition}  
\begin{eqnarray}  
\sum_{Y_2}\sixjn{Z_2}{X_1}{Y_2}{A}{Y_1}{X_2}{3}\sixjn{X_1}{Y_1}{Z_1}{A}{Z_2}{Y_2}{3}(\frac{v^{1/2}_{Y_2}v^{1/2}_{Z_1}v^{1/2}_{X_2}}{v^{1/2}_{Z_2}v^{1/2}_{Y_1}v^{1/2}_{X_1}v^{1/2}_{A}})^{\pm 1}=\sixjn{Y_1}{Z_2}{X_2}{A}{X_1}{Z_1}{3}\label{racah6j3}  
\end{eqnarray}  
\end{proposition}  
 
\Proof  
 
It can  be proved by continuation arguments from the $6j(1)$ case.  
\eoProof  
\medskip  
  
We can prove other identities, which will be useful in the following chapters, and which are direct consequences of the previous ones.   
\begin{proposition}  
The following pentagonal equations are satisfied:  
\begin{eqnarray}  
\sum_{X_2}\sixjn{Y_1}{Z_3}{X_2}{C}{X_1}{Z_2}{3}\!\!\sixjn{X_1}{X_3}{B}{A}{C}{X_2}{1}\!\!\sixjn{Y_1}{Z_1}{X_3}{A}{X_2}{Z_3}{3}\!\!=\!\sixjn{Y_1}{Z_1}{X_3}{B}{X_1}{Z_2}{3}\!\!\sixjn{Z_1}{Z_2}{B}{C}{A}{Z_3}{1},\label{penta36j3}\\  
\sum_{Z_2}\sixjn{Y_1}{Z_1}{X_2}{A}{X_4}{Z_2}{3}\sixjn{Z_2}{Y_1}{X_4}{B}{X_1}{Y_2}{3}\sixjn{Y_2}{Z_1}{X_3}{A}{X_1}{Z_2}{3}=\sixjn{B}{X_2}{X_3}{A}{X_1}{X_4}{1}\sixjn{Z_1}{Y_1}{X_2}{B}{X_3}{Y_2}{3}.\label{penta46j3}  
\end{eqnarray}  
\end{proposition}  
 
\Proof  
 
It is obtained from the pentagonal relations  (\ref{penta1sixj3})(\ref{penta2sixj3}) and the discrete orthogonality relation (\ref{ortho16j3}).  
\eoProof  
  
We also have the hexagonal equation:  
\begin{proposition}  
\begin{eqnarray}  
&&\hskip -1cm\sum_{C}\frac{v^{1/2}_{C}v^{1/2}_{Y_2}}{v^{1/2}_{B}v^{1/2}_{Y_1}}\sixjn{Z_1}{X_1}{Y_1}{C}{Y_3}{X_2}{3}\sixjn{Y_3}{Y_1}{C}{A}{B}{Y_2}{1}\sixjn{X_1}{X_2}{C}{A}{B}{X_3}{1}=\nonumber\\  
&&\hskip 1cm=\sum_{Z_2}\frac{v^{1/2}_{X_2}v^{1/2}_{Z_2}}{v^{1/2}_{X_3}v^{1/2}_{Z_1}}\sixjn{Y_3}{Z_1}{X_2}{A}{X_3}{Z_2}{3}\sixjn{X_1}{Z_1}{Y_1}{A}{Y_2}{Z_2}{3}\sixjn{Z_2}{Y_3}{X_3}{B}{X_1}{Y_2}{3}.\label{hexa6j3}  
\end{eqnarray}  
\end{proposition}  
 
\Proof  
 
It is obtained by applying successively  (\ref{racah6j3})(\ref{penta1sixj3})(\ref{penta46j3})(\ref{racah6j3}).  
\eoProof  
 
 \medskip 
 
Finally the action of the automorphism of the complex line is given as follows:  
\begin{proposition} \label{propsym6j3} 
\begin{equation}  
\sixjn{\underline{Z_1}}{\underline{Y_1}}{\underline{X_2}}{A}  
{\underline{X_1}}{\underline{Y_2}}{3}=\sixjn{Z_1}{Y_1}{X_2}{A}{X_1}{Y_2}{3}  
g(A,X_1,X_2,Z_1,Y_1,Y_2)(-1)^{Y_1-Y_2+X_1-X_2}\label{symbiz6j3}  
\end{equation}  
where $g$ is a fourth root of unit which expression is:  
\begin{equation}  
g=  
\frac{\varphi_{(-Y_2+X_1+Z_1,X_1-X_2+Y_1-Y_2)}  
\varphi_{(X_1+Y_2+Z_1+1,X_1-X_2+Y_2-Y_1)}  
\varphi_{(X_1+Y_2-Z_1,X_1-X_2+Y_2-Y_1)}}  
{\varphi_{(2Y_2+1,1)}  
\varphi_{(2X_2+1,1)}  
\varphi_{(Y_2+Y_1-A,-2A-1)}\varphi_{(X_1+X_2-A,-2A-1)}  
\varphi_{(Y_2+Z_1-X_1,Y_2-Y_1+X_2-X_1)}}.  
\end{equation}  
\end{proposition}  
  
\Proof  
 
It is proved using the inversion relation and the Sears identity. Using the explicit expression for   
$\varphi(\alpha, n)$ with $n$ integer, we obtain that   $g$ is a fourth root of unit, but there is   
no simpler formula for $g.$  
\eoProof  
  
\subsection{Quantum Lorentz group, Principal Unitary  Representations and Harmonic Analysis}  
  
We will  
recall in this section  fondamental results  on  $\mathfrak{U}_q(sl(2,\mathbb{C})_{\RR})$ and on   
the  harmonic analysis \cite{BR1} on $SL_q(2,C)_{\RR}.$  
   
$\mathfrak{U}_q(sl(2,\mathbb{C})_{\RR})$ by definition is the quantum double of  
 $\mathfrak{U}_q(su(2)).$ As a  
result we can write $\mathfrak{U}_q(sl(2,\mathbb{C})_{\RR})=  
 \mathfrak{U}_q(su(2))\otimes \mathfrak{U}_q(su(2))^{*}$ as   
a vector space, where $\mathfrak{U}_q(su(2))^{*}$ denotes the restricted dual of  
$\mathfrak{U}_q(su(2)),$ i.e the Hopf algebra spanned by the matrix elements of representations   
contained in $Irr(\mathfrak{U}_q(su(2))).$  
 
A basis of $\mathfrak{U}_q(su(2))^{*}$ is the set of matrix elements of irreducible  
representations of $\mathfrak{U}_q(su(2))$, which we will denote by  
 ${\buildrel {{}_{B}} \over  g}{}^i_j, B\in \halfinteger, i,j=-B,...,B.$ 
 
It can be shown that $\mathfrak{U}_q(su(2))^{*}$ is isomorphic, as a star Hopf algebra, to the quantum envelopping algebra $\mathfrak{U}_q(an(2))$ where $an(2)$ is the  Lie algebra  of traceless  complex upper triangular 2$\times$2 matrices with real diagonal. 
 
\medskip
 
${\mathfrak U}_q(su(2))$ being a factorizable Hopf algebra, it is possible to give a nice  
generating family  of ${\mathfrak U}_q(su(2))$. Let us   
introduce,   
for each $I\in \halfinteger$ the elements  
 ${\buildrel I \over L}{}^{\!(\pm)}\in { End}(\CC^{d_I} )\otimes {\mathfrak U}_q(su(2))$  
 defined by  
 ${\buildrel I \over L}{}^{\!(\pm)}=({\buildrel I\over \pi} \otimes id)(R^{(\pm)}).$   
The matrix elements of ${\buildrel I \over L}{}^{\!(\pm)}$ when $I$ describes  $\halfinteger$ span the vector space  ${\mathfrak U}_q(su(2))$. 
    
The star Hopf algebra structure on $\mathfrak{U}_q(sl(2,\mathbb{C})_{\RR})$  is described in  details in \cite{BR1}.  Let us simply recall  that we have:  
\begin{eqnarray}  
&&\hskip -1.2cm{\buildrel I \over L}{}^{(\pm)}{}^{i}_{j} {\buildrel J \over L}{}^{(\pm)}{}^{k}_{l}=  
\sum_{Kmn} \Clebphi{I}{J}{K}{i}{k}{m}\!{\buildrel K\over L}{}^{(\pm)}{}^{m}_{n}\!\Clebpsi{I}{J}{K}{j}{l}{n}\;,\;\;\;\;\;\;{\buildrel {IJ} \over {R}}_{12} \; {\buildrel I \over L}{}^{(+)}_1 {\buildrel J \over L}{}^{(-)}_2 =  
 {\buildrel J \over L}{}^{(-)}_2 {\buildrel I \over L}{}^{(+)}_1 \; {\buildrel {IJ} \over {R}}_{12},\label{Rll}\\ 
&&\hskip -1.2cm{\buildrel I \over g}{}^{i}_{j} {\buildrel J \over g}{}^{k}_{l}=\sum_{Kmn}   
\Clebphi{I}{J}{K}{i}{k}{m}\!{\buildrel K \over g}{}^{m}_{n}\!\Clebpsi{I}{J}{K}{j}{l}{n}\;,\;\;\;\;\;\;\;\;\;\;\;\;\;\;\;\;\;\;\;\;\;\;\;  
{\buildrel {IJ} \over  R}{}^{(\pm)}_{12} {\buildrel I \over L}{}^{(\pm)}_1 {\buildrel J \over g}_2=  
{\buildrel J \over g}_2 {\buildrel I \over L}{}^{(\pm)}_1 {\buildrel {IJ} \over  R}{}^{(\pm)}_{12},\label{RLg}  
\end{eqnarray}  
\vskip -0.4cm  
\begin{eqnarray}  
&&\hskip -1cm\Delta({\buildrel I \over L}{}^{(\pm)}{}^a_b)=\sum_{c}{\buildrel I \over L}{}^{(\pm)}{}^c_b   
\otimes {\buildrel I \over L}{}^{(\pm)}{}^a_c \;,\;\;\;\;\;\;\;\;\;\;  
\Delta({\buildrel I \over g}{}^a_b)=\sum_{c}{\buildrel I \over  
g}{}^c_b \otimes {\buildrel I \over g}{}^a_c  ,\label{coproduitD}\\  
&&\hskip -1cm({\buildrel I \over L}{}^{(\pm)}{}^a_b)^{\star}= S^{-1}({\buildrel I \over L}{}^{(\mp)}{}^b_a)\;,\;\;\;\;\;\;\;\;\;\;\;\;\;\;\;\;\;\;\;\;\;\;\;({\buildrel I \over g}{}^a_b)^{\star}=S^{-1}({\buildrel I \over g}{}^b_a).\label{staralg}  
\end{eqnarray}  
The center of $ \mathfrak{U}_q(sl(2,\mathbb{C})_{\mathbb{R}})$ is a polynomial algebra in two variables $\Omega_+, \Omega_-$ and we have $\Omega_{\pm}= 
tr({\buildrel \onehalf \over \mu}{}^{-1}{\buildrel \onehalf \over L}{}^{(\mp)}{}^{-1} {\buildrel \onehalf \over g}).$ 
 
\medskip  
Principal Unitary Representations of   
$\mathfrak{U}_q(sl(2,\mathbb{C})_{\mathbb{R}})$ have been   
classified in \cite{Pu}, and the following description is given in \cite{BR1}.   
They are labelled by a couple of complex numbers $(X_0, X_1) \in {\mathbb S},$ where   
${\mathbb S}=\{(X_0, X_1) \in \CC^2\;\;/\;\;2X_0+1 =(m+i\rho)\;,\;2X_1+1 = (-m+i\rho),\;\;  
\;\; m \in {\frac{1}{2}}\ZZ, \rho \in ]-\frac{\pi}{\hbar}, \frac{ \pi}{\hbar}] \;\;\}.$ For $(X_0, X_1) \in {\mathbb S}$ we will often denote $m_X=X_0-X_1, i\rho_X=X_0+X_1+1.$

Let us denote by  ${\buildrel {\!\!{}_{(X_0X_1)}} \over {\Pi}}$ the principal  representation associated to $(X_0, X_1)$ and  $V(X_0X_1)$ the associated $\mathfrak{U}_q(sl(2,\mathbb{C})_{\mathbb{R}})$ module.  $V(X_0X_1)$ is a Harisch-Chandra module and we have    
 $V(X_0X_1)=\bigoplus_{C,C-\vert m_x\vert \in \NN} {\buildrel {C}\over V}$ as a ${\mathfrak U}_q(su(2))$ module. In term of the basis $\{{\buildrel {{}_{C}} \over e}_r(X_0X_1)={\buildrel {{}_{C}} \over e}_r, C\in \halfinteger, r=-C,\dots,C\}$ of the module $V(X_0X_1)$, the action of the generators of $\mathfrak{U}_q(sl(2,\mathbb{C})_{\mathbb{R}})$ is given by:  
\begin{eqnarray}  
{\buildrel {{}_{B}} \over  L}{}^{(\pm)}{}^i_j \;  {\buildrel {{}_{C}} \over e}_r &=  
&\sum_{n} {\buildrel {{}_{C}} \over e}_n\;\;  
 {\buildrel {{}_{BC}}\over  R}{}^{(\pm)}{}^{in}_{jr},\label{repl}\\  
{\buildrel {{}_{B}} \over  g}{}^i_j  \;  {\buildrel {{}_{C}} \over e}_r &=&   
\sum_{DEpx} {\buildrel {{}_{E}} \over e}_p\;\;   
\Clebphi{E}{B}{D}{p}{i}{x}\Clebpsi{B}{C}{D}{j}{r}{x} \Lambda^{BD}_{EC}(X_0X_1),\label{repg}  
\end{eqnarray}   
where the complex numbers $\Lambda^{BD}_{EC}(X_0X_1)$ have to verify certain constraints explained  in the appendix. It is a fondamental result that these coefficients can be expressed in  terms of $6j(1)$ as follows:  
\begin{eqnarray}  
&&\hskip -1cm\Lambda^{BC}_{AD}(X_0X_1)= \sum_{X_2} \!  
\sixjn{X_0}{X_1}{A}{B}{C}{X_2}{1} \!\sixjn{X_0}{X_1}{D}{B}{C}{X_2}{1}  
\frac{v_{X_2}}{v_{X_1}}   
\frac{v^{1/4}_{A}v^{1/4}_{D}}{v^{1/2}_{B}v^{1/2}_{C}}\label{formlamb6j}\\  
&&= \sum_{X_2} \!\sixjn{X_1}{X_0}{A}{B}{C}{X_2}{1}\!  
\sixjn{X_1}{X_0}{D}{B}{C}{X_2}{1} \frac{v_{X_0}}{v_{X_2}}   
\frac{v^{1/2}_{B}v^{1/2}_{C}}{v^{1/4}_{A}v^{1/4}_{D}}.  
\end{eqnarray} 
The action of the center on the module  $V(X_0X_1)$ is such that  ${\buildrel {\!\!{}_{(X_0X_1)}} \over {\Pi}}(\Omega_{\pm})=\omega_{\pm} id$ where $\omega_{+}=q^{2X_0+1}+q^{-2X_0-1}, \omega_{-}=q^{2X_1+1}+q^{-2X_1-1}.$ 
 
We can endow $V(X_0X_1)$ with a structure of pre-Hilbert space by defining the hermitian form $<.,.>$  such that the basis   $\{{\buildrel {{}_{C}} \over e}_r(X_0X_1), C\in \halfinteger, r=-C,\dots,C\}$ of $V(X_0X_1)$ is orthonormal.   
  
The representation  ${\buildrel {\!\!{}_{(X_0X_1)}} \over {\Pi}}$ is unitary in the sense that $\forall v, w\in V(X_0X_1),\;\; \forall a\in \mathfrak{U}_q(sl(2,\mathbb{C})_{\mathbb{R}}),$ $ <a^{*}v,w>=<v, aw>,$ 
 this last property being equivalent   to the relation: 
$\Lambda^{BC}_{AD}(X_1X_0)={\overline {\Lambda^{BC}_{AD}(X_0X_1)}}.$    
\medskip 
 
We will denote by $H(X_0X_1)$  the separable Hilbert space, completion of  $V(X_0X_1)$ which Hilbertian basis is    
 $\{{\buildrel {{}_{C}} \over e}_r(X_0X_1), C\in \halfinteger, r=-C,\dots,C\}.$  

\medskip 
  
The automorphism of the complex line is now playing a key role because of the following important result: the principal representations associated to  $(X_0, X_1)$ and to $(\underline{X_0}, \underline{X_1})$ are unitary equivalent.  
  
 \medskip

Let us now recall some basic facts about the algebra of functions on $SL_q(2,\CC)_{\RR}$ \cite{PW,Wo2, BR1}. We will use the notations of \cite{BR1}.  
 The space of compact supported functions on the quantum  
 Lorentz group, denoted $Fun_c(SL_q(2,\CC)_{\RR})$ is,  by definition,   
$Fun(SU_q(2)') \otimes \left( \bigoplus_{ I \in \onehalf \ZZ^{+}} End(\CC^{d_I})\right).$ This is a $C^{*}$ algebra without unit. It contains the dense *-subalgebra \cite{BR1}  
 $Fun_{cc}(SL_q(2,\CC)_{\RR})=Pol(SU_q(2)')\otimes \left( \bigoplus_{ I \in \onehalf \ZZ^{+}} End(\CC^{d_I})\right)$ which is a multiplier Hopf algebra \cite{VD}, and which can be understood as being the quantization of the algebra generated by  polynomials functions on $SU(2)$ and compact supported functions on $AN(2).$ 
  
$(\kE{C}{m}{n}{D}{p}{q})_{C,D,m,n,p,q}$ is a vector  basis of  $Fun_{cc}(SL_q(2,\CC)_{\RR})$ which is defined, for example,  by duality from the generators of the envelopping algebra:   
\begin{eqnarray}  
&&<{\buildrel {{}_{A}} \over  L}{}^{(\pm)}{}^i_j \otimes {\buildrel {{}_{B}} \over  g}{}^k_l, \;\kE{C}{m}{n}{D}{p}{q}>=  
{\buildrel {{}_{AC}}\over  R}{}^{(\pm)}{}^{im}_{jn}\delta_{B,D}\delta^{k}_{q}\delta^{p}_{l}.  
\end{eqnarray}  
We can describe completely the structure of the multiplier Hopf algebra in this basis;  
\begin{eqnarray}  
&&\hskip -1.5cm\kE{A}{i}{j}{B}{k}{l} \cdot \kE{C}{m}{n}{D}{p}{q}= \sum_{Frs} \Clebphi{A}{C}{F}{i}{m}{r}\Clebpsi{A}{C}{F}{j}{n}{s} \delta^p_l   
\delta_{B,D} \kE{F}{r}{s}{B}{k}{q}\nonumber\\  
&&\hskip -1.5cm\Delta(\kE{A}{i}{j}{B}{k}{l})={\cal F}^{-1}_{23}\;  
(\!\!\!\!\sum_{\substack{{C,D,m,}\\{p,q,r,s}}}\!\!\!\!\Clebphi{C}{D}{B}{q}{s}{l}\!\!\Clebpsi{C}{D}{B}{p}{r}{k}\!\kE{A}{i}{m}{C}{p}{q} \!\otimes \!  
 \kE{A}{m}{j}{D}{r}{s}\;)\;{\cal F}_{23}\nonumber\\  
&&\hskip -1.5cm\epsilon(\kE{A}{i}{j}{B}{k}{l})=\delta^i_j \delta_{B,0}\;\;\;\;\;\;\;(\kE{A}{i}{j}{B}{k}{l})^{\star}=S^{-1}({\buildrel A \over k}{}^j_i) \! \otimes   
\! {\buildrel B \over E}{}^l_{k}\;\;\;\;\mbox{with}\;\;  
{\cal F}^{-1}_{12}=\sum_{J,x,y}\!\!{\buildrel J \over E}{}^x_y \! \otimes \! S^{-1}\!({\buildrel J \over k}{}^y_x). 
\end{eqnarray}  
The space of right and left invariant  linear forms (also called Haar measures) on $Fun_c(SL_q(2,\CC)_{\RR})$  is a vector space of dimension one and we will pick one element $h$, which is defined by:  
\begin{eqnarray}  
h(\kE{A}{i}{j}{B}{m}{l})= \delta_{A,0}({\buildrel {{}_{B}} \over \mu}{}^{-1})^m_l \qd{B}.\label{haar}  
\end{eqnarray}  
Using the $L^2$ norm, $\vert \vert \; a \; \vert \vert_{L^2}=h(a^* a)^{\onehalf}$, we can complete the space  $Fun_c(SL_q(2,\CC)_{\RR})$   into the Hilbert space of $L^2$ functions on the quantum Lorentz group, denoted $L^2(SL_q(2,\CC)_{\RR})$.  
  
 $Fun_{cc}(SL_q(2,\CC)_{\RR})$ is a multiplier Hopf algebra with basis $(u_\alpha)=((\kE{C}{m}{n}{D}{p}{q})_{C,D,m,n,p,q})$. The restricted dual of $Fun_{cc}(SL_q(2,\CC)_{\RR})$, denoted  ${\tilde \mathfrak U}_q(sl_q(2,\CC)_{\RR}),$   is the vector space spanned by the dual basis 
 $(u^{\alpha})= 
({\buildrel {{}_{C}}\over X }{}^{n}_m \otimes {\buildrel {{}_{D}}\over g}{}^q_p).$  It is also, by duality, a multiplier Hopf algebra and   ${\mathfrak U}_q(sl(2,\CC)_{\RR} )$ is included as an algebra in the  multiplier algebra of  ${\tilde \mathfrak U}_q(sl(2,\CC)_{\RR} )$.    
If $\Pi$ is the  principal representation  of ${\mathfrak U}_q(sl(2,\CC)_{\RR} )$, acting on $V(Z_0Z_1),$ it is possible \cite{BR1} to associate to it a unique representation ${\tilde \Pi}$ of   ${\tilde \mathfrak U}_q(sl(2,\CC)_{\RR} ),$ acting on $V(Z_0Z_1)$, such that  
 ${\tilde \Pi}({\buildrel {{}_{C}}\over X }{}^{n}_m)({\buildrel {{}_{D}} \over e}_r)=\delta_{C,D}\delta^{n}_{r}{\buildrel {{}_{D}} \over e}_m.$ 
  
We define for all $\psi$ element of  $Fun_{cc}(SL_q(2,\CC)_{\RR})$, the operator  $\Pi(\psi)=\sum_I {\tilde\Pi}(u^I)\; h(u_I \; \psi).$ It is easy to show that  $\Pi(\psi)$ is of finite rank.  
  
 \medskip 
 
If $f$ is a function on ${\mathbb S}$, we will sometimes write  $f(m_X,\rho_X)$ instead of   $f(X_0X_1).$    
The Plancherel  formula can be written as:  
\begin{eqnarray}   
&&\hskip -0.6cm \forall \psi\in Fun_{cc}(SL_q(2,\CC)_{\RR}),\;\; \vert \vert \; \psi \;  
 \vert \vert_{L^2}^2  = \; \int d{\cal P}(X_0X_1) \;\;  
{\rm tr} \left( \; {\buildrel {\!\!{}_{(X_0X_1)}}  
 \over {\Pi}}\!\!(\mu^{-1}){\buildrel {\!\!{}_{(X_0X_1)}}   
\over {\Pi}}\!\!\!\!(\psi){\buildrel {\!\!{}_{(X_0X_1)}} \over {\Pi}}\!\!\!\!(\psi)^{\dagger} \; \right),  
\end{eqnarray}  
where  we have denoted    
\begin{eqnarray}  
&&\hskip -0.5cm\int d{\cal P}(X_0X_1) \; f(X_0X_1)=  
\!\!\!\sum_{m \in \onehalf \ZZ} \; \int_{-\frac{\pi}{\hbar}}^{\frac{\pi}{\hbar}} \!\!\!d\rho\;\;  
{\cal P}(m,\rho)\; f(m,\rho)\nonumber\\  
&&\mbox{with}\;\;\;\;\;{\cal P}(m,\rho)=  
\frac{\hbar}{2\pi}(q-q^{-1})^2\;[m+i\rho][m-i\rho].\nonumber  
\end{eqnarray}  
  
The proof of this theorem is purely combinatorial and uses as a central tool the following identity on Fourier coefficients of the Laurent polynomials $\Lambda^{BC}_{AA}$:  
\begin{equation}  
\int d{\cal P}(X_0X_1) \; \Lambda^{BC}_{AA}(X_0X_1)\; = \;\delta_{B,0}\delta_{A,C}\qd{A}\label{intlambda}.  
\end{equation}

A Plancherel Theorem for $L^2$ functions has been proved \cite{BR1}, it follows easily from the previous result and from the following  lemma, which will be useful later on 
    
\begin{lemma}

 The only function  $f:{\mathbb S}\rightarrow \CC$ satisfying the following conditions:  
  
1)$\forall m\in \onehalf \ZZ, \forall \rho\in ]-\frac{\pi}{\hbar},\frac{\pi}{\hbar}] , f(m,\rho)=f(-m,-\rho),$  
  
2)$\forall m\in  \onehalf \ZZ, f(m,.)$ is a $L^2$ function on $]-\frac{\pi}{\hbar},\frac{\pi}{\hbar}], $  
  
3)$\exists m_0\in  \onehalf \NN, \forall m, \vert m \vert > m_0, f(m,.)=0,$  
  
4) $\exists  A,D \in \halfinteger \cap [\vert m_0 \vert, +\infty[,  
 \forall B,C \in \halfinteger,\;\;\int d{\cal P}(X_0X_1)\; f(X_0X_1) \; \Lambda^{BC}_{AD}(X_0X_1)\; = \;0,$   
  
\noindent is the nul function.  
\end{lemma}

\section{Intertwiners associated to unitary representations of the quantum Lorentz group}  

The aim of this section is to give explicit formulae of  the intertwiners between   
the representation   
${\buildrel {\!\!{}_{(X_0X_1)}} \over {\Pi}}\!\!\!\otimes \!\!{\buildrel {\!\!{}_{(Y_0Y_1)}} \over {\Pi}}$ and the representation  ${\buildrel {\!\!{}_{(Z_0Z_1)}} \over {\Pi}},$ in terms  of complex continuations of $6j$ symbols of $\mathfrak{U}_q(su(2)).$  
\begin{proposition}  
Let $\GPhin{X_0X_1}{Y_0Y_1}{Z_0Z_1}:V(X_0X_1)\otimes V(Y_0Y_1)\rightarrow V(Z_0Z_1)$ be a   
$\mathfrak{U}_q(sl(2,\CC)_{\RR})$ intertwiner. We necessarily have:  
\begin{equation}  
\hskip -0.5cm\GPhin{X_0X_1}{Y_0Y_1}{Z_0Z_1}(\ontop{A}{e}_i(X_0X_1)\otimes  
\ontop{B}{e}_j(Y_0Y_1))=\!\sum_{C,k}\ontop{C}{e}_k(Z_0Z_1)\!\Clebpsi{A}{B}{C}{i}{j}{k}\!\ElemRed{A}{B}{C}{X_0X_1}{Y_0Y_1}{Z_0Z_1},\label{defreducedelement}  
\end{equation}  
where the coefficients $\ElemRed{A}{B}{C}{X_0X_1}{Y_0Y_1}{Z_0Z_1}$, called ``reduced elements'',  are complex numbers. \\  
Inversely such a  map defines an intertwiner if and only if the following conditions on the reduced elements are  satisfied:  
$\forall A,B,C\in \onehalf \ZZ^+, Y^{(0)}_{(A,B,C)}=Y^{(1)}_{(A,m_X)}=Y^{(1)}_{(B,m_Y)}=1,$ 
\begin{eqnarray}  
&&\hskip -1cm\sum_{QRSP} \ElemRed{R}{P}{T}{X_0X_1}{Y_0Y_1}{Z_0Z_1}  
\Lambda^{US}_{RA}(X_0X_1)\Lambda^{UQ}_{PB}(Y_0Y_1)  
\sixjn{U}{W}{T}{R}{P}{Q}{0}\sixjn{R}{W}{Q}{B}{U}{S}{0}\!\!\times\nonumber\\  
&&\times\sixjn{B}{W}{S}{U}{A}{C}{0}=\ElemRed{A}{B}{C}{X_0X_1}{Y_0Y_1}{Z_0Z_1}\Lambda^{UW}_{TC}(Z_0Z_1)\label{inter1}.  
\end{eqnarray}  
\end{proposition}  
  
\Proof   
 
The necessary condition comes from the fact that (\ref{defreducedelement}) is equivalent to the property that  
 $\GPhin{X_0X_1}{Y_0Y_1}{Z_0Z_1}$ is an intertwiner of $\mathfrak{U}_q(su(2))$ module. 
Such a map is an intertwiner of $\mathfrak{U}_q(an(2))$ module, if moreover,   
\vskip -0.5cm  
\begin{eqnarray}  
\GPhin{X_0X_1}{Y_0Y_1}{Z_0Z_1}(\Delta({\buildrel U \over g}{}^{m}_{n})    
(\ontop{A}{e}_i(X_0X_1)\otimes \ontop{B}{e}_j(Y_0Y_1)))={\buildrel U \over g}{}^{m}_{n}    
\GPhin{X_0X_1}{Y_0Y_1}{Z_0Z_1}(\ontop{A}{e}_i(X_0X_1)\otimes \ontop{B}{e}_j(Y_0Y_1)).\nonumber  
\end{eqnarray}  
This last condition can be rewritten as (\ref{inter1}) using (\ref{repg})(\ref{ortho1cont})(\ref{vertexIRF}). This concludes the proof.  
\eoProof  

\medskip
  
Remark: A very important point is that $\GPhin{X_0X_1}{Y_0Y_1}{Z_0Z_1}$ maps the algebraic tensor product of the two domains $V(X_0X_1)\otimes V(Y_0Y_1)$ to  $V(Z_0Z_1).$ As a result the sum $(\ref{defreducedelement})$ is finite.  
It can be seen that there are no non zero intertwiner from $V(Z_0Z_1)$ to the algebraic tensor product $V(X_0X_1)\otimes V(Y_0Y_1).$  
  
\begin{lemma}  
The space of $\mathfrak{U}_q(sl(2,\CC)_{\RR})$ intertwiners from the module $V(X_0X_1)\otimes V(Y_0Y_1)$ to the module $V(Z_0Z_1)$ is of dimension $0$ or $1.$  
\end{lemma}  
 
\Proof  
 
An elementary proof is obtained by analyzing the rank of the system of linear equations (\ref{inter1}). Using the isomorphism of module between $V(X_0X_1)$ and $V({\underline {X_0}},{\underline {X_1}})$, we can always assume that $m_X, m_Y, m_Z$ are non negative. 
It is easy to show that we can always assume that $m_Z\leq m_X+m_Y.$ If not, we use  
the fact that  the space of intertwiners from $V(X_0X_1)\otimes V(Y_0Y_1)$ to $V(Z_0Z_1)$  
is in one to one correspondence with the space of intertwiners from  
$V(Z_0Z_1)\otimes V({\bar X_0}{\bar X_1})$ to $V(Y_0Y_1)$ to exchange $m_Y$ and $m_Z$.  This one 
 to one correspondence  easily follows from the isomorphism $V(X_0X_1)^{*}\approx{\bar  V}(X_0X_1) \approx V({\bar X_0}{\bar X_1})$ and where we denoted by $V(X_0X_1)^{*}$ the restricted dual of the Harisch-Chandra module $V(X_0X_1).$ 
 
This system of linear equations is equivalent to the subsystem where we have chosen $U=\onehalf$, because  $\mathfrak{U}_q(an(2))$ is generated as an algebra by  ${\buildrel U \over g}{}^{m}_{n}$ where $U=\onehalf.$ 
This system is therefore  equivalent to the following one: 
$\forall \sigma,\tau \in \{\onehalf,-\onehalf  \},$ 
 
\begin{eqnarray}  
&&\hskip -1cm\sum_{\epsilon,\sigma,\mu,\tau\in\{-\onehalf,\onehalf\}} \ElemRed{\!\!A+$\nu+\rho$}{B+$\epsilon+\mu$}{\!C+$\sigma+\tau$}{X_0X_1}{Y_0Y_1}{Z_0Z_1} 
\Lambda^{\onehalf \;\;A+\nu}_{A+\nu+\rho\;A}(X_0X_1) 
\Lambda^{\onehalf \;\;B+\epsilon}_{B+\epsilon+\mu\;B}(Y_0Y_1)\times\nonumber\\ 
&&\times\sixjn{A\!+\!\nu\!+\!\rho}{C\!+\!\sigma}{B\!+\!\epsilon}{\onehalf}{B\!+\!\epsilon\!+\!\mu}{C\!+\!\sigma\!+\!\tau}{0} 
\sixjn{C+\sigma}{A\!+\!\nu\!+\!\rho}{B\!+\!\epsilon}{\onehalf}{B}{A\!+\!\nu}{0}\!\!\sixjn{B}{C\!+\!\sigma}{A\!+\!\nu}{\onehalf}{A}{C}{0}= 
\nonumber\\ 
&&\hskip 1.5cm= 
\ElemRed{A}{B}{C}{X_0X_1}{Y_0Y_1}{Z_0Z_1} 
\Lambda^{{\onehalf}{C+\sigma}}_{{C+\sigma+\tau}{C}}(Z_0Z_1) 
\label{inter1/2}.  
\end{eqnarray}  
We will denote by  $S(\sigma,\tau)$ this system of equations. 
 
It is easy to show that the system $S(-\onehalf,-\onehalf)$ completely determines the reduced elements at the point $(A,B,C)$ in terms of the reduced element of the points $(A',B',C-1)$ with $A'=A+\epsilon, B'=B+\nu$ with $\epsilon,\nu \in \{\onehalf,-\onehalf  \}$. 
  Therefore the rank of the system is less than the rank of the vectors which components 
 are the reduced elements at the points $(A,B,m_Z).$\\ 
Let $\Delta=\{ (A,B)\in(\onehalf\ZZ^{+})^2, Y^{(0)}_{(A,B,m_Z)}=1\}, P=\{(A,B)\in(\onehalf\ZZ^{+})^2, A-m_X\in \NN\}, 
 Q=\{(A,B)\in\onehalf(\ZZ^{+})^2, B-m_Y\in \NN\}.$ 
$\Delta\cap P\cap Q$ is the intersection of a lattice with a convex set which boundary consits in 3 segments and two half-lines. In the case where $m_Z\leq m_X+m_Y$ one of this segment degenerate to the point $p=(m_X,m_Y).$ 
It is easy to show, by a direct computation, that the system $S(\onehalf,-\onehalf)$ and  
$S(-\onehalf,\onehalf)$ are independent. As a result we can take linear combination of these two  
systems to have linear combinations of reduced elements at $(A,B,m_Z)$ involving only $8$ points and  
not 9. It is easy to show that the use of both of these systems determine uniquely the reduced elements at the point $(A,B,m_Z)$ in terms of the reduced elements at the point $p, p+(1,0), p+(0,1), p+(1,1).$ 
As a result the system is of rank less than four. 
But the reduced elements at the point  $p, p+(1,0), p+(0,1), p+(1,1)$ are solutions of three systems of linear equations which can be shown to be independent:\\ 
$S(\onehalf,-\onehalf)$ at $(m_X,m_Y,m_Z),$ 
$S(-\onehalf,\onehalf)$ at $(m_X,m_Y,m_Z),$ and $S(\onehalf,\onehalf)$ at $(m_X,m_Y,m_Z-1).$ 
As a result the rank of the system is of dimension less than one. 
\eoProof

\begin{theorem}  
Assume that the numbers $\rho_X,\rho_Y,\rho_Z$ and  $\epsilon_X\rho_X+\epsilon_Y\rho_Y+\epsilon_Z\rho_Z$ with $\epsilon_X,\epsilon_Y,\epsilon_Z\in \{-1,1\}$ are non zero.  
 The space of $\mathfrak{U}_q(sl(2,\CC)_{\RR})$ intertwiners from the module $V(X_0X_1)\otimes V(Y_0Y_1)$ to the module $V(Z_0Z_1)$ is of dimension $0$ if and only if $m_X+m_Y+m_Z\notin\ZZ.$\\  
If $m_X+m_Y+m_Z\in\ZZ$ it is a one dimensional space which admits a non zero element, whose reduced  element, given in terms of $6j(1)$ and $6j(3),$ is   
\begin{eqnarray}  
&&\hskip -1cm\ElemRed{R}{P}{T}{\!X_0X_1}{\!Y_0Y_1}{Z_0Z_1}\!\!=\!\!\sum_{X_2}\!\sixjn{\!Y_0}{\!X_0}{\!Z_0\!}{\!T}{\!Z_1}{\!X_2\!}{3}\!\!\!  
\sixjn{\!X_1\!}{\!X_0}{\!R}{\!T\!}{\!P\!}{\!X_2}{1}\!\!\!\sixjn{\!Z_1}{\!X_2}{\!Y_0\!}{\!P}{\!Y_1}{\!X_1\!}{3}\!\!\!\!\!  
\frac{v^{1/4}_{P}v^{1/4}_{T}v^{1/4}_{X_0}v^{1/4}_{X_1}[d_P]^{1/2}}{v^{1/4}_{R}v^{1/2}_{X_2}e^{i\pi P}}. \nonumber\\ 
&& 
\label{form6jprol1}  
\end{eqnarray} 
\vskip -0.3cm  
We will denote by $\GPhin{X_0X_1}{Y_0Y_1}{Z_0Z_1}$ the  associated intertwining operator.  
 
The reduced elements satisfy the condition: 
 
$\forall C, C-\vert m_Z\vert\in \NN, \exists A,B\in \onehalf {\ZZ}^{+},  
\ElemRed{A}{B}{C}{X_0X_1}{Y_0Y_1}{Z_0Z_1}\not=0.$ 
 
\end{theorem}  
  
\Proof  
 
We have assumed that $\rho_X,\rho_Y,\rho_Z$ and $\epsilon_X\rho_X+\epsilon_Y\rho_Y+\epsilon_Z\rho_Z$ with $\epsilon_X,\epsilon_Y,\epsilon_Z\in \{-1,1\}$ are non zero in order that all the $6j(1)$ and $6j(3)$ are well defined in the expression (\ref{form6jprol1}). But this condition can nevertheless be removed by normalizing the reduced elements by a function $F(Z_0Z_1,Y_0Y_1,Z_0Z_1)$ which remove the singularities  of this expression. 
 
Using only  polynomial identities on continuation of  $6j$, we will show that  the left handside of (\ref{inter1}) and the righthandside of (\ref{inter1}) are equal if we take the ansatz (\ref{form6jprol1}) for the reduced element.   
  
The lefthandside of (\ref{inter1})=  
\begin{eqnarray*}  
&&\hskip -0.6cm=e^{-i\pi P}\sqrt{[d_P]}  
\!\!\!\!\!\sum_{QRSPX_2X_3Y_2}\!\!\!\!  
\frac{v^{1/2}_{P}v^{1/4}_{T}  
v^{1/4}_{A}v^{1/4}_{B}v^{1/4}_{X_0}v_{X_3}v_{Y_2}}  
{v_{U}v^{1/2}_{X_2}v^{1/2}_{S}v^{1/2}_{Q}v^{3/4}_{X_1}v_{Y_1}}  
\sixjn{Y_0}{X_0}{Z_0}{T}{Z_1}{X_2}{3}  
\!\!\sixjn{X_1}{X_0}{R}{T}{P}{X_2}{1}\nonumber\\  
&&\times\!\!\sixjn{Z_1}{X_2}{Y_0}{P}{Y_1}{X_1}{3}  
\!\!\sixjn{X_0}{X_1}{R}{U}{S}{X_3}{1}  
\!\!\sixjn{X_0}{X_1}{A}{U}{S}{X_3}{1}  
\!\!\sixjn{U}{W}{T}{R}{P}{Q}{0}  
\!\!\sixjn{R}{W}{Q}{B}{U}{S}{0}\nonumber\\  
&&\times\!\!\sixjn{B}{W}{S}{U}{A}{V}{0}  
\!\!\sixjn{Y_0}{Y_1}{P}{U}{Q}{Y_2}{1}  
\!\!\sixjn{Y_0}{Y_1}{B}{U}{Q}{Y_2}{1}.  
\end{eqnarray*}  
We apply first the pentagonal identity (\ref{penta3cont}) to the second and sixth $6j$, then   
the pentagonal identity (\ref{penta3cont}) to the fourth and the seventh $6j$ and finally we realize the sum on $R$ by applying the orthogonality (\ref{ortho1cont}), to obtain  
\vskip -0.3cm  
\begin{eqnarray*}  
&&\hskip -0.6cm=e^{-i\pi P}\sqrt{[d_P]}  
\!\!\!\!\!\sum_{QSPX_2X_3Y_2X_4}\!\!\!\!  
\frac{v^{1/2}_{P}v^{1/4}_{T}  
v^{1/4}_{A}v^{1/4}_{B}v^{1/4}_{X_0}v_{X_3}v_{Y_2}}  
{v_{U}v^{1/2}_{X_2}v^{1/2}_{S}v^{1/2}_{Q}v^{3/4}_{X_1}v_{Y_1}}  
\sixjn{Y_0}{X_0}{Z_0}{T}{Z_1}{X_2}{3}  
\!\!\sixjn{X_1}{X_4}{Q}{U}{P}{X_2}{0}\nonumber 
\end{eqnarray*}  
\vskip -0.3cm  
\begin{eqnarray*}  
&&\times\!\!\sixjn{Z_1}{X_2}{Y_0}{P}{Y_1}{X_1}{3}  
\!\!\sixjn{X_0}{X_2}{T}{U}{W}{X_4}{0}  
\!\!\sixjn{X_0}{X_1}{A}{U}{S}{X_3}{1}  
\!\!\sixjn{X_1}{X_4}{Q}{B}{U}{X_3}{1}  
\!\!\sixjn{X_0}{X_3}{S}{B}{W}{X_4}{1}\nonumber\\  
&&\times\!\!\sixjn{B}{W}{S}{U}{A}{V}{0}  
\!\!\sixjn{Y_0}{Y_1}{P}{U}{Q}{Y_2}{1}  
\!\!\sixjn{Y_0}{Y_1}{B}{U}{Q}{Y_2}{1}.  
\end{eqnarray*}  
Now we apply the symetries (\ref{symcont})(\ref{sym6j3}) and transform the sum over $P$ of the second, the third, and the ninth $6j$ using the hexagonal identity (\ref{hexa6j3}), to obtain  
\vskip -0.3cm  
\begin{eqnarray*}  
&&\hskip -0.6cm=e^{-i\pi Q}\sqrt{[d_Q]}  
\!\!\!\!\!\sum_{QSZ_3X_2X_3Y_2X_4}\!\!\!\!  
\frac{v^{1/4}_{T}  
v^{1/4}_{A}v^{1/4}_{B}v^{1/4}_{X_0}v_{X_3}v^{1/2}_{Y_2}v^{1/2}_{Z_3}}  
{v_{U}v^{1/2}_{X_4}v^{1/2}_{S}v^{1/2}_{Z_1}v^{3/4}_{X_1}v^{1/2}_{Y_1}}  
\sixjn{Y_0}{X_0}{Z_0}{T}{Z_1}{X_2}{3}  
\!\!\sixjn{Y_0}{Z_1}{X_2}{U}{X_4}{Z_3}{3}  
\nonumber 
\end{eqnarray*}  
\vskip -0.3cm  
\begin{eqnarray*} 
&&\times\!\!\sixjn{X_1}{Z_1}{Y_1}{U}{Y_2}{Z_3}{3}  
\!\!\sixjn{X_0}{X_2}{T}{U}{W}{X_4}{0}  
\!\!\sixjn{X_0}{X_1}{A}{U}{S}{X_3}{1}  
\!\!\sixjn{X_1}{X_4}{Q}{B}{U}{X_3}{1}  
\!\!\sixjn{X_0}{X_3}{S}{B}{W}{X_4}{1}\nonumber\\  
&&\times\!\!\sixjn{B}{W}{S}{U}{A}{V}{0}  
\!\!\sixjn{Z_3}{Y_0}{X_4}{Q}{X_1}{Y_2}{3}  
\!\!\sixjn{Y_0}{Y_1}{B}{U}{Q}{Y_2}{1}  
(e^{i \pi}q)^{(X_4+Y_1-X_1-Y_0)}\frac{\nu_1(d_{X_1})\nu_1(d_{Y_0})}{\nu_1(d_{X_4})\nu_1(d_{Y_1})}.  
\end{eqnarray*}  
We realize the sum over $Q$ of the sixth, ninth and tenth $6j$ using the pentagonal identity (\ref{penta1sixj3}) and symetries (\ref{symcont}), the sum over $X_2$ of the first, second and fourth $6j$ using the pentagonal identity (\ref{penta1sixj3}) and transform the sum over $S$ of the fifth, seventh and eighth $6j$ according to the hexagonal identity (\ref{hexacont}) to obtain  
\vskip -0.3cm  
\begin{eqnarray*}  
&&\hskip -0.6cm=e^{-i\pi B}\sqrt{[d_B]}  
\!\!\!\!\!\sum_{Z_3X_5X_3Y_2X_4}\!\!\!\!  
\frac{v^{1/4}_{T}v^{1/2}_{V}  
v^{1/4}_{B}v^{1/4}_{X_0}v^{1/2}_{X_3}v^{1/2}_{Y_2}v^{1/2}_{Z_3}}  
{v^{1/2}_{W}v^{1/4}_{A}v_{U}v^{1/2}_{X_5}v^{1/2}_{Z_1}v^{1/4}_{X_1}v^{1/2}_{Y_1}}  
\sixjn{X_1}{Z_1}{Y_1}{U}{Y_2}{Z_3}{3}  
\!\!\sixjn{Z_3}{Y_1}{X_3}{U}{X_1}{Y_2}{3}\nonumber\\  
&&\times\!\!\sixjn{Z_3}{Y_0}{X_4}{B}{X_3}{Y_1}{3}  
\sixjn{Y_0}{Z_3}{X_4}{W}{X_0}{Z_0}{3}  
\!\!\sixjn{Z_0}{Z_1}{T}{U}{W}{Z_3}{1}  
\!\!\sixjn{B}{X_4}{X_3}{U}{X_1}{X_5}{1}  
\!\!\sixjn{X_0}{X_1}{A}{B}{V}{X_5}{1}\nonumber\\  
&&\times\!\!\sixjn{X_0}{X_5}{V}{U}{W}{X_4}{1}  
(e^{i \pi}q)^{(X_4+2Y_1-X_1-Y_0-Y_2)}  
\frac{\nu_1(d_{X_1})\nu_1(d_{Y_0})\nu_1(d_{Y_2})}{\nu_1(d_{X_4})\nu_1(d_{Y_1})^2}.  
\end{eqnarray*}  
Then, we realize the sum over $Y_2$ of the first and second $6j$ by using the Racah identity (\ref{racah6j3}) and symetries (\ref{sym6j3}) to find  
\vskip -0.3cm  
\begin{eqnarray*}  
&&\hskip -0.6cm=e^{-i\pi B}\sqrt{[d_B]}  
\!\!\!\!\!\sum_{Z_3X_5X_3X_4}\!\!\!\!  
\frac{v^{1/4}_{T}v^{1/2}_{V}v_{Z_3}  
v^{1/4}_{B}v^{1/4}_{X_0}v^{1/4}_{X_1}}  
{v^{1/2}_{U}v^{1/4}_{A}v_{Z_1}v^{1/2}_{W}v^{1/2}_{X_5}}  
\sixjn{Y_1}{Z_3}{X_3}{U}{X_1}{Z_1}{3}  
\!\!\sixjn{Z_3}{Y_0}{X_4}{B}{X_3}{Y_1}{3}\nonumber\\  
&&\times\!\!  
\sixjn{Y_0}{Z_3}{X_4}{W}{X_0}{Z_0}{3}  
\!\!\sixjn{Z_0}{Z_1}{T}{U}{W}{Z_3}{1}  
\!\!\sixjn{B}{X_4}{X_3}{U}{X_1}{X_5}{1}  
\!\!\sixjn{X_0}{X_1}{A}{B}{V}{X_5}{1}  
\sixjn{X_0}{X_5}{V}{U}{W}{X_4}{1}  
\nonumber\\  
&&\times\!  
(e^{i \pi}q)^{(X_4+Y_1+Z_1-X_1-Y_0-Z_3)}  
\frac{\nu_1(d_{X_1})\nu_1(d_{Y_0})\nu_1(d_{Z_3})}{\nu_1(d_{X_4})\nu_1(d_{Y_1})\nu_1(d_{Z_1})}.  
\end{eqnarray*}  
Now, we realize the sum over $X_3$ of the first, second and fifth $6j$ by using the pentagonal identity (\ref{penta2sixj3}) and symetries (\ref{sym6j3}) to find  
\vskip -0.3cm  
\begin{eqnarray*}  
&&\hskip -1.7cm=e^{-i\pi B}\sqrt{[d_B]}  
\!\sum_{Z_3X_5X_4}\!  
\frac{v^{1/4}_{T}v^{1/2}_{V}v_{Z_3}  
v^{1/4}_{B}v^{1/4}_{X_0}v^{1/4}_{X_1}}  
{v^{1/2}_{U}v^{1/4}_{A}v_{Z_1}v^{1/2}_{W}v^{1/2}_{X_5}}  
\sixjn{Y_0}{Z_1}{X_5}{U}{X_4}{Z_3}{3}  
\!\!\sixjn{Z_1}{Y_1}{X_1}{B}{X_5}{Y_0}{3}\nonumber\\  
&&\times  
\sixjn{Y_0}{Z_3}{X_4}{W}{X_0}{Z_0}{3}  
\!\!\sixjn{Z_0}{Z_1}{T}{U}{W}{Z_3}{1}  
\!\!\sixjn{X_0}{X_1}{A}{B}{V}{X_5}{1}  
\sixjn{X_0}{X_5}{V}{U}{W}{X_4}{1}.  
\end{eqnarray*}  
\vskip -0.2cm  
Finally we realize the sum over $X_4$ of the first, third and sixth $6j$ using the pentagonal equation (\ref{penta36j3}) to conclude  
\begin{eqnarray*}  
&&\hskip -0.9cm=e^{-i\pi B}\sqrt{[d_B]}\sum_{X_5}  
\frac{v^{1/4}_{V}  
v^{1/4}_{B}v^{1/4}_{X_0}v^{1/4}_{X_1}}  
{v^{1/4}_{A}v^{1/2}_{X_5}}  
\sixjn{Y_0}{Z_1}{X_5}{V}{X_0}{Z_0}{3}  
\!\!\sixjn{X_0}{X_1}{A}{B}{V}{X_5}{1}  
\!\!\sixjn{Z_1}{Y_1}{X_1}{B}{X_5}{Y_0}{3}\nonumber\\  
&&\times\!\!\sum_{Z_3} \frac{v^{1/4}_{T}v^{1/4}_{V}v_{Z_3}}  
{v^{1/2}_{U}v^{1/2}_{W}v_{Z_1}}  
\!\!\sixjn{Z_0}{Z_1}{T}{U}{W}{Z_3}{1}\sixjn{Z_1}{Z_0}{V}{W}{U}{Z_3}{1}\nonumber\\  
&&\hskip -0.9cm=\mbox{the righthandside of (\ref{inter1}).}\nonumber  
\end{eqnarray*}  
This concludes the proof that the expression of the reduced elements defines an intertwiner operator. 
 
Let us now prove the final part of the theorem. 
Let us fix $C$ such that $(C-\vert m_Z\vert) \in \NN$ and assume that $\forall A, B, \;\ElemRed{A}{B}{C}{X_0X_1}{Y_0Y_1}{Z_0Z_1}\!\!=0.$ In this event, by multiplying the reduced element by $\sixjn{X_1}{X_0}{A}{C}{B}{X_3}{1}$, and  
summing over $A$, we would obtain, after the use of the orthogonality relation on $6j(1),$ $\forall B, \forall X_2,\;\; \sixjn{Y_0}{X_0}{Z_0}{C}{Z_1}{X_2}{3}\!\!\sixjn{Z_1}{X_2}{Y_0}{B}{Y_1}{X_1}{3}\!\!\!=0.$ 
But from the orthogonality relation on $6j(3),$ $\exists X_2 \; / \;\sixjn{Y_0}{X_0}{Z_0}{C}{Z_1}{X_2}{3}\!\!\!\not=0.$ As a result we  obtain that $\sixjn{Z_1}{X_2}{Y_0}{B}{Y_1}{X_1}{3}\!\!=0$ for all $B$ subject to the selection rules. 
From the behaviour  of $\sixjn{Z_1}{X_2}{Y_0}{B}{Y_1}{X_1}{3}$ when $B$ is large (see Eq. (\ref{asymptotic6j3})), we obtain a contradiction. 
As a result the statement of the theorem holds true and implies in particular that the intertwiner is non zero. 
\eoProof  
  
\section{Alternative Construction of Intertwiners in terms of the Quantum Haar Measure}  
  
Let us define the linear map   
${\hat \Phi}[X_0X_1,Y_0Y_1]:V(X_0X_1)\otimes V(Y_0Y_1)\rightarrow \int^{\oplus} d(Z_0Z_1)H(Z_0Z_1)$ where   
$\forall w\in  V(X_0X_1)\otimes V(Y_0Y_1)  
,{\hat \Phi}[X_0X_1,Y_0Y_1](w)$ is the family of functions defined by   
  
\noindent ${\hat \Phi}[X_0X_1,Y_0Y_1](w)(Z_0Z_1)=N^{{}_{(X_0X_1)\;(Y_0Y_1)}}_{\;\;\;{}_{(Z_0Z_1)}}\GPhin{X_0X_1}{Y_0Y_1}{Z_0Z_1}(w),$  
and $N^{{}_{(X_0X_1)\;(Y_0Y_1)}}_{\;\;\;{}_{(Z_0Z_1)}}$ are  complex numbers depending on $X_0, X_1, Y_0, Y_1, Z_0, Z_1.$  
 
We want to find $N$ such that ${\hat \Phi}[X_0X_1,Y_0Y_1]$ is an isometry.  
As explained in the introduction, this is a delicate  problem which requires  another description of the space of intertwiners, where this isometry property is a direct consequence of Plancherel theorem. 
  
\begin{theorem}  
Let $l\in V(X_0X_1), l'\in V(Y_0Y_1), v''\in V(Z_0Z_1)$,   
the following operator is well defined, and is an   intertwiner of  $\mathfrak{U}_q(sl(2,\CC)_{\RR})$ module:  
\begin{eqnarray}  
&&\hskip -1cm \Upsilon_{\;\;{}_{(Z_0Z_1)}}^{{}_{(X_0X_1)}{}_{(Y_0Y_1)}}[l,l';v'']:V(X_0X_1)\otimes V(Y_0Y_1)\rightarrow V(Z_0Z_1),\nonumber\\  
&&\hskip -1cm\Upsilon_{\;\;{}_{(Z_0Z_1)}}^{{}_{(X_0X_1)}{}_{(Y_0Y_1)}}[l,l';v'']= \sum_{IJK} {{\buildrel {\!\!{}_{(Z_0Z_1)}} \over {\Pi}}}(u^{K})\dagger \vert v''>  <l\otimes l'\vert{\buildrel {\!\!{}_{(X_0X_1)}} \over {\Pi}}  (u^I) \otimes {\buildrel {\!\!{}_{(Y_0Y_1)}} \over {\Pi}}(u^J)  
 \;h(u_Iu_Ju_K).  
\end{eqnarray}  
 
Its matrix elements satisfy the relation:  
\begin{eqnarray}  
&&\hskip -4.5cm \sum_{IJK}{\buildrel {\!\!{}_{(X_0X_1)}} \over {\Pi}}\!{}^{A',i'}_{A,i}\!(u^I){\buildrel {\!\!{}_{(Y_0Y_1)}} \over {\Pi}}\!{}^{B',j'}_{B,j}\!(u^J){\buildrel {\!\!{}_{(Z_0Z_1)}} \over {\Pi}}\!{}^{\dagger}{}^{C,k}_{C',k'}(u^K)\;\;h(u_Iu_Ju_K)=\nonumber 
\end{eqnarray}  
\vskip -0.8cm 
\begin{eqnarray} 
&&\hskip 3.5cm =\sum_{v}\Clebpsi{A}{B}{C}{i}{j}{k}{\buildrel C \over \mu}  
{}_{k'}^{v}\!\Clebphi{A'}{B'}{C'}{i'}{j'}{v}\!\!\Entrela{A'}{B'}{C'}{A}{B}{C}{X_0X_1}{Y_0Y_1}{Z_0Z_1}  
\nonumber 
\end{eqnarray}  
\vskip -0.4cm 
\begin{eqnarray} 
&&\hskip -1.5cm\mbox{with}\Entrela{A'}{B'}{C'}{A}{B}{C}{X_0X_1}{Y_0Y_1}{Z_0Z_1}\!\!=\!\!\!\sum_{KLMN}\frac{\qd{N}\qd{K}}{\qd{C}\qd{C'}}\Lambda^{KL}_{A'A}(X_0X_1)\Lambda^{KM}_{B'B}(Y_0Y_1)\Lambda^{KN}_{C'C}(Z_1Z_0)\times  
\nonumber\\  
&&\hskip 5cm\times\sixjn{B}{N}{L}{K}{A}{C}{0}\!\!\sixjn{K}{N}{C'}{A'}{B'}{M}{0}\!\!\sixjn{A'}{N}{M}{B}{K}{L}{0}\label{form3lambda}.  
\end{eqnarray}  
This infinite series  converges absolutely and uniformly in $\rho_X,\rho_Y,\rho_Z. $  
Its square is a continuous function of $\rho_X,\rho_Y,\rho_Z. $  
\end{theorem}  
  
\Proof 
 
We will use the convention of summation of repeated up and low small indices.  
\begin{eqnarray}  
&&\hskip -3.5cm \sum_{IJK}{\buildrel {\!\!{}_{(X_0X_1)}} \over {\Pi}}\!{}^{A',i'}_{A,i}\!(u^I){\buildrel {\!\!{}_{(X_0X_1)}} \over {\Pi}}\!{}^{B',j'}_{B,j}\!(u^J){\buildrel {\!\!{}_{(X_0X_1)}} \over {\Pi}}\!{}^{\dagger}{}^{C,k}_{C',k'}(u^K)\;h(u_Iu_Ju_K)\;=\nonumber 
\end{eqnarray}  
\vskip -0.4cm  
\begin{eqnarray}  
&&\hskip -0.5cm=\!\!\!\!\!\!\!\!\sum_{KLK'MK''N}\!\!\!\!\!\!\!\! \Lambda^{KL}_{A'A}(X_0X_1)  
\Clebphi{A'}{K}{L}{r}{s}{t}\!\Clebpsi{K}{A}{L}{u}{i}{t}  
\Lambda^{K'M}_{B'B}(Y_0Y_1)  
\Clebphi{B'}{K'}{M}{a}{b}{c}\!\Clebpsi{K'}{B}{M}{d}{j}{c}\!\!  
\times \nonumber\\  
&&\hskip -0.5cm\times   
\Lambda^{K''N}_{C'C}(Z_1Z_0) \Clebpsi{C'}{K''}{N}{m}{p}{n}\!\Clebphi{K''}{C}{N}{l}{k}{n}  
{\buildrel {C'} \over w}{}_{qk'}({\buildrel {C'} \over w}{}^{-1})^{mo}  
h(\kE{A'}{i'}{r}{K}{u}{s}\;  
 \kE{B'}{j'}{a}{K'}{d}{b} \;\kE{C'}{q}{o}{K''}{p}{l})\nonumber  
\end{eqnarray}  
\vskip -0.4cm  
\begin{eqnarray}  
&&\hskip -0.5cm=\!\!\!\sum_{KLMN} \!\!\!\frac{e^{2i \pi K}\qd{K}}{\qd{C'}}  
\Lambda^{KL}_{A'A}(X_0X_1)\Lambda^{KM}_{B'B}(Y_0Y_1)\Lambda^{KN}_{C'C}(Z_1Z_0)  
\Clebphi{A'}{K}{L}{r}{s}{t}\!\Clebpsi{K}{A}{L}{u}{i}{t}\!\!  
\Clebphi{B'}{K}{M}{a}{p}{c}\nonumber\\  
&&\hskip -0.5cm\times  
\Clebpsi{K}{B}{M}{s}{j}{c}\!\!  
 \Clebpsi{C'}{K}{N}{m}{p}{n}\!\Clebphi{K}{C}{N}{l}{k}{n}  
\Clebpsi{A'}{B'}{C'}{r}{a}{m} ({\buildrel {K} \over \mu}{}^{-1})^{u}_{l}  
\Clebphi{A'}{B'}{C'}{i'}{j'}{v}{\buildrel {C'} \over \mu}{}_{k'}^{v}\nonumber  
\end{eqnarray}  
after the use of formulae (\ref{haar}) and (\ref{normaclebsch})  
\begin{eqnarray}  
&&\hskip -0.5cm=\!\!\!\sum_{KLMN} \!\!\!\frac{e^{2i \pi K}\qd{K}}{\qd{C'}}  
\Lambda^{KL}_{A'A}(X_0X_1)\Lambda^{KM}_{B'B}(Y_0Y_1)\Lambda^{KN}_{C'C}(Z_1Z_0)  
\sixjn{K}{N}{C'}{A'}{B'}{M}{0}  
\Clebphi{A'}{K}{L}{r}{s}{t}\!  
\nonumber\\  
&&\hskip -0.5cm\times \Clebpsi{K}{A}{L}{u}{i}{t}\!\! \Clebpsi{A'}{M}{N}{r}{c}{n}  
\!\!\Clebpsi{K}{B}{M}{s}{j}{c}\!\!  
 \!\Clebphi{K}{C}{N}{l}{k}{n}  
 ({\buildrel {K} \over \mu}{}^{-1})^{u}_{l}  
\Clebphi{A'}{B'}{C'}{i'}{j'}{v}{\buildrel {C'} \over \mu}{}_{k'}^{v}\nonumber  
\end{eqnarray}  
\vskip -0.4cm  
\begin{eqnarray}  
&&\hskip -0.5cm=\!\!\!\sum_{KLMN} \!\!\!\frac{e^{2i \pi K}\qd{K}}{\qd{C'}}  
\Lambda^{KL}_{A'A}(X_0X_1)\Lambda^{KM}_{B'B}(Y_0Y_1)\Lambda^{KN}_{C'C}(Z_1Z_0)  
\sixjn{K}{N}{C'}{A'}{B'}{M}{0}  
\sixjn{A'}{N}{M}{B}{K}{L}{0}  
\nonumber\\  
&&\hskip -0.5cm\times \Clebpsi{K}{A}{L}{u}{i}{t}\!\!   
\Clebpsi{L}{B}{N}{t}{j}{n}  
 \!\Clebphi{K}{C}{N}{l}{k}{n}  
 ({\buildrel {K} \over \mu}{}^{-1})^{u}_{l}  
\Clebphi{A'}{B'}{C'}{i'}{j'}{v}{\buildrel {C'} \over \mu}{}_{k'}^{v}\nonumber  
\end{eqnarray}  
after having used twice the formula (\ref{vertexIRF})  
\begin{eqnarray}  
&&\hskip -0.5cm=\!\!\!\sum_{KLMN} \!\!\!\frac{e^{2i \pi K}\qd{K}}{\qd{C'}}  
\Lambda^{KL}_{A'A}(X_0X_1)\Lambda^{KM}_{B'B}(Y_0Y_1)\Lambda^{KN}_{C'C}(Z_1Z_0)  
\sixjn{K}{N}{C'}{A'}{B'}{M}{0}  
\sixjn{A'}{N}{M}{B}{K}{L}{0}  
\nonumber\\  
&&\hskip -0.5cm\times \frac{e^{2i \pi (N-C)}\qd{N}}{\qd{C}}  
\Clebphi{K}{L}{A}{l}{t}{i}\Clebpsi{L}{B}{N}{t}{j}{n}\!\Clebpsi{K}{N}{C}{l}{n}{k}  
 ({\buildrel {K} \over \mu}{}^{-1})^{u}_{l}  
\Clebphi{A'}{B'}{C'}{i'}{j'}{v}{\buildrel {C'} \over \mu}{}_{k'}^{v}\nonumber\\  
&&\hskip -0.5cm= \mbox{righthandside of (\ref{form3lambda})},\nonumber  
\end{eqnarray}  
where we have used formulae (\ref{contra}) and (\ref{vertexIRF}) to conclude. 
  
This proof holds true as soon as we have shown that  these  series are absolutely convergent.   
Proof of the convergence of the last series  can  be obtained using asymptotic properties of coefficients $\Lambda^{BC}_{AD},$ as well as asymptotic properties of $6j(0)$. Indeed, using the selection rules of $6j(0),$ we know that in the   four sums  
\begin{eqnarray}  
&&\sum_{Klmn}\frac{\qd{K+n}\qd{K}}{\qd{C}\qd{C'}}\Lambda^{K\;K+l}_{A'A}(X_0X_1)\Lambda^{K\;K+m}_{B'B}(Y_0Y_1)\Lambda^{K\;K+n}_{C'C}(Z_1Z_0)\times  
\nonumber\\  
&&\hskip 1cm\times\sixjn{B}{K+n}{K+l}{K}{A}{C}{0}\!\!\sixjn{K}{K+n}{C'}{A'}{B'}{K+m}{0}\!\!\sixjn{A'}{K+n}{K+m}{B}{K}{K+l}{0},  
\end{eqnarray}  
for a fixed $K,$ the range of the sum over $l,m,n$ is finite and fixed by $\vert l \vert \leq {\rm min}(A,A'),\vert m \vert \leq {\rm min}(B,B'), \vert n \vert \leq {\rm min}(C,C'),\vert l-n \vert \leq B, \vert n-m \vert \leq A'.$ Thus, in order to show the absolute convergence  of this series, it is sufficient to bound  the general term in $K$ by a geometric series.  
This can be proved using the behaviour of the coefficients $\Lambda^{KK+l}_{AD}$ when $K$ is large, derived in the appendix.  
 Precisely, using (\ref{asymptlambda1}) and (\ref{asymptot1})(\ref{asymptot2}), we have the following property: $\forall l,m,n,A,B,C,A',B',C' \in \onehalf \ZZ,$ there exist $Q,Q' \in \RR,$ such that $\forall \rho_X, \rho_Y, \rho_Z \in ]-\frac{\pi}{\hbar},\frac{\pi}{\hbar}]$ : 
\begin{eqnarray}  
&&\hskip -1cm\forall K,\left\vert  \frac{\qd{K+n}\qd{K}}{\qd{C}\qd{C'}}\Lambda^{K\;K+l}_{A'A}(X_0X_1)\Lambda^{K\;K+m}_{B'B}(Y_0Y_1)\Lambda^{K\;K+n}_{C'C}(Z_1Z_0) \right\vert \; \leq \; Q K^3 q^{2K}  
\nonumber\\  
&&\hskip -1cm\forall K,\left\vert  \sixjn{B}{K+n}{K+l}{K}{A}{C}{0}\!\!\sixjn{K}{K+n}{C'}{A'}{B'}{K+m}{0}\!\!\sixjn{A'}{K+n}{K+m}{B}{K}{K+l}{0} \right\vert \; \leq \; Q'.  
\end{eqnarray}  
The series is  therefore absolutely convergent and  uniformly convergent in   
$\rho_X ,\rho_Y,\rho_Z.$  
From this last result and using the property that 
 ${\cal N}^{(A)}(q^{2X_1+1},m_X){\cal N}^{(D)}(q^{2X_1+1},m_X)^{-1}\Lambda^{BC}_{AD}(X_0X_1)$ is a  continuous functions of $\rho_X,$  
 the square of the  matrix coefficients of the operator $\Upsilon_{\;\;{}_{(Z_0Z_1)}}^{{}_{(X_0X_1)}{}_{(Y_0Y_1)}}$ are continuous functions of the three variables $\rho_X, \rho_Y, \rho_Z.$  
Note that in the case where $A=A', B=B', C=C'$ the expression (\ref{form3lambda}) is itself a continuous function of $\rho_X, \rho_Y, \rho_Z.$ 
 
\medskip 
 
We still have to show that the linear map $\Upsilon_{\;\;{}_{(Z_0Z_1)}}^{{}_{(X_0X_1)}{}_{(Y_0Y_1)}}[l,l';v'']$ is an intertwiner.  
We have   
\begin{eqnarray}  
&&\Upsilon_{\;\;{}_{(Z_0Z_1)}}^{{}_{(X_0X_1)}{}_{(Y_0Y_1)}}[l,l';v'']=\nonumber\\  
&& =\sum_{IJK} {{\buildrel {\!\!{}_{(Z_0Z_1)}} \over {\Pi}}}(u^{K})\dagger \vert v''>  <l\otimes l'\vert{\buildrel {\!\!{}_{(X_0X_1)}} \over {\Pi}}  (u^I) \otimes {\buildrel {\!\!{}_{(Y_0Y_1)}} \over {\Pi}}(u^J)  
 \;h(u_Iu_Ju_K)\nonumber\\  
&&= \sum_{IJK} {{\buildrel {\!\!{}_{(Z_0Z_1)}} \over {\Pi}}}(S^{-1}(u^{K})) \vert v''> <l\otimes l'\vert{\buildrel {\!\!{}_{(X_0X_1)}} \over {\Pi}}  (u^I) \otimes {\buildrel {\!\!{}_{(Y_0Y_1)}} \over {\Pi}}(u^J)\;h(u_Iu_Ju_K)\nonumber  
\end{eqnarray}  
Using the identity  verified by any Haar measure on a co-quasitriangular Hopf algebra $A$ \cite{Wo1,C}  
\begin{eqnarray}  
h(ab)=\sum_{(a)}h(ba_{(2)})<\mu,a_{(1)}><\mu,a_{(3)}>, \;\;\;\forall a,b \in A,   
\end{eqnarray}  
and simple transformations using properties of the antipode and of the group-like element $\mu$, we obtain that $\Upsilon_{\;\;{}_{(Z_0Z_1)}}^{{}_{(X_0X_1)}{}_{(Y_0Y_1)}}[l,l';v'']$ is related to the operator:  
\begin{eqnarray}  
&& \Xi_{\;\;{}_{(Z_0Z_1)}}^{{}_{(X_0X_1)}{}_{(Y_0Y_1)}}[l,l';v'']:V(X_0X_1)\otimes V(Y_0Y_1)\rightarrow V(Z_0Z_1),\\  
&&\Xi_{\;\;{}_{(Z_0Z_1)}}^{{}_{(X_0X_1)}{}_{(Y_0Y_1)}}[l,l';v'']= \sum_{IJK} {{\buildrel {\!\!{}_{(Z_0Z_1)}} \over {\Pi}}}(S(u^{K})\mu) \vert v''>  <l\otimes l'\vert{\buildrel {\!\!{}_{(X_0X_1)}} \over {\Pi}}  (u^I) \otimes {\buildrel {\!\!{}_{(Y_0Y_1)}} \over {\Pi}}(u^J)  
 \;h(u_Ku_Iu_J),\nonumber  
\end{eqnarray}  
as follows:  
  
$\Upsilon_{\;\;{}_{(Z_0Z_1)}}^{{}_{(X_0X_1)}{}_{(Y_0Y_1)}}[l,l';v'']=  
\sum_{C,r}\Xi_{\;\;{}_{(Z_0Z_1)}}^{{}_{(X_0X_1)}{}_{(Y_0Y_1)}}[l,l';{\buildrel {C} \over e}_r]  
<{\buildrel {C} \over \mu}{\buildrel {C} \over e}_r, v''>.  
$ 
 
As a result it is equivalent to show that $\Xi_{\;\;{}_{(Z_0Z_1)}}^{{}_{(X_0X_1)}{}_{(Y_0Y_1)}}[l,l';v'']$ is an intertwiner operator.  
This is a simple consequence of the absolute convergence of the series and of the right invariance of the Haar measure.  
  \eoProof  

 The following proposition gives the link between intertwiners constructed using the Haar measure and the reduced elements.
 
\begin{proposition}  
There exists a real number ${\cal M}_{\;\;{}_{(Z_0Z_1)}}^{{}_{(X_0X_1)}{}_{(Y_0Y_1)}}$ such that  
\begin{eqnarray}  
&&\hskip -0.6cm\Entrela{A'}{B'}{C'}{A}{B}{C}{X_0X_1}{Y_0Y_1}{Z_0Z_1}=   
{\cal M}_{\;\;{}_{(Z_0Z_1)}}^{{}_{(X_0X_1)}{}_{(Y_0Y_1)}}  
\ElemRed{A}{B}{C}{X_0X_1}{Y_0Y_1}{Z_0Z_1} \; \overline{\ElemRed{A'}{B'}{C'}{X_0X_1}{Y_0Y_1}{Z_0Z_1}}.  
\label{mainformula}  
\end{eqnarray}  
\end{proposition}  
  
\Proof  
 
The operator $\Xi_{\;\;{}_{(Z_0Z_1)}}^{{}_{(X_0X_1)}{}_{(Y_0Y_1)}}[l,l';v'']$ being an intertwiner,the theorem 1 ensures the existence of  a complex number $\lambda (l,l';v'')$ such that:  
 $\Xi_{\;\;{}_{(Z_0Z_1)}}^{{}_{(X_0X_1)}{}_{(Y_0Y_1)}}[l,l';v'']=\lambda (l,l';v'')  
\GPhin{X_0X_1}{Y_0Y_1}{Z_0Z_1}.  
$  
We will denote $\Pi={\buildrel {\!\!{}_{(X_0X_1)}} \over {\Pi}}, \Pi'={\buildrel {\!\!{}_{(Y_0Y_1)}} \over {\Pi}},\Pi''={\buildrel {\!\!{}_{(Z_0Z_1)}} \over {\Pi}}$  
and we have  the following sequence of equalities:  
\begin{eqnarray}  
&&\hskip -0.6cm \overline{\lambda(l, l';v'')}\;\;\overline{<l''\vert \Phi^{\Pi \;\Pi'}_{\;\;\Pi''} \vert v \otimes v'>}=\nonumber\\  
&&\hskip -0.6cm=\sum_{IJK}<v''\vert\Pi''((S(e^I)\mu)^*)\vert l''><v\otimes v' \vert \Pi(e^J{}^*) \otimes \Pi'(e^K{}^*)\vert l\otimes l'> \overline{h(e_Ie_Je_K)}=\nonumber\\  
&&\hskip -0.6cm=\sum_{IJK}<v''\vert\Pi''(\mu S^{-1}(e^I{}^*))\vert l''><v\otimes v' \vert \Pi(e^J{}^*) \otimes \Pi'(e^K{}^*)\vert l\otimes l'> h(e^*_Ke^*_Je^*_I)=\nonumber\\  
&&\hskip -0.6cm=\sum_{IJK}<v''\vert\Pi''(\mu S^{-1}(e^I))\vert l''><v\otimes v' \vert \Pi(e^J) \otimes \Pi'(e^K)\vert l\otimes l'> h(S^{-1}(e_K)S^{-1}(e_J)S^{-1}(e_I))=\nonumber\\  
&&\hskip -0.6cm=\sum_{IJK}<v''\vert\Pi''(S(e^I)\mu)\vert l''><v\otimes v' \vert \Pi(e^J) \otimes \Pi'(e^K)\vert l\otimes l'> h(e_Ie_Je_K)=\nonumber\\  
&&\hskip -0.6cm= \lambda(v, v';l'')\;<v'' \vert \Phi^{\Pi \;\Pi'}_{\;\;\Pi''} \vert l \otimes l'>.\nonumber  
\end{eqnarray}  
As a result   
\begin{eqnarray}  
\lambda(v, v';l'')= {\cal M}_{\;\;\Pi''}^{\Pi \;\Pi'} \; {\overline { <l''\vert \Phi^{\Pi \;\Pi'}_{\;\;\Pi''} \vert v \otimes v'>}},\nonumber  
\end{eqnarray}  
from which we get  
\begin{eqnarray}  
&&\hskip -1cm\sum_{IJK} <l''\vert  {\Pi''}(S(e^{I})\mu) \vert v''> <l\vert \Pi(e^J)\vert v>  
 <l'\vert \Pi'(e^K)\vert v'> \;h(e_Ie_Je_K)=\nonumber\\  
&&\hskip 1cm= {\cal M}_{\;\;\Pi''}^{\Pi \;\Pi'}  
 <l''\vert \Phi^{\Pi \;\Pi'}_{\;\;\Pi''} \vert v \otimes v'> \overline{<v'' \vert \Phi^{\Pi \;\Pi'}_{\;\;\Pi''} \vert l \otimes l'>}.  
\end{eqnarray}  
Using the relation between the operator $\Xi_{\;\;{}_{(Z_0Z_1)}}^{{}_{(X_0X_1)}{}_{(Y_0Y_1)}}[l,l';v'']$ and $\Upsilon_{\;\;{}_{(Z_0Z_1)}}^{{}_{(X_0X_1)}{}_{(Y_0Y_1)}}[l,l';v'']$ we obtain the announced result.  \eoProof  
  
In order to make a precise connection, later on, between matrix elements of   
${\hat \Phi}[X_0X_1,Y_0Y_1]$ and Askey-Wilson polynomials we need an explicit formula for  ${\cal M}.$  
  
\begin{proposition}  
The normalization factor ${\cal M}_{\;\;{}_{(Z_0Z_1)}}^{{}_{(X_0X_1)}{}_{(Y_0Y_1)}}$ is a positive number and is given by the following expression:  
\begin{eqnarray}  
&&\hskip -1cm{\cal M}_{\;\;{}_{(Z_0Z_1)}}^{{}_{(X_0X_1)}{}_{(Y_0Y_1)}}=  
\vert\frac{ q^{2(X_0-X_1+Y_0-Y_1+Z_0-Z_1)}(1-q^2)(1)^2_{\infty}  
}{\nu_1(2X_0\!+\!1)^2\nu_1(2Y_0\!+\!1)^2\nu_1(2Z_0\!+\!1)^2}\vert \times\nonumber\\  
&&\hskip 0.8cm\times\frac{\vert \xi(2X_1\!+\!1)\xi(2Y_1\!+\!1)\xi(2Z_1\!+\!1)\vert   
}{\vert \xi(X_1\!+\!Y_1\!-\!Z_1\!+\!1)\xi(X_1\!-\!Y_1\!+\!Z_1\!+\!1)\xi(Y_1\!+\!Z_1\!-\!X_1+\!1)\xi(X_1\!+\!Y_1\!+\!Z_1+\!2)\vert }  \label{Mformula} 
\end{eqnarray}  
where we have defined the function $\xi$, by $\xi(z)=(z)_{\infty}(1-z)_{\infty}.$  
\end{proposition}  
  
\Proof  
 
We will identify certain asymptotics of the lefthandside and the righthandside of (\ref{mainformula}) in order to compute the normalization ${\cal M}_{\;\;\;\;{}_{(Z_0Z_1)}}^{{}_{(X_0X_1)}{}_{(Y_0Y_1)}}.$ Let us first find the behaviour of   
\begin{eqnarray}  
&&\hskip -1cm f(T)=\Entrela{P+T}{P}{T}{P+T}{P}{T}{X_0X_1}{Y_0Y_1}{Z_0Z_1}\;\;\;\; \mbox{when}\;\;T \rightarrow +\infty, \;\;P=Y_0-Y_1 >0.  
\end{eqnarray}  
 
It can be seen from the properties (\ref{asymptot1})(\ref{asymptot2})(\ref{asymptlambda2}),  
  that the unique leading term in the explicit expression : 
\begin{eqnarray}  
&&\hskip -1cm f(T)=  
\sum_{KlMn}\frac{\qd{T+n}\qd{K}}{\qd{T}^2}\Lambda^{K\;P+T+l}_{P+T}(X_0X_1)  
\Lambda^{K\;M}_{\;\;P}(Y_0Y_1)\Lambda^{K\;T+n}_{\;\;T}(Z_1Z_0)\!\times  
\nonumber\\  
&&\times\sixjn{P}{T\!+\!n}{P\!+\!T\!+\!l}{K}{P\!+\!T}{T}{0}\!\!\sixjn{P\!+\!T}{T\!+\!n}{M}{K}{P}{T}{0}  
\!\!\sixjn{P\!+\!T}{T\!+\!n}{M}{P}{K}{P\!+\!T\!+\!l}{0},\nonumber  
\end{eqnarray}  
 is  the term corresponding to $l=n=0.$  
 The  behaviour when $T \rightarrow +\infty$ of $f(T)$ is therefore given by  
\begin{eqnarray}  
f(T) \sim q^{2T} (1-q^2)\sum_{KM} \qd{K}\Clebpsi{K}{\;P}{\;M}{0}{-P}{-P}\Clebphi{\;P}{K}{\;M}{-P}{0}{-P} \Lambda^{K\;M}_{\;\;P}(Y_0Y_1).\nonumber  
\end{eqnarray}  
The series in the right handside can be exactly computed  using the formula (\ref{formsumlambda}) proved in the appendix. We finally obtain, for $P=Y_0-Y_1 >0:$  
\begin{eqnarray}  
\Entrela{P+T}{P}{T}{P+T}{P}{T}{X_0X_1}{Y_0Y_1}{Z_0Z_1} \sim_{{}_{T \rightarrow +\infty}} q^{2T} (1-q^2) \frac{\hbar \qd{\vert Y_0-Y_1\vert }}{2\pi {\cal P}(Y_0Y_1)}.  
\end{eqnarray}  
From the fact that the righthandside of the last expression is positive, and that the right handside of (\ref{mainformula}) is in this situation a modulus square we obtain that  
 ${\cal M}_{\;\;{}_{(Z_0Z_1)}}^{{}_{(X_0X_1)}{}_{(Y_0Y_1)}}$ is a positive number.  
Let us now compute the behaviour of the righthandside of (\ref{mainformula})  
for $T \rightarrow +\infty, P=m_Y$ (this computation has been  done for $0 \leq m_Y,$ but for  
$m_Y < 0$ a similar proof can be done). In the appendix, we have proved the following result:   
\begin{eqnarray}  
&&\hskip -0.8cm\ElemRed{$\!\!\!T\!+\!Y_0\!\!-\!\!Y_1$}{$\!\!Y_0\!\!-\!\!Y_1$}{T}{X_0X_1}{Y_0Y_1}{Z_0Z_1}  
 {\sim}\;  
e^{-i \pi 2T}q^{-T(Y_0+Y_1)}\;e^{i \pi (\frac{X_0-X_1}{2})}   
q^{-Z_0Z_1}\sqrt{\qd{Y_0-Y_1}}\!\!\times\nonumber\\ 
&&\hskip -0.5cm\times \nu_1(2X_1\!+\!1)\nu_1(2Z_0\!+\!1)   
\frac{\varphi_{(2X_0,P\!+\!T\!+\!X_0\!-\!X_1)}\varphi_{(Z_1\!+\!X_1\!-\!Y_1,Z_1\!-\!Z_0\!+\!X_1\!-\!X_0\!+\!Y_0\!-\!Y_1)}}{\varphi_{(Z_0\!+\!Z_1,T)}}  
\times\nonumber\\  
&&\hskip -0.5cm\times \;  
\frac{\nu_{\infty}(Y_0\!+\!X_0\!-\!Z_0\!+\!1)  
\nu_{\infty}(Y_0\!+\!X_0\!+\!Z_0\!+\!2)\nu_{\infty}(Y_0\!-\!X_0\!-\!Z_0)  
\nu_{\infty}(Y_0\!+\!Z_0\!-\!X_0\!+\!1)}{  
q^{\frac{3}{2}X_0^2+\onehalf X_0-\onehalf Z_0^2+\onehalf Z_0}  
(1)_{\infty}  
\nu_{\infty}(2Y_1\!+\!1)\nu_{\infty}(2Y_0\!+\!1)}\times\nonumber\\  
&&\hskip -0.5cm\times \;  
\frac{\nu_{\infty}(Y_1\!+\!X_1\!-\!Z_1\!+\!1)\nu_{\infty}(Y_1\!+\!X_1\!+\!Z_1\!+\!2)  
\nu_{\infty}(Y_1\!-\!X_1\!-\!Z_1)\nu_{\infty}(Y_1\!+\!Z_1\!-\!X_1\!+\!1)}{  
q^{-\frac{3}{2}X_1^2-\onehalf X_1-\onehalf Z_1^2 -\onehalf Z_1}  
\nu_{\infty}(2X_0\!+\!1)  
\nu_{\infty}(-\!2X_0)\nu_{\infty}(Z_0\!+\!Z_1+1)\nu_{\infty}(-\!Z_0\!-\!Z_1)}.  
\end{eqnarray}  
From  (\ref{mainformula}), and after elementary algebraic relations on q-factorials, we finally get the announced expression for ${\cal M}_{\;\;{}_{(Z_0Z_1)}}^{{}_{(X_0X_1)}{}_{(Y_0Y_1)}}$.   
\eoProof  
 
\medskip 
  
The formula $(\ref{mainformula})$ is an identity where the left hansdside is a complicated series, and the right handside is the  infinite product (${\cal M}$), multiplied by square roots of rational fractions.  
Let us give an example of this formula in its simplest case, which is  achieved when $A=A'=B=B'=C=C'=0.$ We necessarily have  $m_X=m_Y=m_Z=0,$ i.e  $X_0=X_1, Y_0=Y_1, Z_0=Z_1.$  
In this case we can compute  ${\cal M}_{\;\;{}_{(Z_0Z_1)}}^{{}_{(X_0X_1)}{}_{(Y_0Y_1)}}$ directly. Indeed:  
\begin{equation}  
\Entrela{0}{0}{0}{0}{0}{0}{X_0X_1}{Y_0Y_1}{Z_0Z_1}=  
\sum_{B}[d_B]^2 \Lambda^{BB}_{0}(X_0X_1) \Lambda^{BB}_{0}(Y_0Y_1)\Lambda^{BB}_{0}(Z_0Z_1)\label{Sommem0}. 
\end{equation}  
Using the identity   
$\Lambda^{BB}_{0}(X_0X_0)=\frac{[(2X_0+1)(2B+1)]}{[2X_0+1][2B+1]},$  
 the formula (Exercice 12 of \cite{WW}) 
\begin{equation}  
4\sum_{n=1}^{+\infty}q^{n}\frac{sin(2nz)}{1-q^{2n}}=  
\frac{\theta'_{4}(z,\tau)}{\theta_4(z,\tau)}  
\end{equation}  
where $q=e^{-\hbar}=e^{i\pi\tau}$, 
and the Jacobi's fundamental formulae (\cite{WW} 21.22):  
\begin{equation} 
2\theta_4(x')\theta_4(y')\theta_4(z')\theta_4(w')= 
\sum_{j=1}^4(-1)^{\delta_{j,2}}\theta_j(x)\theta_j(y)\theta_j(z)\theta_j(w) 
\end{equation} 
where $2w'=-w+x+y+z$ (and circular permutation to define x',y',z'), 
 we obtain  that the series (\ref{Sommem0}) is equal to:

\begin{equation}  
\frac{(1-q^2)q^{-3/4}(1)_{\infty}^3\theta_1(2\lambda_x)\theta_1(2\lambda_y)\theta_1(2\lambda_z)(sin(2\lambda_x)sin(2\lambda_y)sin(2\lambda_z))^{-1}}  
{8\theta_4(\lambda_x+\lambda_y+\lambda_z)\theta_4(-\lambda_x+\lambda_y+\lambda_z)\theta_4(\lambda_x-\lambda_y+\lambda_z)  
\theta_4(\lambda_x+\lambda_y-\lambda_z)}\label{Mtheta}  
\end{equation}  
where we have chosen $2X_0+1=i\frac{4\lambda_x}{\hbar}, 2Y_0+1=i\frac{4\lambda_y}{\hbar},   
2Z_0+1=i\frac{4\lambda_z}{\hbar},$ and $\theta_j, (j=1,..4)$ are the usual Theta functions, as defined in \cite{WW}.  
Note that this expression is exactly the coefficient ${\cal M}_{\;\;{}_{(Z_0Z_0)}}^{{}_{(X_0X_0)}{}_{(Y_0Y_0)}},$ once we have expressed the   
Theta functions in term of infinite products.  
However such a simple derivation of ${\cal M}_{\;\;{}_{(Z_0Z_1)}}^{{}_{(X_0X_1)}{}_{(Y_0Y_1)}},$ does not generalize so easily to the case where $m_X, m_Y, m_Z$ are non zero. \\ 
Formula (\ref{Mtheta}) is the quantum deformation of a formula first   
found in \cite{ARRW}.  
 
\bigskip 
  
Remark: We will need in the next section the following identity, which can be easily proved using $(\ref{Mformula})$. ${\cal M}_{\;\;{}_{(Z_0Z_1)}}^{{}_{(X_0X_1)}{}_{(Y_0Y_1)}},$ can also be written as: 
\begin{eqnarray} 
&&\hskip -0.6cm {\cal M}_{\;\;{}_{(Z_0Z_1)}}^{{}_{(X_0X_1)}{}_{(Y_0Y_1)}}= 
\psi_1\times 
 \frac{q^{Z_0-3 Z_1+Y_0-3 Y_1+X_0-3 X_1-3}(1)_{\infty}^2(1-q^2)}{\nu_1(2X_0\!+\!1)\nu_1(-2X_1\!-\!1)\nu_1(2Y_0\!+\!1) 
\nu_1(-2Y_1\!-\!1)\nu_1(2Z_0\!+\!1)\nu_1(-2Z_1\!-\!1)}\times\nonumber\\ 
&&\times \frac{\xi(2X_1\!+\!1,2Y_1\!+\!1,2Z_1\!+\!1)}{\xi(Y_1\!+\!X_1\!-\!Z_1\!+\!1,X_1\!+\!Z_1\!-\!Y_1\!+\!1,Y_1\!+\!Z_1\!-\!X_1\!+\!1,X_1\!+\!Y_1\!+\!Z_1\!+\!2)}\label{Mformula2} 
\end{eqnarray} 
where we have denoted $\xi(a_1,\cdots, a_n)=\prod_{k=1}^{n}\xi(a_k)$ and where $\psi_1$ is the fourth root of unit given by: 
\begin{eqnarray} 
&&\hskip -0.6cm \psi_1=\frac{\varphi_{(2Y_1,2Y_1-2Y_0-1)}\times (\mbox {\small cycl. perm. on}\;\; X,Y,Z))}{\varphi_{(Y_1\!+\!X_1\!-\!Z_1\!+\!1, Y_1\!-\!Y_0\!+\!X_1\!-\!X_0\!-\!Z_0\!+\!Z_1)}\times (\mbox{\small cycl. perm. on}\;\; X,Y,Z)\times \varphi_{(Y_1\!+\!X_1\!+\!Z_1\!+\!2,Y_1\!-\!Y_0\!+\!X_1\!-\!X_0\!+\!Z_1\!-\!Z_0)}}\nonumber 
\end{eqnarray}

\section{Decomposition Theorem of the Tensor Product of Principal Representations and Askey-Wilson Polynomials}  
\subsection{Decomposition Theorem of the Tensor Product of Principal Representations} 
Let us now prove that, chosing $N^{{}_{(X_0X_1)\;(Y_0Y_1)}}_{\;\;\;{}_{(Z_0Z_1)}}\!\!=({\cal M}_{\;\;{}_{(Z_0Z_1)}}^{{}_{(X_0X_1)}{}_{(Y_0Y_1)}})^{1/2}\!>0$, ${\hat \Phi}[X_0X_1,Y_0Y_1]$ is an isometry.

\begin{theorem}  \label{isometrytheorem} 
${\hat \Phi}[X_0X_1,Y_0Y_1]: V(X_0X_1)\otimes V(Y_0Y_1) \rightarrow   
\int^{\oplus}d(Z_0Z_1)H(Z_0Z_1)$ is an isometry.  
This is a direct consequence of the following identity:  
\begin{equation}\forall C\in \onehalf \ZZ^{+},   
\int d{\cal P}(Z_0Z_1)\Entrela{A'}{B'}{C}{A}{B}{C}{X_0X_1}{Y_0Y_1}{Z_0Z_1}=\delta_{A,A'}\delta_{B,B'}Y^{(1)}_{(A,X_0-X_1)}Y^{(1)}_{(B,Y_0-Y_1)}.  
\label{completeness}  
\end{equation}  
\end{theorem}  
  
\Proof  
 
$ V(X_0X_1)\otimes V(Y_0Y_1)$ is a ${\mathfrak U}_q(su(2))$ module, let us denote by   
$( V(X_0X_1)\otimes V(Y_0Y_1))^{[C]}$ its isotypic component of spin $C$.  
From the definition of $ {\hat \Phi}[X_0X_1,Y_0Y_1]$ and the choice $N=({\cal M}_{\;\;{}_{(Z_0Z_1)}}^{{}_{(X_0X_1)}{}_{(Y_0Y_1)}})^{1/2},$ it is simple to show that the restriction of $ {\hat \Phi}[X_0X_1,Y_0Y_1]$ to $( V(X_0X_1)\otimes V(Y_0Y_1))^{[C]}$ is an isometry if and  only if $(\ref{completeness})$ is satisfied for this $C.$  
We  already  proved that $\Entrela{A'}{B'}{C}{A}{B}{C}{X_0X_1}{Y_0Y_1}{Z_0Z_1}$ is a piecewise  continuous function of $\rho_Z$, as a result it is integrable, and   
using formula (\ref{form3lambda}) and  the uniform convergence properties of this series, we can invert the integral and the sum:  
\begin{eqnarray}  
&&\hskip -2cm\mbox{lefthandside of} \;\;(\ref{completeness})=\nonumber\\ 
&&\hskip -2cm=\!\!\!\sum_{KLMN}\!\!\frac{\qd{N}\qd{K}}{\qd{C}^2}\Lambda^{KL}_{A'A}(X_0X_1)\Lambda^{KM}_{B'B}(Y_0Y_1)\int d{\cal P}(Z_0Z_1)\Lambda^{KN}_{A'A}(Z_1Z_0)\times\nonumber\\  
&&\hskip 2cm\times\sixjn{B}{N}{L}{K}{A}{C}{0}\!\!\sixjn{K}{N}{C}{A'}{B'}{M}{0}\!\!\sixjn{A'}{N}{M}{B}{K}{L}{0} \nonumber 
\end{eqnarray}  
Using the  Plancherel formula (\ref{intlambda}), we obtain that this series is equal to:  
\begin{eqnarray}  
&&\hskip -1cm=\delta_{A,A'}\delta_{B,B'}\Lambda^{0A}_{AA}(X_0X_1)\Lambda^{0B}_{BB}(Y_0Y_1)\sixjn{B}{C}{A}{0}{A}{C}{0}\!\!\sixjn{0}{C}{C}{A}{B}{B}{0}\!\!\sixjn{A}{C}{B}{B}{0}{A}{0}=\nonumber\\  
&&\hskip -1cm=\mbox{righthandside of} \;\;(\ref{completeness})\nonumber.  
\end{eqnarray}  
 \eoProof  
 
\medskip  
 
 \noindent ${\hat \Phi}[X_0X_1,Y_0Y_1]$ being an isometry from $ V(X_0X_1)\otimes V(Y_0Y_1)$ to the Hilbert space $\int^{\oplus}d(Z_0Z_1)H(Z_0Z_1)$, there exists a  unique isometry   
$ \Phi^{\#}[X_0X_1,Y_0Y_1]$ from the tensor product of Hilbert space  $ H(X_0X_1)\otimes H(Y_0Y_1)$ to   
 $\int^{\oplus}d(Z_0Z_1)H(Z_0Z_1)$ which restriction to $ V(X_0X_1)\otimes V(Y_0Y_1)$ gives back ${\hat \Phi}[X_0X_1,Y_0Y_1].$

\begin{theorem}  
$\Phi^{\#}[X_0X_1,Y_0Y_1]:  H(X_0X_1)\otimes H(Y_0Y_1) \rightarrow   
\int^{\oplus}d(Z_0Z_1)H(Z_0Z_1)$ is an invertible isometry of Hilbert space.  
\end{theorem}  
  
\Proof  
 
It is trivial that  $\Phi^{\#}[X_0X_1,Y_0Y_1]$ is an isometry, as a result it  
is injective. We therefore have only to show that   
 $\Phi^{\#}[X_0X_1,Y_0Y_1]$  
is surjective.   
The proof goes along the same lines as the proof of the surjectivity of the  
Fourier    
 transform in \cite{BR1}.    
It is sufficient to show that  
$ E=\Phi^{\#}[X_0X_1,Y_0Y_1](V(X_0X_1)\otimes V(Y_0Y_1))$ is dense in 
$\int^{\oplus}d(Z_0Z_1)H(Z_0Z_1),$ which is equivalent to show that  
$E^{\perp}=\{0\}.$  
Let $f\in \int^{\oplus}d(Z_0Z_1)H(Z_0Z_1),$ and assume that  
 $<f, \Phi[X_0X_1,Y_0Y_1](\ontop{A}{e}_i(X_0X_1)\otimes  
\ontop{B}{e}_j(Y_0Y_1))>=0, \forall A,B,i,j.$ 
This in particular implies that  
\begin{equation} 
<f, \Phi[X_0X_1,Y_0Y_1](\Delta(a)(\ontop{A}{e}_i(X_0X_1)\otimes 
\ontop{B}{e}_j(Y_0Y_1)))>=0 
\end{equation} 
for any element $a$ in the center of ${\mathfrak U}_q(sl(2,\CC))_{\RR}.$ 
As a result, by taking $a=\Omega_+^{p}\Omega_-^{r}$ and using the property that  
$\GPhin{X_0X_1}{Y_0Y_1}{Z_0Z_1}$ is an intertwiner, we get the identity: 
$\forall p,r\in \NN, \forall A,B,C\in \onehalf \ZZ^{+}, \forall k=-C,..,C$ 
\begin{equation} 
\int d(Z_0Z_1) <f\vert \ontop{C}{e}_k(Z_0Z_1)> 
\omega_+(Z_0Z_1)^p\omega_-(Z_0Z_1)^r 
 g(X_0X_1,Y_0,Y_1,Z_0,Z_1)=0,\label{centerorthogonality} 
\end{equation} 
where $g(X_0X_1,Y_0,Y_1,Z_0,Z_1)= 
\ElemRed{A}{B}{C}{X_0X_1}{Y_0Y_1}{Z_0Z_1} N^{{}_{(X_0X_1)\;(Y_0Y_1)}}_{\;\;\;{}_{(Z_0Z_1)}}.$ 
$\vert g\vert^2$ being a  continuous  function in $\rho_Z,$   
$<f\vert \ontop{C}{e}_k> g$ is an $L^2$ function in the variable $\rho_Z.$ 
We can now apply  the argument of the proof of Plancherel Formula in \cite{BR1}and conclude that  the set of equations (\ref{centerorthogonality}) for all $p,r$ implies that $<f\vert \ontop{C}{e}_k> g=0$ as a $L^2$ function. From theorem 1, we can always find for every $C$, a couple $A,B$ such that $g\not=0$. From the explicit form of $g$ we know that the zeroes of $g$ considered as a $\rho_Z$ function are finite, we obtain that  $<f\vert \ontop{C}{e}_k>=0,$  which implies  $f=0.$ 
 \eoProof

\subsection{R-matrix in the Tensor Product of Infinite dimensional Representations} 
In this section we will compute the $R$-matrix in the tensor product of two irreducible infinite  
dimensional 
representations $ ({\buildrel {_{(X_0X_1)}}\over \Pi}\otimes  
 {\buildrel {_{(Y_0Y_1)}}\over \Pi})$. We will in particular show that the 
expression of $({\buildrel {_{(X_0X_1)}}\over \Pi}\otimes  
 {\buildrel {_{(Y_0 Y_1)}}\over \Pi})({\cal R})$ is an operator which domain 
contains all the vectors of the form  
${\buildrel A\over e}_m(X_0X_1)\otimes {\buildrel B\over e}_n(Y_0Y_1)$, and more precisely we have the following proposition: 
\begin{proposition} 
The expression of the $R$ matrix of  $\mathfrak{U}_{q}(sl(2,\mathbb{C})_{\RR})$ represented on  
${\buildrel {_{(X_0X_1)}}\over V}\otimes {\buildrel {_{(Y_0Y_1)}}\over V}$ is given by 
the following action: 
\begin{equation} 
({\buildrel {_{(X_0X_1)}}\over \Pi}\otimes  {\buildrel {_{(Y_0Y_1)}}\over \Pi})({\cal R}) 
({\buildrel B\over e}_m \otimes {\buildrel C\over e}_n)= 
\sum_{DF}{\buildrel B\over e}_j\otimes {\buildrel F\over e}_p 
\Clebphi{F}{B}{D}{p}{j}{x}\Clebpsi{B}{C}{D}{m}{n}{x} \coefclambda{B}{D}{F}{C}{X_0}{X_1}. 
\end{equation} 
In particular this expression is a finite sum although the universal formula for ${\cal R}$ is  an infinite sum. 
\end{proposition} 
\Proof 
 
The element ${\buildrel A\over L}$ are multipliers of ${\tilde \mathfrak{U}}{}_{q}(sl(2,\mathbb{C})_{\mathbb{R}})$ and we have the trivial relation  
${\buildrel A\over L}{}^{(\pm)}{}^i_j= 
\sum_{B}^{\oplus}{\buildrel {AB}\over R}{}^{ik}_{jl} {\buildrel B\over X}{}^l_k.$ 
The $R$ matrix of ${\mathfrak U}_q(sl(2,\CC)_{\RR})$ is written as: 
${\cal R}=\sum_{A}{\buildrel A\over X}{}^i_j\otimes 1\otimes 1\otimes {\buildrel A\over g}{}^j_i.$ 
From the expression of a representation $\tilde{\Pi}$ of  ${\tilde \mathfrak{U}}{}_{q}(sl(2,\mathbb{C})_{\mathbb{R}})$, associated to  
$(X_0,X_1)$ a straightforward  computation of the representation of the R-matrix gives the formula of the statement. 
  \eoProof 
 
\bigskip

Remark. This situation is in sharp contrast with the case of $SU_{q}(1,1)$ for $q$ real. The structure of 
Hopf algebra on $\mathfrak{U}_{q}(su(1,1))$ is the same  as $\mathfrak{U}_{q}(su(2)),$ the only 
difference is in the definition of the star structure: $J_z^{\star}=J_z, J_+^{\star}=-J_-,J_-^{\star}=-J_+.$ 
It is easy to classify the irreducible unitary representations, the principal unitary representation of  
$\mathfrak{U}(su(1,1))$ being now easily quantized as follows: 
\begin{eqnarray} 
&&\hskip -1cmJ_z.e_m=m e_m\;\;\;\;\;\;J_+.e_m=[m-\tau+\epsilon]e_{m+1}\;\;\;\;\;\;J_-.e_m=-[m+\tau+\epsilon]e_{m+1}\nonumber 
\end{eqnarray} 
where as usual $\epsilon\in\{0,\onehalf\}, \tau=i\nu-\onehalf, \nu\in \mathbb{R},$ and $e_m$ is an orthonormal 
basis of the representation. 
The universal $R$ matrix is still formally defined by  
$R=q^{2 J_z \otimes J_z}e_{q^{-1}}^{(q-q^{-1})\;(q^{J_z}J_+ \otimes J_-q^{-J_z})}.$ 
It is easy to see that $(\pi\otimes \pi')(R)(e_m\otimes e_n')$ has no meaning in the $l^2$ sense except in the 
trivial  case where $q=1.$ 
Indeed a straightforward  computation shows that a formal expansion of $(\pi\otimes \pi')(R)(e_m\otimes e_n')$ gives $(\pi\otimes \pi')(R)(e_m\otimes e_n')=\sum_{p=0}^{+\infty}a_{p}(e_{m+p}\otimes e_{n-p}')$ with $a_{p}\sim (q-q^{-1})\gamma q^{-2p^2+\alpha p+\beta}$ for $0<q<1$ with $\gamma\not=0.$ 
 
\subsection{Connection with Askey-Wilson Polynomials.} 
In the sequel, we will show that we can reexpress  $N^{{}_{(X_0X_1)\;(Y_0Y_1)}}_{\;\;\;{}_{(Z_0Z_1)}}\ElemRed{A}{B}{C}{X_0X_1}{Y_0Y_1}{Z_0Z_1}$ in terms of q-Racah polynomials and Askey-Wilson polynomials, and that the relation (\ref{completeness}) is an orthogonality property mixing q-Racah polynomials and Askey-Wilson polynomials in a non trivial way.

\begin{lemma}  
The following three identities are satisfied: 
\begin{eqnarray}  
&&\hskip -1.5cm 1)\;\;\;\overline{\ElemRed{A}{B}{C}{X_0X_1}{Y_0Y_1}{Z_0Z_1}}=\psi_2\ElemRed{A}{B}{C}{\underline{X_1}\underline{X_0}}{\underline{Y_1}\underline{Y_0}}{\underline{Z_1}\underline{Z_0}}\nonumber\\  
&&\hskip -1.5cm 2)\;\;\;\ElemRed{A}{B}{C}{\underline{X_1}\underline{X_0}}{\underline{Y_1}\underline{Y_0}}{\underline{Z_1}\underline{Z_0}}=\psi_3\ElemRed{A}{B}{C}{X_1X_0}{Y_1Y_0}{Z_1Z_0}\nonumber\\ 
&&\hskip -1.5cm 3)\;\;\;\ElemRed{A}{B}{C}{X_1X_0}{Y_1Y_0}{Z_1Z_0}=e^{-i\pi B}\sqrt{[d_B]}\;  
\sum_{X_2}\sixjn{Y_0}{X_0}{Z_0}{T}{Z_1}{X_2}{3}\!\!  
\sixjn{X_1}{X_0}{A}{C}{B}{X_2}{1}\times\nonumber\\  
&&\hskip 5cm\times\sixjn{Z_1}{X_2}{Y_0}{B}{Y_1}{X_1}{3}\!\!  
\frac{v^{1/4}_{A}v^{1/2}_{X_2}}{v^{1/4}_{B}v^{1/4}_{C}v^{1/4}_{X_0}v^{1/4}_{X_1}}  
\end{eqnarray}  
where $\psi_2=e^{i \pi (B+C+m_X+A)}$, and $\psi_3$ is such that 
\begin{eqnarray}  
&&\psi_1 \psi_2 \psi_3=e^{i \pi (m_X+m_Y-m_Z)}\varphi_{(2X_1+1,-1)}\varphi_{(2Y_1+1,-1)}\varphi_{(2Z_1+1,-1)}\nonumber 
\end{eqnarray}  
\end{lemma}  
  
\Proof  
 
The identity $1)$ follows from a very simple computation. The identity $2)$ follows from the use of propositions (\ref{propsym6j1}) and (\ref{propsym6j3}). Identity $3)$ is a direct consequence of a Racah and a pentagonal equation on $6j(3)$ symbols.  
 \eoProof  

\medskip   
Using these identities theorem (\ref{isometrytheorem}) can be reformulated as follows 
\begin{proposition} \label{askeyproof}  
\begin{eqnarray}  
&&\hskip -1cm\int \!d(Z_0Z_1)\Gamma^{\!{}_{(X_0X_1)(Y_0Y_1)}}_{\;\;\;{}_{(Z_0Z_1)}}\!\sixjn{Y_0}{X_0}{Z_0}{A}{Z_1}{X_2}{3}\!\!\sixjn{Y_0}{X_0}{Z_0}{A}{Z_1}{X_3}{3}\!\!\sixjn{Z_1}{X_2}{Y_0}{B}{Y_1}{X_1}{3}\!\!\sixjn{Z_1}{X_3}{Y_0}{C}{Y_1}{X_1}{3}\!\!=\nonumber\\   
&&=\delta_{X_2,X_3}\delta_{B,C}Y^{(1)}_{(A,X_0-X_2)}Y^{(1)}_{(B,X_1-X_3)}Y^{(1)}_{(B,Y_0-Y_1)}\label{O3}.  
\end{eqnarray}  
with  
\begin{eqnarray}  
&&\Gamma^{\!{}_{(X_0X_1)(Y_0Y_1)}}_{\;\;{}_{(Z_0Z_1)}}=  
\frac{q^{2Z_0-2 Z_1+2Y_0-2 Y_1}(1)_{\infty}^2(1-q^2)[d_B]}{(2X_1\!+\!1)_1(2Y_0\!+\!1)_1(2Z_0\!+\!1)_1}\times\nonumber\\ 
&&\times \frac{\xi(2X_1\!+\!1,2Y_1\!+\!1,2Z_1\!+\!1)} 
{\xi(Y_1\!+\!X_1\!-\!Z_1\!+\!1,X_1\!+\!Z_1\!-\!Y_1\!+\!1,Y_1\!+\!Z_1\!-\!X_1\!+\!1,X_1\!+\!Y_1\!+\!Z_1\!+\!2)} 
\end{eqnarray}  
\end{proposition}  
  
\Proof  
 
Straightforward using a Racah-Wigner symmetry on $6j(3)$ and an orthogonality relation on $6j(1)$ symbols.  
 \eoProof  
  
\begin{theorem}  
Using the notations (\ref{defaskeypol}) we have   
\begin{eqnarray}  
&&\hskip -1cm\sqrt{{\cal P}(Z_0Z_1)}\sqrt{\Gamma^{\!{}_{(X_0X_1)(Y_0Y_1)}}_{\;\;{}_{(Z_0Z_1)}}}\sixjn{Z_1}{X_2}{Y_0}{B}{Y_1}{X_1}{3}\!\!\!= \!\sqrt{\hbar}  
\sqrt{\frac{w^{(AW)}(z;a,b,c,d)}{h^{(AW)}_n(a,b,c,d)}}p^{(AW)}_n(\tau(z);a,b,c,d) \label{expressionaskey} 
\end{eqnarray}  
with $n=B+Y_0-Y_1, a=q^{2X_1+2Y_1+3}, b=q^{2Y_1-2X_1+1}, c=q^{2X_2-2Y_0+1}, d=q^{-2Y_0-2X_2-1}, 2z=2 Z_1+1.$   
\end{theorem}  
  
\Proof  
 
Using two Sears transformations we can turn the parameters of the hypergeometric functions of the left handside into those of the right handside. Checking now the identity (\ref{expressionaskey}) is therefore reduced to straightforward but tedious manipulations on $\nu_{\infty}$ functions. 
 \eoProof  
 
\medskip 
 
Note that the relation (\ref{O3}) when $A=0$ reduces to the identity: 
\begin{equation} 
\int d\rho_Z \;{\cal P}(0,\rho_Z)\;\Gamma^{(X_0X_1)\;(Y_0Y_1)}_{(Z_0Z_0)} \sixjn{Z_0}{X_0}{Y_0}{B}{Y_1}{X_1}{3}\sixjn{Z_0}{X_0}{Y_0}{C}{Y_1}{X_1}{3}=\delta_{B,C} 
\end{equation} 
which is exactly, using (\ref{expressionaskey}), the orthogonality condition on Askey-Wilson polynomials for the family of parameters $n=B+Y_0-Y_1, a=q^{2X_1+2Y_1+3}, b=q^{2Y_1-2X_1+1}, c=q^{2X_0-2Y_0+1}, d=q^{-2Y_0-2X_0-1}, 2z=2 Z_0+1.$ 
 
Although one would suspect that the relation (\ref{O3}) can be proved for any $A$ by disentangling the sum over $m_Z$ and the integration over $\rho_Z,$ we have not been able to prove it in this way.

\section{Conclusion}  
  
In this work we have given exact formulae for the Clebsch-Gordan coefficients of the tensor  
product of two principal unitary representations of the quantum Lorentz Group. We have found explicit expressions of the intertwiners in terms of q-Racah polynomials and Askey-Wilson polynomials.  
A consequence of this relation is the proof of the proposition (\ref{askeyproof}) which contains in its simplest case the proof of orthogonality of Askey-Wilson polynomials.  
One should generalize our work to  obtain explicit expressions for intertwiners of two representations in the set of principal and complementary series. This should as well gives non trivial relations on orthogonal polynomials.  
  
A very interesting question is the generalization  to other quantization of complex semi-simple  Lie algebras. In the classical case the structure of the tensor product of two principal representations is known \cite{Wi} and  there also exists polynomials in several variables generalizing Askey-Wilson polynomials \cite{Ko}.   
  
Results \cite{BR2} which are nevertheless at hand are the construction of 6j(6) in terms of non terminating basic hypergeometric functions, the study of their properties and their understanding as matrix elements of fusion matrices.  
Using these $6j(6),$ we will directly obtain expressions for the 6j of the principal representations of the quantum Lorentz group as well as very interesting  pentagonal equations satisfied by them.

\section{Appendix}  
  
\subsection{Formulae on Basic Hypergeometric Functions.}

In the sequel we shall frequently use the following notations:  
  
$\forall x \in \CC, \forall k \in \NN,$  
\begin{eqnarray*}   
&&\hskip -1.2cm [x]=\frac{q^x-q^{-x}}{q-q^{-1}}\;,\;\;\;\;d_x=2x+1\;,\;\;\;\;[x]_k=\prod_{n=1}^{ k} [x+n-1]\;,\;\;\;\;[k]!=[1]_k,\;\;\;\;\;v_x^{1/4}=e^{i\frac{\pi}{2}x}q^{-\frac{1}{2}x(x+1)}.  
\end{eqnarray*}  
In this  article, the square root of a complex number is defined  by:  
\begin{eqnarray}  
\forall x\in \CC, \sqrt{x}=\sqrt{\vert x \vert}e^{i \frac{Arg(x)}{2}},  
\mbox{where}\;\; x=\;\vert x \vert e^{i  Arg(x)}, Arg(x) \in ]-\pi,\pi],   
\end{eqnarray}  
and for all complex number $z$ with non zero imaginary part,  we define $\epsilon(z)=\;\mbox{sign}\; (Im(z)).$  

 We will define  the following basic functions: $\forall \alpha,\beta,\gamma \in \CC, \forall n \in {\mathbb Z},$  
\begin{eqnarray*}  
&&(\alpha)_{\infty}=(q^{2\alpha},q^2)_{\infty}=  
\prod_{k=0}^{+\infty}(1-q^{2\alpha+2k})\;,\;\;\;\;\;\;\;\;\;\;\;\;  
(\alpha)_{n}=\frac{(\alpha)_{\infty}}{(\alpha+n)_{\infty}},\\  
&&\nu_{\infty}(\alpha)=\prod_{k=0}^{+\infty}\sqrt{1-q^{2\alpha+2k}}  
\;,\;\;\;\;\;\;\;\;\;\;\;\;\;\;\;\;\;\;\;\;\;\;\;\;\;\;\;\;\;\;\;\;\nu_{n}(\alpha)=\frac{\nu_{\infty}(\alpha)}{\nu_{\infty}(\alpha+n)},  
\end{eqnarray*}  
\begin{eqnarray*}  
&&\;\;\;\;\;  
\omega(\alpha;\beta,\gamma)=\frac{\nu_{\infty}(\alpha+\beta+\gamma+2)\nu_{\infty}(\alpha+\beta-\gamma+1)\nu_{\infty}(\alpha-\beta+\gamma+1)}{\nu_{\infty}(-\alpha+\beta+\gamma+1)}.  
\end{eqnarray*}  
  
The q-factorials satisfy the following relations:  
 $\forall \alpha \in \CC, \forall n \in \ZZ,$  
\begin{eqnarray}  
&&(\alpha)_{-n}(1-\alpha)_{n}=(-1)^{n}q^{-2n\alpha +n(n+1)} \;\;\mbox{(inversion relation)},\label{inversionrelation}\\  
&&\lim_{\alpha \rightarrow +\infty}(\alpha)_{n}=1,\;\;\;\;\;\;\;\;\;\;\;\;  
(\alpha)_{n}\sim_{\alpha \rightarrow -\infty} (-1)^{n}q^{2n\alpha+n(n-1)}.  
\end{eqnarray}

In order to control the signs of the expressions in our article, we are led to introduce the following functions $\varphi:\CC\times{\onehalf}\ZZ\rightarrow \CC$, defined by:    
\begin{equation}  
\varphi_{(\alpha,n)}=\frac{\nu_{\infty}(\alpha-n+1)\nu_{\infty}(n-\alpha)}{\nu_{\infty}(\alpha+1)\nu_{\infty}(-\alpha)}q^{-n\alpha+\onehalf n(n-1)}\;\;,\forall \alpha\in \CC, \forall n\in {\onehalf}\ZZ.  
\end{equation}  
This function  $\varphi$ satisfies the two following relations:   
\begin{eqnarray}  
&&\varphi_{(\alpha, n)}\varphi_{(\alpha-n, p)}=\varphi_{(\alpha, n+p)}, \forall\alpha \in \CC, \forall n,p\in {\onehalf}\ZZ\\  
&&\frac{\varphi_{(\alpha,n)}^2}{\varphi_{(\alpha,p)}^2}=(-1)^{n-p}, \forall \alpha \in \CC,\forall n,p \in \onehalf\ZZ\;, (n-p) \in \ZZ,   
\end{eqnarray}  
The last equation implies in particular that $\varphi_{(\alpha,n)}^2=(-1)^{n},   
\forall \alpha \in \CC,\forall n\in \ZZ.$ Note that when $n$ is  half an odd integer,  
 $ \varphi_{(\alpha,n)}$ is not a fourth root of unity.  
  
With our choice of square root, we have the more precise  results:  
\begin{eqnarray}  
&&\varphi_{(\alpha,n)}=e^{-i\frac{\pi}{2}\; n \; \epsilon(q^{\alpha})}, \forall n \in\ZZ, \forall \alpha\in\CC,\\  
&&\vert \varphi_{(-1+i\rho, n)} \vert^2 = e^{i\pi (p+1/2+\epsilon(q^{i\rho}))}\varphi_{(i\rho, -p)}^2,\forall n,p \in \onehalf+\ZZ \;, \forall \rho \in \RR.  
\end{eqnarray}   
  
We will also make an extensive use of the following basic hypergeometric functions, associated to complex numbers $\alpha_0, \cdots, \alpha_n, \beta_1, \cdots, \beta_n$ and defined by:  
\vskip -0.5cm  
\begin{eqnarray*}  
&&\forall Z\in \CC, \;\;\;\;\;\;\;\;\;\;\;\;\;\;\;{}_{n+1}\!\Phi_n  
\left[\!\!\begin{array}{ccccc} \alpha_0 &\!\!\alpha_1 &\!\!\cdots &\!\!\!\!\alpha_n &\\  
&\!\!\!\!\!\!\!\!\beta_1 &\!\!\!\!\!\!\!\!\cdots &\!\!\!\!\!\!\!\beta_n & \; ; Z\end{array} \!\!\right]\!=\sum_{k=0}^{+\infty}\frac{\prod_{i=0}^{n}(\alpha_i)_k}{(1)_k\prod_{i=1}^{n}(\beta_i)_k}q^{2kZ},\nonumber\\  
&&{}_{n+3}W_{n+2}(\alpha_0\; ; \;\alpha_1, \cdots, \alpha_n \; ; \; Z)\!= \sum_{k=0}^{+\infty}\frac{(1-q^{2\alpha_0+4k})}{(1-q^{2\alpha_0})}\frac{\prod_{i=0}^{n}(\alpha_i)_k}{(1)_k\prod_{i=1}^{n}(1+\alpha_0-\alpha_i)_k}q^{2kZ}.  
\end{eqnarray*}  
In an expression involving  ${}_{n+1}\!\Phi_n$ or ${}_{n+1}\!W_n,$ if $Z$ is not specified it will mean that $Z=1$. \\  
  
Let us recall some  properties of basic hypergeometric functions which are proved for example in \cite{GR}:  
\begin{eqnarray}  
&&\hskip -1cm {}_{2}\!\Phi_1  
\left[\!\!\begin{array}{ccc} \alpha_0 &\!\!\alpha_1, &\\  
&\!\!\!\!\!\!\!\!\beta_1 & \; ; \beta_1\!-\!\alpha_0\!-\!\alpha_1\end{array} \!\!\right]=  
\frac{(\beta_1\!-\!\alpha_0)_{\infty}(\beta_1\!-\!\alpha_1)_{\infty}}{(\beta_1)_{\infty}(\beta_1\!-\!\alpha_0\!-\!\alpha_1)_{\infty}}\label{heine},\\  
&&\hskip -1cm{}_{4}\!\Phi_3\left[\!\!\begin{array}{cccc} -n &\!\!\alpha_1 &\!\!\alpha_2 &\!\!\!\!\alpha_3 \\&\!\!\!\!\!\!\!\!\beta_1 &\!\!\!\!\!\!\!\!\beta_2 &\!\!\!\!\!\!\!\beta_3 \end{array} \!\!\right]\!=  
\frac{ {}_{4}\!\Phi_3\left[\!\!\begin{array}{cccc} -n &\;\;\;\;\;\alpha_1 &\;\;\;\;\;\beta_1\!-\!\alpha_2 &\;\;\;\;\;\beta_1\!-\!\alpha_3 \\ &\!\!\!\!\!\!\!\!\beta_1 &\!\!1\!-\!n\!+\!\alpha_1\!-\!\beta_2 &1\!-\!n\!+\!\alpha_1\!-\!\beta_3 \end{array} \!\!\right]  
}{q^{2n\alpha_1} (\beta_2-\alpha_1)^{-1}_{n}(\beta_3-\alpha_1)^{-1}_{n}  (\beta_2)_{n}(\beta_3)_{n}}  
\!\label{sears},\\  
&&\hskip -1cm{}_{6}W_{5}(\alpha_0\; ; \;\alpha_1, \alpha_2, -n \; ; \; 1\!+\!n\!+\!\alpha_0\!-\!\alpha_1\!-\!\alpha_2 )=  
\frac{(1+\alpha_0)_{\infty}(1\!+\!\alpha_0\!-\!\alpha_1\!-\!\alpha_2)_{\infty}}{(1\!+\!\alpha_0\!-\!\alpha_1)_{\infty}(1\!+\!\alpha_0\!-\!\alpha_2)_{\infty}}.\label{phi65}  
\end{eqnarray}  
Relation (\ref{heine}) is called Heine formula, relation (\ref{sears}) is the Sears  
 Identity.

Let us now recall the definitions  of  $q-$Racah polynomials $p_n^{(R)}$, of the $q-$Racah discrete measure $w^{(R)}$, and of the square of the norm of these polynomials \cite{W}\cite{AW1}:  
\begin{eqnarray}  
&&\hskip -0.5cm \forall x\in \CC, \forall \alpha,\beta,\gamma,\delta\in \CC,  
\mu(x)=q^{-2x}\!+\!q^{2x+2+2\gamma+2\delta},\nonumber\\  
&&\hskip -0.5cm p_n^{(R)}(\mu(x);q^{2\alpha},q^{2\beta},q^{2\gamma},q^{2\delta})={}_{4}\!\Phi_3  
\left[\!\!\begin{array}{cccc} -n & n\!+\!1\!+\!\alpha\!+\!\beta & -x  & \gamma\!+\!\delta\!+\!x\!+\!1 \\  
&\!\!\!\!\!\!\!\!\!\alpha\!+\!1 &\!\!\!\!\!\!\!\!\beta\!+\!\delta\!+\!1 &\!\!\!\!\!\!\!\gamma\!+\!1\end{array} \!\!\right],\nonumber\\  
&&\hskip -0.5cm w^{(R)}(x;q^{2\alpha},q^{2\beta},q^{2\gamma},q^{2\delta})=\frac{(\gamma\!+\!\delta\!+\!1,\alpha\!+\!1,\beta\!+\!\delta\!+\!1,\gamma\!+\!1)_x (1-q^{2(2x+\gamma+\delta+1)})}{q^{2x(\alpha+\beta+1)}(1,\gamma\!+\!\delta\!+1-\!\alpha,\gamma\!-\!\beta\!+\!1,\delta\!+\!1)_x(1-q^{2(\gamma+\delta+1)})}\nonumber,\\  
&&\hskip -0.5cm h^{(R)}_n(q^{2\alpha},q^{2\beta},q^{2\gamma},q^{2\delta})=  
\frac{q^{2n(\gamma+\delta+1)}  
(1-q^{2(\alpha+\beta+1)})  
(1,\beta\!+\!1,\alpha\!-\!\delta\!+\!1,\alpha\!+\!\beta\!-\!\gamma\!+\!1)_n}  
{(1-q^{2(2n+\alpha+\beta+1)})  
(\alpha\!+\!\beta\!+\!1,\alpha\!+\!1,\beta\!+\!\delta\!+\!1,\gamma\!+\!1)_n}  
\times\nonumber\\  
&&\hskip 3.5cm \times \frac{(\gamma\!+\!\delta\!+\!2,\gamma\!-\!\alpha\!-\!\beta,\delta\!-\!\alpha, -\beta)_{\infty}}  
{(\gamma\!+\!\delta\!-\!\alpha\!+\!1,\gamma\!-\!\beta\!+\!1,\delta\!+\!1,\!-\!\alpha\!-\!\beta\!-\!1)_{\infty}}\label{defracahpol}.
\end{eqnarray}  
The definition of  Askey-Wilson polynomials $p_n^{(AW)},$  of the Askey-Wilson measure $w^{(AW)}$ and of the square of the norm of these polynomials is \cite{AW2}:  
\begin{eqnarray}  
&&\hskip -1cm \forall z\in \CC, \forall \alpha,\beta,\gamma,\delta\in \CC,\tau(z)=\frac{1}{2}(q^{2z}+q^{-2z})\nonumber\\  
&&\hskip -1cm \frac{p_n^{(AW)}(\tau(z);q^{2\alpha},q^{2\beta},q^{2\gamma},q^{2\delta})}{  
q^{-2n\alpha}(\alpha\!+\!\beta,\alpha\!+\!\gamma,\alpha\!+\!\delta)_n   
}=\;{}_{4}\!\Phi_3 \!  
\left[\!\!\begin{array}{cccc} -n & n\!+\!\alpha\!+\!\beta\!+\!\gamma\!+\!\delta\!-\!1 & \alpha\!+\!z  & \alpha\!-\!z  \\  
&\!\!\!\!\!\!\!\!\!\!\!\!\!\!\!\!\!\!\!\!\!\!\!\!\!\!\!\!\!\!\!\!\!\!\!\!\alpha\!+\!\beta &\!\!\!\!\!\!\!\!\!\!\!\!\!\!\!\!\!\!\!\!\!\!\!\!\!\!\!\!\!\!\!\!\!\!\!\!\!\!\!\!\!\!\alpha\!+\!\gamma\! &\!\!\!\!\!\!\!\!\!\!\!\!\!\!\!\!\alpha\!+\!\delta\end{array} \!\!\right],  
\nonumber 
\end{eqnarray} 
\vskip -0.3cm 
\begin{eqnarray} 
&&\hskip -1cm w^{(AW)}(z;q^{2\alpha},q^{2\beta},q^{2\gamma},q^{2\delta})=  
\frac{  
(2z, -2z)_{\infty}}  
{(\alpha+z,\alpha\!-\!z,\beta\!+\!z,\beta\!-\!z,\gamma\!+\!z,\gamma\!-\!z,\delta\!+\!z,\delta\!-\!z)_{\infty}},\label{defaskeypol}\\  
&&\hskip -1cm h^{(AW)}_n(q^{2\alpha},q^{2\beta},q^{2\gamma},q^{2\delta})=  
\frac{(\alpha\!+\!\beta\!+\!n)_{\infty}^{-1}(n+1)_{\infty}^{-1}  
(\alpha\!+\!\beta\!+\!\gamma\!+\!\delta\!+\!2n)_{\infty}  
(\alpha\!+\!\beta\!+\!\gamma\!+\!\delta\!+\!n\!-\!1)_{n}}  
{(\alpha\!+\!\gamma\!+\!n,  
\alpha\!+\!\delta\!+\!n,\beta\!+\!\gamma\!+\!n,\beta\!+\!\delta\!+\!n,\gamma\!+\!\delta\!+\!n)_{\infty}}.\nonumber 
\end{eqnarray}  
  
\subsection{Asymptotics of Intertwiners}  
  
Our aim is to compute the behaviour of $\ElemRed{\!\!P+T}{\;\;P}{T}{X_0X_1}{Y_0Y_1}{Z_0Z_1}$   
(see (\ref{form6jprol1}) for the explicit expression) for $P=Y_0-Y_1\geq 0, T\rightarrow +\infty.$ The sum   
on $X_2$ being finite, we can work  term by term.\\  
The asymptotics of  $\sixjn{X_1}{X_0}{P+T}{T}{P}{\;\;X_2}{1}$ is very easy to obtain because, after a Racah-Wigner symetry (\ref{symcont}), the hypergeometric part of the result tends to one when $T$ is large.  The  behaviour is therefore given by:  
\begin{eqnarray}  
&&\hskip -1cm\sixjn{X_1}{X_0}{P+T}{T}{P}{\;\;X_2}{1}\!\!\sim   
\frac{\varphi_{(X_0+X_2-T,P+X_2-X_1)}\nu_1(2X_2+1)\omega(P;X_1,X_2)}  
{q^{\frac{1}{2}(P-3X_2-X_1-1)(P+X_2-X_1)} \nu_{\infty}(2P+1)\nu_{\infty}(1)}\nonumber  
\end{eqnarray}  
The behaviour of  $\sixjn{Y_0}{X_0}{Z_0}{T}{Z_1}{X_2}{3}$ requires  more attention, but can be achieved using the following computation. 
 
 By using first  the Sears identity (\ref{sears}), and then taking the limit when $T$ goes to infinity, we easily conclude  using Heine summation formula (\ref{heine}):   
\begin{eqnarray}  
&&\hskip -0.5cm\frac{(Z_1\!+\!Y_0\!-\!X_0\!-\!T\!+\!1)_{\infty}{}_4\Phi_3\!\left[\!\!\begin{array}{cccc} \mbox{\it \small $Z_1\!-\!Z_0-\!$T} &\; \mbox{\it \small $Z_1\!+\!Z_0-\!$T$\;+1$} & \;\;\;\;\;\mbox{\it \small $X_2\!-\!X_0-\!$T} & \;\;\;\mbox{\it \small $-X_2\!-\!X_0\!-\!$T$-\!$1}\\ & \hskip -3.5cm\mbox{\it \small $-\!$2T} & \hskip -3.5cm\mbox{\it \small $\!-\!X_0\!-\!Y_0-\!$T$\;+\!Z_1$} & \hskip -1.8cm \mbox{\it \small $Z_1+Y_0\!-\!X_0\!-\!$T$\;+1$} \end{array} \!\!\right]}  
{(Z_0\!+\!Y_0\!-\!X_0\!+\!1,Y_0\!-\!X_0\!-\!Z_0,2Z_1\!+\!2)_{\infty}}=\nonumber 
\end{eqnarray} 
\begin{eqnarray}  
&&\hskip 0.5cm = \frac{(2T\!+\!1)_{\infty}{}_4\Phi_3\!\left[\!\!\begin{array}{cccc} \mbox{\it \small $Z_1\!-\!Z_0-\!$T} &\; \mbox{\it \small $Z_1\!+\!Z_0-\!$T$\;+$1} & \;\;\;\;\;\mbox{\it \small $Z_1\!-\!X_2-\!Y_0$} & \;\;\;\mbox{\it \small $X_2\!+\!Z_1\!-\!Y_0\!+\!1$}\\ & \hskip -3.5cm\mbox{\it \small $2Z_1\!+\!2$} & \hskip -3.5cm\mbox{\it \small $\!-\!X_0\!-\!Y_0-\!$T$\;+\!Z_1$} & \hskip -1.8cm \mbox{\it \small $Z_1+X_0\!-\!Y_0\!-\!$T$\;+1$} \end{array} \!\!\right]  
}  
{(T\!+\!Y_0\!-\!Z_1\!-\!X_0,T\!+\!Z_1\!+\!Z_0\!+\!2,Z_1\!+\!T\!-\!Z_0+1  
)_{\infty}}\!\!\nonumber\\  
&&\hskip 1.5cm{\sim}_{{}_{T \rightarrow +\infty}} \;\;{}_2\Phi_1\!\left[\!\!\begin{array}{lll} \mbox{\it \small $Z_1\!-\!X_2-\!Y_0$} & \;\;\;\mbox{\it \small $X_2\!+\!Z_1\!-\!Y_0\!+\!1$} &\\ &  
 \hskip -0.5cm\mbox{\it \small $2Z_1\!+\!2$} & \; ;2Y_0+1\end{array} \right]  
\nonumber\\  
&&\hskip 1.5cm{\sim}_{{}_{T \rightarrow +\infty}} \;\;\frac{(Z_1\!+\!Y_0\!+\!X_2\!+2,Z_1\!+\!Y_0\!-\!X_2\!+1)_{\infty}}{(2Y_0\!+\!1,2Z_1\!+\!2)_{\infty}}.\nonumber  
\end{eqnarray}  
This gives us the following behaviour  
\begin{eqnarray}  
&&\hskip -1cm\sixjn{Y_0}{X_0}{Z_0}{T}{Z_1}{X_2}{3}\!\!\!\sim e^{i \pi (T+X_0-X_2)}q^{-2TY_0}\frac{\varphi_{(X_0+X_2-T,-T)}\varphi_{(Z_0+Z_1-T,-T)}\omega(Y_0;X_2,Z_1)\omega(Y_0;X_0,Z_0)}{q^{(Z_1-Z_0)(X_0+Z_0-Y_0)+(X_2-X_0)(Z_1+X_0-Y_0+1)}(1,2Y_0\!+\!1)_{\infty}}\times\nonumber\\  
&&\hskip 1cm\times \frac{\nu_1(2X_2\!+\!1)\nu_1(2Z_0\!+\!1)(Y_0\!-\!X_0\!-\!Z_0,X_0\!+\!Z_0\!-\!Y_0\!+\!1)_{\infty}  
}{\nu_{\infty}(-\!X_0\!-\!X_2)\nu_{\infty}(X_0\!+\!X_2\!+\!1)\nu_{\infty}(-\!Z_0\!-\!Z_1)\nu_{\infty}(Z_0\!+\!Z_1\!+\!1)}\label{asymptotic6j3} 
\end{eqnarray}  
The expression of $\sixjn{Z_1}{X_2}{Y_0}{P}{Y_1}{X_1}{3}$ being drastically simplified for $P=Y_0-Y_1,$ there is no problem anymore to obtain the behaviour of the whole expression,  except  that in order to have a simple answer, we have to take the sum over $X_2.$ Precisely, denoting $\alpha=2X_1\!+\!2Y_1\!-\!2Y_0\!+\!1, \;m=\!2Y_0\!-\!2Y_1\!, \;k=\!Y_0\!-\!Y_1\!+\!X_2\!-\!X_1$ we already have obtained  the following behaviour when $T \rightarrow +\infty$  
\begin{eqnarray}  
&&\hskip -1.2cm\ElemRed{$\!\!\!T\!+\!Y_0\!\!-\!\!Y_1$}{$\!\!Y_0\!\!-\!\!Y_1$}{T}{X_0X_1}{Y_0Y_1}{Z_0Z_1}  
\!\!{\sim}\;  
q^{-T(Y_0+Y_1)}\;e^{i \pi (\frac{X_0-X_1}{2})}   
q^{-Z_0Z_1}\!\!\times\nonumber\\  
&&\hskip -0.8cm\times \sqrt{\qd{Y_0-Y_1}}\nu_1(2X_1\!+\!1)\nu_1(2Z_0\!+\!1)  
\frac{\varphi_{(2X_0,P\!+\!T\!+\!X_0\!-\!X_1)}\varphi_{(Z_1\!+\!X_1\!-\!Y_1,Z_1\!-\!Z_0\!+\!X_1\!-\!X_0\!+\!Y_0\!-\!Y_1)}}{\varphi_{(Z_0\!+\!Z_1,T)}}  
\times\nonumber\\  
&&\hskip -0.8cm\times \;  
\frac{\nu_{\infty}(Y_0\!+\!X_0\!-\!Z_0\!+\!1)  
\nu_{\infty}(Y_0\!+\!X_0\!+\!Z_0\!+\!2)\nu_{\infty}(Y_0\!-\!X_0\!-\!Z_0)  
\nu_{\infty}(Y_0\!+\!Z_0\!-\!X_0\!+\!1)}{  
q^{\frac{3}{2}X_0^2+\onehalf X_0-\onehalf Z_0^2+\onehalf Z_0}  
(1)_{\infty}  
\nu_{\infty}(2Y_1\!+\!1)\nu_{\infty}(2Y_0\!+\!1)}\times\nonumber\\  
&&\hskip -0.8cm\times \;  
\frac{\nu_{\infty}(Y_1\!+\!X_1\!-\!Z_1\!+\!1)\nu_{\infty}(Y_1\!+\!X_1\!+\!Z_1\!+\!2)  
\nu_{\infty}(Y_1\!-\!X_1\!-\!Z_1)\nu_{\infty}(Y_1\!+\!Z_1\!-\!X_1\!+\!1)}{  
q^{-\frac{3}{2}X_1^2-\onehalf X_1-\onehalf Z_1^2 -\onehalf Z_1}  
\nu_{\infty}(2X_0\!+\!1)  
\nu_{\infty}(-\!2X_0)\nu_{\infty}(Z_0\!+\!Z_1+1)\nu_{\infty}(-\!Z_0\!-\!Z_1)}\times A\label{}  
\end{eqnarray} where   
 
$A=\frac{1}{(1+\alpha)_m}\sum_{k} \frac{(1-q^{2\alpha+4k})  
(\alpha)_k(-m)_k}{(1-q^{2\alpha})(1)_k(\alpha+m+1)_k} q^{2k^2+2k(m+\alpha)}.$  
  
Using the formula  (\ref{phi65}) and taking in it the limit where   $\alpha_1, \alpha_2 \rightarrow -\infty,$ we obtain that $A=1.$  
This gives the behaviour of $\ElemRed{$\!\!\!T\!+\!Y_0\!\!-\!\!Y_1$}{$\!\!Y_0\!\!-\!\!Y_1$}{T}{X_0X_1}{Y_0Y_1}{Z_0Z_1}$ when $T$ goes to infinity.  
  
\subsection{Formulae on coefficients $\Lambda^{BC}_{AD}(X_0X_1)$}  
  
Let us recall some properties of the coefficients $\Lambda^{BC}_{AD}(X_0X_1)$ which were  proved in   
\cite{BR1}.  
We will say that the coefficients $\Lambda^{BC}_{AD}$ are off-diagonal if $A \not= D$ and diagonal  if $A=D$ (in the latter case  we will use the notation  $\Lambda^{BC}_{A}=\Lambda^{BC}_{AA}$). The coefficients $\Lambda^{BC}_{A}$ are said to be on the boundary if $B\!+\!C\!=\!A,\; A\!+\!B\!=\!C, \;$ or $A\!+\!C\!=\!B.$  
The coefficients  $\Lambda^{BC}_{AD}(X_0X_1)$ are proportional to  $Y^{(0)}_{(A,B,C)}Y^{(0)}_{(D,B,C)}Y^{(1)}_{(A,X_0-X_1)}Y^{(1)}_{(D,X_0-X_1)},$  so  they vanish according to these selection rules.   
 The following relation holds true (it is in fact equivalent to the property that ${\buildrel {\!\!{}_{(X_0X_1)}}\over {\Pi}}$ is a representation of ${\frak U}_q(an(2))$):  
\begin{eqnarray}  
&&\hskip -1cm \coefblambda{A}{C}{F}{G}\coefblambda{B}{D}{G}{H}=\!\!\sum_{KU} \frac{e^{i \pi (G+U)}(\qn{d_U}\qn{d_G})^{\frac{1}{2}}\;\;\coefblambda{K}{U}{F}{H}}{e^{i \pi (C+D)}(\qn{d_C}\qn{d_D})^{\frac{1}{2}}}  
\sixjn{A}{D}{U}{B}{C}{G}{0}\!\!\!\sixjn{B}{U}{C}{F}{A}{K}{0}\!\!\!\sixjn{H}{U}{K}{A}{B}{D}{0}.\label{systan}  
\end{eqnarray}   
  
As a result, by taking $F=H$ in this expression, the off-diagonal coefficients are obtained, up to a sign, from the diagonal coefficients.  \medskip
  
These diagonal coefficients can be computed in various ways. They satisfy   the following system of linear equations:  
 \begin{eqnarray}   
&&\hskip -1.2cm \forall (A,B,C)\in \halfinteger\times\halfinteger\times  {\onehalf}\mathbb{Z}^{+*}\;/\;B\!+\!A\!-\!C,\; B\!+\!C\!-\!A, \;A\!+\!C\!-\!B \in \NN^*\nonumber\\  
&&\hskip -1.3cm   
 \frac{[\mbox{\small \it A$+$B$+$C$+1$}]}{[\mbox{\small \it $2$A$+1$}]}\coefalambda{B}{\;C}{\;\;A}\!+\!    
\frac{[\mbox{\small \it A$+$C$-$B$+1$}]}{[\mbox{\small \it $2$A$+1$}]}\coefalambda{B\!-\!1}{\;C}{\;\;\;A}\!\!-\!\frac{[\mbox{\small \it C$+$B$-$A$-1$}]}{[\mbox{\small \it $2$A$+1$}]}\coefalambda{B\!-\!1}{\;C\!-\!1}{\;\;\;\;A}  
\!+\!\frac{[\mbox{\small \it B$+$A$-$C$+1$}]}{[\mbox{\small \it $2$A$+1$}]}\coefalambda{B}{\;C\!-\!1}{\;\;A}=\nonumber\\   
&&\hskip -1.4cm =  
\frac{[\mbox{\small \it A$+$B$-$C$-1$}]}{[\mbox{\small \it $2$A$-$}1]}\coefalambda{B\!-\!1}{\;\;C}{\;\;\;A\!-\!1}\!+\!\frac{[\mbox{\small \it A$+$C$-$B$-1$}]}{[\mbox{\small \it $2$A}-1]}\coefalambda{B}{\;\;C\!-\!1}{\;\;A\!-\!1}\!+\!  
\frac{[\mbox{\small \it A$+$B$+$C$-1$}]}{[\mbox{\small \it $2$A}-1]}\coefalambda{B\!-\!1}{\;C\!-\!1}{\;\;\;A\!-\!1}  
\!\!-\!\frac{[\mbox{\small \it C$+$B$-$A$+$}1]}{[\mbox{\small \it $2$A$-$}1]}\coefalambda{B}{\;C}{A\!-\!1},  
\label{linearANexpli}  
\end{eqnarray}   
  
as well as  $\forall A,B,C \in \halfinteger / Y^{(0)}_{(A,B,C)}=1,$  
\begin{eqnarray}  
&&\hskip -0.5cm [B\!\!+\!\!C\!\!-\!\!A\!\!+\!\!1][A\!\!+\!\!B\!\!+\!\!C\!\!+\!\!2]q^{\pm(C-B)}\Lambda^{B+\onehalf\;C+\onehalf}_{\;\;A}  
-[A\!\!+\!\!C\!\!-\!\!B][A\!\!+\!\!B\!\!-\!\!C\!\!+\!\!1]q^{\mp(C+B+1)}\Lambda^{B+\onehalf\;C-\onehalf}_{\;\;A}+\nonumber\\  
&&\hskip -0.4cm+[A\!\!+\!\!C\!\!+\!\!B\!\!+\!\!1][B\!\!+\!\!C\!\!-\!\!A]q^{\mp(C-B)}\Lambda^{B-\onehalf\;C-\onehalf}_{\;\;A}  
-[A\!\!+\!\!B\!\!-\!\!C][A\!\!+\!\!C\!\!-\!\!B\!\!+\!\!1]q^{\pm(C+B+1)}\Lambda^{B-\onehalf\;C+\onehalf}_{\;\;A}\nonumber\\  
&&\hskip -0.5cm=\omega_{\pm}[2C\!\!+\!\!1][2B\!\!+\!\!1]\Lambda^{B\;C}_{\;\;A},\!\!\!\!\label{lambdacasimir}  
\end{eqnarray}  
where we have denoted   $\omega_+=q^{2X_0+1}\!+\!q^{-2X_0-1},\;\omega_-=q^{2X_1+1}\!+\!q^{-2X_1-1}.$  
 
\medskip  
 
These two linear systems have a unique normalized solution.  
The value of the diagonal coefficients on the boundary have the following simple expression:  
  
\vskip -0.4cm  
\begin{eqnarray}  
&&\hskip -1cm\Lambda^{B\;\;\;C}_{B+C}(X_0X_1)= \sum_{k=-B}^{B}\frac{q^{-2k(X_0+X_1+1)}\binom{B+C+X_0-X_1}{B+k}_q\binom{B+C-X_0+X_1}{B-k}_q}{(-1)^{2B}\binom{2B+2C}{2B}_q},\label{lambdaboundary}\\  
&&\hskip -1cm\Lambda^{B\;A+B}_{\;\;\;A}(X_0X_1)=\!\!\!\sum_{k=-B}^{B}\frac{q^{-2k(X_0+X_1+1)}\binom{A+X_0-X_1+B-k}{B-k}_q\binom{A-X_0+X_1+B+k}{B+k}_q}{\binom{2A+2B+1}{2B}_q}\nonumber\\ 
&&= \frac{v_{X_0}}{v_{X_1}}\Lambda^{A+B\;B}_{\;\;\;A}(X_1X_0)\label{symetryLambda}. 
\end{eqnarray}  
  
Different explicit formulae for the coefficients $\Lambda^{BC}_{AD}$ can be obtained, the fundamental  one is expressed in terms of 6j(1) (\ref{formlamb6j}).  This last expression can  be simplified, using the universal shifted cocycle associated  to $6j(1),$  and we get:   
\begin{eqnarray}  
&&\hskip -1.8cm\Lambda^{BC}_{AD}(X_0,X_1)  
= \sum_{\sigma_1 \sigma_2}\frac{v_A^{1/4}}{v_D^{1/4}}  
\frac{{\cal N}^{(D)}(q^{2X_1+1},m_X)}{{\cal N}^{(A)}(q^{2X_1+1},m_X)}\Clebpsi{C}{B}{A}{\sigma_2}{\sigma_1}{m_X}\;  
q^{- 2i\sigma_1 \rho_X}  
\Clebphi{B}{C}{D}{\sigma_1}{\sigma_2}{m_X}\label{lambda3j}  
\end{eqnarray}  
where ${\cal N}^{(A)}(q^{2X+1},\sigma)$ are normalization factors, given in \cite{BR1}, which values will not be useful in our paper. This last expression is simpler than   (\ref{formlamb6j}) in many aspects: it is a sum of a product of basic hypergeometric of 3-2 type rather than 4-3 type, and moreover it explicitely shows that $\Lambda^{BC}_{A}$ are Laurent polynomials in $q^{i\rho_X}.$ 
 
\medskip 
  
Remark: From the  system of constraints satisfied by the coefficients $\Lambda^{BC}_{AD}$, it is easy to show that $(\Lambda^{BC}_{AD})^2$ is a Laurent polynomials in both variables $q^{2X_0+1}, q^{2X_1+1}.$ 

\medskip 
  
The constraint systems described before, as well as explicit expressions, allow us to derive most of the asymptotic properties of these coefficients when some of the variables are large.  
  
\begin{proposition}  
The coefficients $\Lambda^{BC}_{AD}$ satisfy the following inequality:  
\begin{equation}  
\forall A,D \in \halfinteger, \forall R \in {\onehalf}\ZZ, \forall m  \in {\onehalf}\ZZ, \exists \;{\cal C}, {\cal B}\;>0,  
\forall\rho\in \RR, \forall B > {\cal B}, \vert \Lambda^{B \;B+R}_{A\;\;D}(m,\rho) \vert \leq  
{\cal C} Bq^{2B}\label{asymptlambda1}.  
\end{equation}  
\end{proposition}  
  
\Proof  
 
The proof is divided in two steps. We first show this inequality when $A=D$, by a direct computation. Then we use the system  (\ref{systan}) and an induction argument to show this inequality when $A\not=D.$  
We can always assume that $R\geq 0$, because the other case is deduced from this one and the identity (\ref{symetryLambda}). 
Using the relation (\ref{lambda3j}), we have $\vert \Lambda^{BB+R}_{A}(m,\rho)\vert \leq \sum_{j}u(j,m-j;B), $ where $u(j,k;B)=\vert \Clebpsi{B+R}{B}{A}{\;\;\;j}{k}{m}\;  
\Clebpsi{B+R}{B}{A}{\;\;-j}{\!-k}{-m}\vert.$

Using the formula (4) of (\cite{KV}, section 14.3.5), we easily obtain that:  
\begin{equation}  
\hskip -0.3cm\frac{\vert u(j,m\!-\!j;B)\vert }{q^{(A-R)(2B+R+A+1)}}=  
\frac{[2A]![2A+1]![A\!+\!B\!-\!j]![B\!+\!A\!+\!j]![2B\!+\!R\!-\!A]!F(j,m\!-\!j;B)F(-j,j\!-\!m;B)}{[A\!+\!R]! [A\!-\!R]![A\!+\!m]! [A\!-\!m]! [B\!+\!m\!-\!j]![B\!-\!m\!+\!j]![2B\!\!+\!\!R\!\!+\!\!A\!\!+\!\!1]!} \nonumber  
\end{equation}  
where   
$F(j,m\!-\!j;B)={}_{3}\!\Phi_2\left[\!\!\begin{array}{ccc} m\!-\!A &-2B\!-\!R\!-\!A\!-\!1 &\!\! R\!-\!A\\&\!\!\!\!\!\!\!\!\!\!\!\!\!\!\!\!\!\!\!\!\!\!\!\!\!\!\!\!\!\!-B\!-\!A\!+\!j &\!\!\!\!\!\!\!\!\!\!\!\!\!\!\!\!\!\!\!\!\!\!\!\!\!\!\!\!\!\!\!\!\!\!-2A  \end{array} \!\!\right]\!.$  
Using the fact that $ -2B\!-\!A\!-\!R\leq -B\!-\!A\!+\!j  
\leq R\!-\!A,$ we can choose $C_1$ such that $\forall j, \forall B, \vert F(j,m-j;B)\vert \leq C_1 q^{-2Bn(m)}$ where $n(m)=inf(A-m,A-R)$.  
As a result, there exists a constant $C_1$ such that   
\begin{eqnarray}  
\forall j,B,&& u(j,m-j;B)\leq C_1 q^{2B(A-R-n(m)-n(-m))-\vert B+j\vert (A+m)-\vert B-j\vert (A-m)-2B(2A+1)}.\nonumber 
\end{eqnarray}  
From this, and the various inequalities between $B, j, m, R$,  it is easy to show that we can find $C_2$ such that 
\begin{equation} 
\vert \Lambda^{BB+R}_{A}(m,\rho)\vert \leq C_2(2B+2R+1)q^{2B}.\label{bounddiagonalLambda} 
\end{equation} 
We can now prove this bound for the non diagonal case using an induction argument on $\vert A-D\vert.$ Indeed from (\ref{systan}), 
\begin{equation} 
\hskip -0.2cm\vert  \coefblambda{B}{B\!+\!R}{A}{D\!+\!1} 
\coefblambda{\onehalf}{D\!+\!\onehalf}{D\!+\!1}{D}\vert \! \leq \! 
\sum_{KU}\!\frac{(\qn{d_{U}}\qn{d_{D+1}})^{\frac{1}{2}}\vert\coefblambda{K}{U}{A}{D}\vert} 
 {(\qn{d_{B+R}}\qn{d_{D+\onehalf}})^{\frac{1}{2}}}\vert 
\sixjn{A}{D\!+\!\onehalf}{U}{\onehalf}{\!B\!+\!R}{D\!+\!1\!}{0}\!\!\!\!\sixjn{\onehalf}{U}{\!B\!+\!R}{A}{B}{K}{0}\!\!\!\!\sixjn{D}{U}{K}{B}{\onehalf}{D\!+\!\onehalf}{0}\!\!\!\vert \nonumber. 
\end{equation} 
Note that, due to the $Y$ functions, the previous sums contain at most four terms: $U=B+R\pm\onehalf, K=B\pm \onehalf.$ 
From the asymptotics property (\ref{bounddiagonalLambda}) when $B$ goes to $+\infty,$ 
and the induction hypothesis, we obtain that there exists a constant $C_3,$ with  
 $\forall B, \forall \rho,\vert  \coefblambda{B}{B+R}{A}{D+1}(m,\rho) 
\coefblambda{\onehalf}{D+\onehalf}{D+1}{D}(m,\rho)\vert \leq  
C_3 B q^{2B}.$ 
From the explicit formulae of the coefficients  
$\coefblambda{\onehalf}{D+\onehalf}{D+1}{D}(m,\rho),$ in \cite{BR1}, there exists $C_4>0$ such that  
$\forall\rho\in \RR , \forall B,\vert \coefblambda{\onehalf}{D+\onehalf}{D+1}{D}(m,\rho),\vert >C_4,$ 
as a result we get that the induction hypothesis is also true for $A, D+1.$ 
\eoProof  
  
We will also need the following asymptotics:  
\begin{proposition}  
The coefficients $\Lambda^{BC}_{A}$ own the following asymptotic behaviour: 
\begin{eqnarray}  
&&\hskip -2cm\forall B \in \halfinteger,\forall R \in {\onehalf}\ZZ, B \pm R \in \NN, \forall (X_0X_1) \in {\mathbb S}, \exists {\cal C} \in \RR \nonumber\\  
&&\hskip -2cm\Lambda^{B \;\;K+R}_{\;\;K} \sim_{K\rightarrow +\infty} q^{2K \vert R \vert} \; {\cal C} \; Y^{(1)}_{(B,R)} \;\;\mbox{with}\;\;(\; {\cal C}=1\;\;\mbox{when}\;\; R=0).\label{asymptlambda2}  
\end{eqnarray}  
\end{proposition}  
  
\Proof  
 
We use the formula (\ref{lambda3j}) to express $\Lambda^{B K+R}_{K}$. The asymptotic when $K$ is large is easily obtained  using the formula (1) of  (\cite{KV} section 14.3.5). The details are left to the reader and we easily obtain the statement of the proposition. 
\eoProof  
  
\medskip 

In the section $5$ we will need the following identity on the following series  of weighted diagonal coefficients:   
\begin{proposition}  
If $P>0$ the following series is absolutely and uniformly convergent in $\rho_X$ and we have: 
\begin{eqnarray}  
&&\hskip -1.3cm\sum_{KM}\qd{K} \Clebpsi{K}{P}{M}{\sigma}{-P}{\sigma-P}\Clebphi{P}{K}{M}{-P}{\sigma}{\sigma-P}\Lambda^{KM}_{P}(X_0X_1) = \frac{\hbar\qd{P}q^{-2i \sigma \rho_X}e^{i\pi \sigma}\delta_{m_X,P}}{2\pi{\cal P}{(X_0X_1)}} .\label{formsumlambda}  
\end{eqnarray}  
When $P=0$ we still have : 
\begin{equation} 
\sum_{K}[d_K]\Lambda^{KK}_0(X_0X_1)= \frac{\hbar\qd{P}\delta_{m_X,P}}{2\pi{\cal P}{(X_0X_1)}}. 
\end{equation} 
\end{proposition}  
  
\Proof \\  
When $P=0$, we see from (\ref{lambda3j}) that $\Lambda^{KK}_0(X_0X_1)$ is non zero only when $m_X=0.$ 
The series $\sum_{K}[d_K]\Lambda^{KK}_0(X_0X_1)$ is in this case very simple to compute because $\Lambda^{KK}_0(X_0X_1)=\frac{[(2X_0+1)(2K+1)]}{[2X_0+1][2K+1]},$ and we obtain the announced statement. Note that in this case the series is not uniformly convergent and is even divergent at $\rho_Z=0.$ 
When $P>0$ we will  first show that the series is absolutely and uniformly convergent.  
Using  the analogue of the Van-Der-Waerden formula for $3j$ \cite{KV}, as well as relation (\ref{symmetry3j}) we can choose $C_1$ such that: 
\begin{eqnarray}  
\vert \Clebpsi{K}{P}{K\!+\!R}{\sigma}{\!-P}{\sigma\!-\!P}\Clebphi{P}{K}{K\!+\!R}{\!\!-P}{\sigma}{\sigma\!-\!P} \vert \!\!\!\! &=&\!\!\!\! \vert \frac{[d_{K+R}][2P]![2K\!\!+\!\!R\!\!-\!\!P]![K\!\!+\!\!R\!\!+\!\!P\!\!-\!\!\sigma]![K\!\!+\!\!\sigma]!  
}{[P\!\!-\!\!R]![R\!\!+\!\!P]![2K\!\!+\!\!R\!\!+\!\!P\!+\!\!1]![K\!\!+\!\!R\!\!-\!\!P\!\!+\!\!\sigma]![K\!\!-\!\!\sigma]!}\vert\nonumber\\  
&\leq& C_1 \; q^{2KP}\nonumber 
\end{eqnarray} 
Using $\Lambda^{K K+R}_{P}(X_0X_1) \leq C_2 \;K q^{2K}$, we obtain that the series is absolutely and uniformly convergent. 
Using  (\ref{systan})(\ref{intlambda})(\ref{contra})(\ref{lambda3j}) we have  
\begin{eqnarray}  
&&\hskip -1cm\int d{\cal P}(X_0X_1)\;\left(\sum_{KM}\qd{K} \Clebpsi{K}{\;P}{\;\;M}{\sigma}{-P}{\sigma-P}\Clebphi{\;P}{K}{\;\;M}{-P}{\sigma}{\sigma-P}\Lambda^{KM}_{P}(X_0X_1)\right) \Lambda^{AC}_{P}(X_0X_1)=\nonumber\\  
&&\hskip -1cm=\sum_{KM}\qd{K} \Clebpsi{K}{\;P}{\;\;M}{\sigma}{-P}{\sigma-P}\Clebphi{\;P}{K}{\;\;M}{-P}{\sigma}{\sigma-P} \int d{\cal P}(X_0X_1) \Lambda^{KM}_{P}(X_0X_1) \Lambda^{AC}_{P}(X_0X_1)=\nonumber 
\end{eqnarray}  
\vskip -0.4cm  
\begin{eqnarray}  
&&\hskip -1cm=\sum_{KM}\qd{K} \Clebpsi{K}{\;P}{\;\;M}{\sigma}{-P}{\sigma-P}\Clebphi{\;P}{K}{\;\;M}{-P}{\sigma}{\sigma-P}   
\sum_{JU} \frac{(\qn{2U\!+\!1}\qn{2P\!+\!1})^{\frac{1}{2}}}{(\qn{2C\!+\!1}\qn{2M\!+\!1})^{\frac{1}{2}}e^{i \pi (C+M-P-U)}}\times\nonumber\\  
&&\times\sixjn{A}{M}{U}{K}{C}{P}{0}\!\!\!\sixjn{K}{U}{C}{P}{A}{J}{0}\!\!\!\sixjn{P}{U}{J}{A}{K}{M}{0}.  
\int d{\cal P}(X_0X_1) \Lambda^{JU}_{P}(X_0X_1) =\nonumber 
\end{eqnarray}  
\vskip -0.4cm  
\begin{eqnarray}  
&&\hskip -1cm=\sum_{KM}\qd{K} \Clebpsi{K}{\;P}{\;\;M}{\sigma}{-P}{\sigma-P}\Clebphi{\;P}{K}{\;\;M}{-P}{\sigma}{\sigma-P}  
e^{2i\pi A} \delta_{A,K}\frac{\qd{P}}{\qd{A}} \sixjn{A}{M}{P}{A}{C}{P}{0}=\nonumber\\  
&&\hskip -1cm=\sum_{M} \qd{M} \Clebphi{\;\;M}{A}{\;P}{P-\sigma}{\sigma}{P}\Clebpsi{A}{\;\;M}{\;P}{\sigma}{P-\sigma}{P}  
 \sixjn{A}{M}{P}{A}{C}{P}{0}=\nonumber 
\end{eqnarray}  
\vskip -0.4cm  
\begin{eqnarray}  
&&\hskip -1cm=\sum_{M}\qd{M} \sixjn{A}{M}{P}{A}{C}{P}{0} \frac{\hbar}{2\pi}\int^{\frac{\pi}{\hbar}}_{-\frac{\pi}{\hbar}} d\rho \; q^{2i \sigma \rho} \Lambda^{AM}_{P}(P,\rho).\nonumber  
\end{eqnarray}  
Then we use  the relation  
\begin{eqnarray}  
&&\hskip -1cm \sum_{M}\qd{M}\sixjn{A}{M}{P}{A}{C}{P}{0}\!\!\!\!\Lambda^{AM}_{P}(X_0X_1)=e^{-2i\pi A}\qd{P}\;\Lambda^{AC}_{P}(X_1X_0).\nonumber  
\end{eqnarray}  
which has been proved in \cite{BR1} and which is the key point to prove  unitarity of the representation and  
the relation   
\begin{eqnarray}  
\hskip -1cm\int^{\frac{\pi}{\hbar}}_{-\frac{\pi}{\hbar}}  d\rho \; q^{2i \sigma \rho} \Lambda^{AC}_{P}(-P,\rho)=  
\int^{\frac{\pi}{\hbar}}_{-\frac{\pi}{\hbar}}  d\rho \; q^{-2i \sigma \rho} \Lambda^{AC}_{P}(P,\rho)  
\end{eqnarray}  
which is proved using twice (\ref{lambda3j}),(\ref{contra}) and (\ref{intlambda}). 
Combining these two steps, we obtain  
\begin{eqnarray}  
&&\hskip -1cm\int d{\cal P}(X_0X_1)\;\left(\sum_{KM}\qd{K} \Clebpsi{K}{\;P}{\;\;M}{\sigma}{-P}{\sigma-P}\Clebphi{\;P}{K}{\;\;M}{-P}{\sigma}{\sigma-P}\Lambda^{KM}_{P}(X_0X_1)\right) \Lambda^{AC}_{P}(X_0X_1)=\nonumber\\  
&&=\int d{\cal P}{(X_0X_1)} \;\frac{\hbar\;e^{2i \pi \sigma} \qd{P}q^{-2 \sigma (X_0+X_1+1)}\delta_{X_0-X_1,P}}{2\pi\; {\cal P}{(X_0X_1)}} \Lambda^{AC}_{P}(X_0X_1).\nonumber  
\end{eqnarray}  
Finally using Lemma 1 we conclude this  proof.  
\eoProof  
  
\bibliographystyle{unsrt}

\end{document}